%% file: main.tex
\documentclass[english,ruled]{article}
\usepackage[T1]{fontenc}
\usepackage[latin9]{inputenc}
\usepackage{verbatim}
\usepackage{subcaption}
\usepackage[linesnumbered,ruled]{algorithm2e}
\usepackage{amsmath}
\usepackage{amsthm}
\usepackage{etoolbox}
\usepackage{amssymb}
\usepackage{graphicx}
\usepackage{xcolor}
\usepackage{mathtools}
\usepackage{multicol}
\usepackage{multirow}
\usepackage{geometry}
\usepackage{booktabs}
\usepackage{enumitem}
\usepackage{setspace}
\makeatletter
\usepackage[toc,page,header]{appendix}
\usepackage{minitoc}
\usepackage{ifthen}
\newboolean{doublecolumn}
\setboolean{doublecolumn}{false}
\newboolean{arxiv}
\setboolean{arxiv}{true}
\usepackage{pgfplots}
\usepackage{pdflscape}
\usepackage[pagebackref=true]{hyperref}
\usepackage{cleveref}
\usepackage{authblk}
\usepackage{todonotes}
\pgfplotsset{compat=1.15}

\newboolean{proofinpaper}
\setboolean{proofinpaper}{false}

\newcommand{\inputinpaper}[1]{%
  \ifbool{proofinpaper}{\begin{proof}
						\input{#1}{}%
						\end{proof}}}

\newcommand{\inputinappendix}[1]{%
\ifbool{proofinpaper}{}{\input{#1}}%
}

\makeatletter
\renewcommand\paragraph{%
  \@startsection{paragraph}{4}{\z@}%
    {-1.5ex\@plus -0.5ex \@minus -0.2ex}
    {-.5em \@plus -.1em}%
    {\normalfont\normalsize\bfseries}%
}
\makeatother
\newcommand{\assign}{\coloneqq}
\newcommand{\nin}{\not\in}

\newcommand{\tmop}[1]{\ensuremath{\operatorname{#1}}}
\newcommand{\tmtextbf}[1]{\text{{\bfseries{#1}}}}
\newcommand{\tmtextit}[1]{\text{{\itshape{#1}}}}
\newcommand{\tmem}[1]{\textit{#1}}
\newcommand{\mathd}{\mathrm{d}}
\definecolor{deepskyblue}{rgb}{0.0, 0.75, 1.0}

\input{macros.tex}

\renewcommand*\backref[1]{\ifx#1\relax \else (cited on #1) \fi}
\theoremstyle{plain}
\newtheorem{lem}{\protect\lemmaname}[section]
\theoremstyle{remark}

\theoremstyle{plain}
\newtheorem{thm}{\protect\theoremname}[section]
\theoremstyle{plain}
\newtheorem{prop}{\protect\propositionname}[section]
\providecommand{\corollaryname}{Corollary}
\theoremstyle{plain}

\theoremstyle{plain}
\newtheorem{exple}{\protect\examplename}[section]
\theoremstyle{plain}
\newtheorem{definition}{\protect\definitionname}[section]

\providecommand{\lemmaname}{Lemma}
\providecommand{\remarkname}{Remark}
\providecommand{\theoremname}{Theorem}
\providecommand{\examplename}{Example}
\providecommand{\propositionname}{Proposition}
\providecommand{\definitionname}{Definition}

\newcommand{\ogd}{{\texttt{OGD}}}
\newcommand{\hdm}{{\texttt{HDM}}}
\newcommand{\adagrad}{{\texttt{AdaGrad}}}
\newcommand{\osgm}{{\texttt{OSGM}}}
\newcommand{\osgmrx}{\texttt{OSGM-R}}
\newcommand{\osgmhandsonrx}{\texttt{Lookahead OSGM-R}}
\newcommand{\osgmhandsoffrx}{\texttt{Vanilla OSGM-R}}

\newcommand{\osgmhx}{\texttt{OSGM-H}}
\newcommand{\osgmhandsonhx}{\texttt{Monotone Lookahead OSGM-H}}
\newcommand{\osgmhandsoffhx}{\texttt{Monotone OSGM-H}}

\newcommand{\overleq}[1]{\mathrel{\overset{#1}{\leq}}}

\crefdefaultlabelformat{#2\textbf{#1}#3} 
\crefname{section}{\textbf{section}}{\textbf{sections}}
\Crefname{section}{\textbf{Section}}{\textbf{Sections}}
\crefname{thm}{\textbf{Theorem}}{\textbf{theorems}}
\Crefname{thm}{\textbf{Theorem}}{\textbf{Theorems}}
\crefname{lem}{\textbf{Lemma}}{\textbf{lemmas}}
\Crefname{lem}{\textbf{Lemma}}{\textbf{Lemmas}}
\crefname{prop}{\textbf{proposition}}{\textbf{propositions}}
\Crefname{prop}{\textbf{Proposition}}{\textbf{Propositions}}
\crefname{algorithm}{\textbf{algorithm}}{\textbf{algorithms}}
\Crefname{algorithm}{\textbf{Algorithm}}{\textbf{Algorithms}}
\crefname{coro}{\textbf{Corollary}}{\textbf{corollaries}}
\Crefname{coro}{\textbf{Corollary}}{\textbf{corollaries}}
\crefname{definition}{\textbf{Definition}}{\textbf{definitions}}
\Crefname{definition}{\textbf{Definition}}{\textbf{definitions}}
\crefname{table}{\textbf{Table}}{\textbf{tables}}
\Crefname{table}{\textbf{Table}}{\textbf{tables}}
\crefname{figure}{\textbf{Figure}}{\textbf{figures}}
\Crefname{figure}{\textbf{Figure}}{\textbf{figures}}
\crefname{exple}{\textbf{Example}}{\textbf{examples}}
\Crefname{exple}{\textbf{Example}}{\textbf{examples}}

\begin{document}

\title{Gradient Methods with Online Scaling \\ \vspace{10pt}
\Large Part I. Theoretical Foundations \\
}

\author[1]{Wenzhi Gao\thanks{gwz@stanford.edu, equal contribution}}
\author[2]{Ya-Chi Chu\thanks{ycchu97@stanford.edu, equal contribution}}
\author[1,3]{Yinyu Ye\thanks{yyye@stanford.edu}}
\author[1,3]{Madeleine Udell\thanks{udell@stanford.edu}}
\affil[1]{ICME, Stanford University}
\affil[2]{Department of Mathematics, Stanford University}
\affil[3]{Department of Management Science and Engineering, Stanford University}

\maketitle

\begin{abstract}
This paper establishes the theoretical foundations of the online scaled gradient methods ({\osgm})\footnote{This paper extends two previous works \cite{gao2024gradient} and \cite{chu2025provable}.}, a framework that utilizes online learning to adapt stepsizes and provably accelerate first-order methods.  
{\osgm} quantifies the effectiveness of a stepsize by a feedback function motivated from a convergence measure and uses the feedback to adjust the stepsize through an online learning algorithm.
Consequently, instantiations of {\osgm} achieve convergence rates that are asymptotically no worse than the optimal stepsize. 
{\osgm} yields desirable convergence guarantees on smooth convex problems, including 1) trajectory-dependent global convergence on smooth convex objectives; 2) an improved complexity result on smooth strongly convex problems, and 3) local superlinear convergence. Notably, {\osgm} constitutes a new family of first-order methods with non-asymptotic superlinear convergence, joining the celebrated quasi-Newton methods. Finally, {\osgm} explains the empirical success of the popular hypergradient-descent heuristic in optimization for machine learning.
\end{abstract}

\input{sec_intro.tex}

\input{sec_framework.tex}
\input{sec_landscape.tex}
\input{sec_action.tex}
\input{sec_oco.tex}
\input{sec_algodemo.tex}
\input{sec_practical.tex}
\input{sec_exp.tex}
\input{sec_literature.tex}

\input{conclusion.tex}

\renewcommand \thepart{}
\renewcommand \partname{}

\bibliography{ref.bib}
\bibliographystyle{plain}

\doparttoc
\faketableofcontents
\part{}

\newpage
\appendix

\addcontentsline{toc}{section}{Appendix}
\part{Appendix}
\parttoc

\newpage
\input{app_landscape.tex}

\input{app_actions.tex}
\input{app_oco.tex}
\input{app_algodemo.tex}
\input{app_moreexps.tex}

\end{document}

%% file: macros.tex
\global\long\def\vertiii#1{\left\vert \kern-0.25ex  \left\vert \kern-0.25ex  \left\vert #1\right\vert \kern-0.25ex  \right\vert \kern-0.25ex  \right\vert }%

\global\long\def\diam{\mathrm{diam}}%

\global\long\def\argmin{\operatornamewithlimits{arg\,min}}%

\global\long\def\diag{\mathrm{Diag}}%

\global\long\def\and{\mathrm{and}}%

\global\long\def\Rbb{\mathbb{R}}%

\global\long\def\Acal{\mathcal{A}}%

\global\long\def\Mcal{\mathcal{M}}%

\global\long\def\Ocal{\mathcal{O}}%

\global\long\def\Pcal{\mathcal{P}}%

\global\long\def\Xcal{\mathcal{X}}%

%% file: sec_intro.tex
\section{Introduction} \label{sec:intro}
Consider the $L$-smooth and $\mu$-strongly convex optimization problem
$\min_{x\in \mathbb R^n}f(x)$.
When $f$ has a large condition number $\kappa \coloneqq L/\mu$, vanilla gradient descent with constant scalar stepsize $1/L$ converges slowly.
Instead of using a constant scalar stepsize, preconditioned gradient descent chooses a preconditioner $P_k \in \Rbb^{n\times n}$, a matrix stepsize, to scale the gradient and accelerate convergence at iteration $k$:
\begin{equation} \label{eqn:pgd}
   x^{k + 1} = x^k - P_k \nabla f (x^k).
\end{equation}
We can locally measure the quality of stepsize $P_k$ by the contraction ratio of suboptimality at $x^k$:
\[ r_{x^k} (P_k):= \tfrac{f(x^{k+1}) - f^{\star}}{f(x^k) - f^{\star}} =  \tfrac{f (x^k - P_k \nabla f (x^k)) - f^{\star}}{f (x^k) - f^{\star}}. \]
The suboptimality after $K$ iterations of preconditioned gradient descent \eqref{eqn:pgd} is the product of these ratios $\{r_{x^k} (P_k)\}$ and can be further bounded by their cumulative sum using the arithmetic-geometric mean inequality:
\begin{equation} \label{eqn:intro-reduction}
\textstyle \tfrac{f (x^{K + 1}) - f^{\star}}{f (x^1) - f^{\star}} = \prod_{k =
   1}^K \tfrac{f (x^{k + 1}) - f^{\star}}{f (x^k) - f^{\star}} = \prod_{k = 1}^K r_{x^k} (P_k) \leq ( \tfrac{1}{K} \sum_{k = 1}^K
   r_{x^k} (P_k) )^K.
\end{equation}
Given this analysis, how should we choose the stepsize $P_k$? 
Fast convergence is achieved if we prespecify stepsizes $\{P_k\}$ to minimize $\sum_{k = 1}^K r_{x^k} (P_k)$.
However, solving this problem is hard for a given instance $f$, as the function $r_{x^k} (P)$ depends on the previous sequence of stepsizes and the optimization landscape. However, this problem is ideally suited to online learning, which optimizes a cumulative sum of functions with provable regret guarantees even when the functions are adversarially chosen. 
For example, updating $\{P_k\}$ with online gradient descent $P_{k + 1} = P_k - \eta \nabla r_{x^k} (P_k)$ guarantees sublinear regret with respect to any fixed stepsize $\hat{P}$ \cite{orabona2019modern}:
\begin{equation} \label{eqn:intro-regret}
 \textstyle \tfrac{1}{K} \sum_{k = 1}^K r_{x^k} (P_k) \leq \tfrac{1}{K} \sum_{k = 1}^K r_{x^k} (\hat{P}) +\mathcal{O} ( \tfrac{1}{\sqrt{K}}).
\end{equation} 
Given the freedom to choose $\hat{P}$, we may take $\hat{P}$ equal to $P^{\star}_r$ that achieves the optimal condition number $\kappa^{\star} < \kappa$ and hence the ratio $r_{x^k} (P^{\star}_r) \leq 1 - \tfrac{1}{\kappa^{\star}}$. The regret guarantee \eqref{eqn:intro-regret} together with \eqref{eqn:intro-reduction} implies a convergence rate
\begin{equation} \label{eqn:intro-optkappa}
   \textstyle f (x^{K + 1}) - f^{\star} \leq [f (x^1) - f^{\star}] ( 1 - \tfrac{1}{\kappa^{\star}}+\mathcal{O} ( \tfrac{1}{\sqrt{K}} ) )^K,
\end{equation}
as well as an asymptotic iteration complexity $\Ocal(\kappa^\star \log(1/\varepsilon))$. Moreover, if $f$ is (locally) quadratic with Hessian matrix $H \in \mathbb{S}^n_{++}$, then $\hat{P} = H^{-1}$ yields perfect conditioning $\kappa^\star = 1$ and hence superlinear convergence:
\[ \textstyle f (x^{K + 1}) - f^{\star} \leq [f (x^1) - f^{\star}]  ( \Ocal(\tfrac{1}{\sqrt{K}} )^K).
\]
This concise proof establishes a new problem-dependent acceleration mechanism for gradient methods, along with a notable local superlinear convergence guarantee. 
More generally, our \textit{online scaled gradient methods} (\osgm) is a family of first-order methods that updates the stepsize on the fly using online learning on a convergence measure.

\subsection{Contributions} 
This paper establishes the theoretical foundations of {\osgm} and showcases several instantiations of {\osgm}. Notable features of {\osgm} include:

\paragraph{Non-asymptotic complexity improvement on gradient descent.}  
One instance of {\osgm} (\Cref{sec:handson-rx}) is guaranteed to find an $\varepsilon$-optimal solution in
\begin{equation} \label{eqn:complexity}
   K_\varepsilon = \Big\lceil \min \big\{ \kappa^\star L^2 \| P_1 - P_r^\star \|^2_F + \kappa^\star \log \big( \tfrac{f(x^1) - f^{\star}}{\varepsilon} \big), ~\kappa \log(\tfrac{f(x^1) - f^{\star}}{\varepsilon}) \big\}\Big\rceil
\end{equation}
iterations, where $P_r^\star$ is the stepsize that achieves the optimal condition number $\kappa^\star$. The complexity in \eqref{eqn:complexity} is no worse than $\Ocal{(\kappa \log(1/\varepsilon))}$, the complexity of vanilla gradient descent. When $\varepsilon \rightarrow 0$, \eqref{eqn:complexity} becomes competitive with $\Ocal{(\kappa^\star \log(1/\varepsilon))}$, the complexity achieved by $P_r^{\star}$. In general, {\osgm} guarantees a problem-dependent complexity $\Ocal(\kappa^{\star}\log(1/\varepsilon))$ on smooth strongly convex problems, improving on classical accelerated methods when $\kappa^{\star} < \sqrt{\kappa}$ and $\varepsilon \rightarrow 0$. {\osgm} can be viewed as an acceleration mechanism of gradient descent that uses preconditioning to achieve acceleration. It differs from the existing momentum-based schemes \cite{nesterov1983method,nesterov2013introductory,necoara2019linear}.

\paragraph{Non-asymptotic local superlinear convergence.}  
To our knowledge, {\osgm} represents the second family of first-order methods that achieve local non-asymptotic superlinear convergence after the celebrated quasi-Newton methods \cite{nocedal1999numerical,fletcher2000practical}. 
The non-asymptotic superlinear convergence rate of {\osgm} matches or surpasses that of quasi-Newton methods, which were analyzed recently \cite{rodomanov2021greedy,rodomanov2021new,rodomanov2022rates,jin2023non,jin2024non,jiang2023online,jiang2024online}. Our superlinear convergence analysis is simple and can be of independent interest.

\paragraph{Analysis of hypergradient descent heuristic.}
An instantiation of {\osgm} (\Cref{sec:monotone-hx}) simplifies to hypergradient descent \cite{schraudolph1999local, gunes2018online}, a popular heuristic that has received limited analysis \cite{kunstner2024searching}. {\osgm} provides a rigorous analysis of hypergradient descent \cite{schraudolph1999local, gunes2018online} through the lens of both online learning (\Cref{thm:global-conv-handson-hx}) and provides new analysis tool using potential function (\Cref{thm:potential-handson-rx}, \Cref{thm:potential-handson-hx}).
Our analysis explains the empirical success of hypergradient descent and motivates its practical extensions. \\

{\osgm} also complements the fruitful line of research in adaptive gradient methods {\cite{kingma2014adam,duchi2011adaptive,hinton2012neural}}, parameter-free first-order methods {\cite{lan2023optimal,li2023simple}}, and parameter-free online learning algorithms {\cite{orabona2019modern,orabona2016coin}}.   In particular, {\osgm} offers a new perspective to model stepsize selection as an online learning problem in which the sequence of decisions lies in the \textit{stepsize space}.
This perspective contrasts with the existing literature on adaptive gradient methods: online learning applied to the \emph{primal space} and the stepsizes are chosen by a fixed rule to minimize a regret upper bound \cite{duchi2011adaptive,mcmahan2010adaptive,hazan2007adaptive}.
This shift in the perspective delivers sharp guarantees for smooth convex objectives. To our knowledge, the only prior work in a similar spirit is \cite{zhuang2019surrogate}, which adapts the stepsize to stochastic noise in nonconvex optimization. \\

\Cref{sec:osgm} introduces the {\osgm} framework along with a roadmap for the rest of the paper. \Cref{sec:literature} surveys the related works in detail.

\paragraph{Not in this paper.} {\osgm} also demonstrates excellent practical convergence. However, discussion of the setup most relevant to practical problems does not fit into this paper. We refer the interested readers to a preliminary conference version of our paper \cite{chu2025provable} and Part II of this paper.

\subsection{Notations}
We use $\| \cdot \|$ to denote vector Euclidean norm or matrix operator norm and $\langle \cdot, \cdot
\rangle$ to denote Euclidean or Frobenius inner product.
Letters $A, a$ denote matrices and scalars. 
$\| A \|_F \assign \sqrt{\sum_{i j} a_{i j}^2}$ denotes the matrix Frobenius norm. Given a positive definite matrix $A \in \mathbb{S}^n_{++}$, we define $\|\cdot\|_A \coloneqq \sqrt{\langle x, A x \rangle}$. Given two symmetric matrices $A, B$, $A \succeq B$ if $A - B$ is positive semidefinite. 
$\Pi_{\mathcal{P}} [\cdot]$ denotes the orthogonal projection onto a closed convex set $\mathcal{P}$. Given a vector $d\in \Rbb^n$, $\text{Diag}(d)$ denotes the diagonal matrix with elements of $d$ on its diagonal.
We use $\mathcal{X}^{\star} = \{ x : f (x) = f^{\star} \}$ to denote the optimal
set of $f$; $\tmop{dist} (P, \mathcal{P}) \assign \| P - \Pi_{\mathcal{P}}
[P] \|_F$ denotes the distance between a point $P$ and a closed convex set
$\mathcal{P}$; $\tmop{diam} (\mathcal{P}) = \max_{X, Y \in \mathcal{P}}  \| X
- Y \|_F$ denotes the diameter of set $\Pcal$ in Frobenius norm. 
A function $f$ is $L$-smooth (has $L$-Lipschitz continuous gradient) if it satisfies $\| \nabla f(x) - \nabla f(y) \| \leq L \|x - y\|$ for all $x, y \in \Rbb^n$; a function $f$ has $H$-Lipschitz Hessian if $\| \nabla^2 f(x) - \nabla^2 f(y) \| \leq H \|x - y\|$ for all $x, y \in \Rbb^n$. We use superscript $x^k$ to index algorithm iterates and subscript $P_k$ to index
stepsize sequence.

%% file: sec_framework.tex
\section{Online scaled gradient methods} \label{sec:osgm}

{\osgm} models stepsize selection as a sequential decision-making problem in a
collaborative environment. 
There are two agents: a stepsize scheduling agent (Scheduler) and a landscape agent (Landscape). Scheduler chooses stepsize from a closed convex candidate set $\Pcal \subseteq \Rbb^{n\times n}$. An iteration of {\osgm} has three steps (\Cref{alg:osgm}):

\begin{enumerate}[leftmargin=40pt, label={{\textsf{Step \arabic*}.}}, itemsep=0pt]
	\item Scheduler makes decision $P_k \in \Pcal$ and proposes an update $x^{k + 1 / 2} = x^k - P_k \nabla f (x^k)$.
	\item Landscape queries the objective function $f$ and 1) chooses the next iterate $x^{k + 1} =\mathcal{M} (x^k, x^{k + 1 / 2})$ based on current $x^k$ and the proposal $x^{k+1/2}$; 2) evaluates the quality of $P_k$ and provides feedback $\ell_{x^k} (P_k)$ to the scheduler.
	\item Scheduler updates the stepsize using an online learning algorithm $P_{k + 1} =\mathcal{A} (P_k, \{ \ell_{x^j} \}_{j \leq k})$.
\end{enumerate}
Feedback $\ell_x$, landscape action $\mathcal{M}$, and online algorithm $\mathcal{A}$ combine to yield different variants of {\osgm}. 
\begin{algorithm}[h]
{\textbf{input} Initial point $x^1$, initial stepsize $P_1$, feedback $\ell_x$, online algorithm $\mathcal{A}$, landscape action $\Mcal$}\\
\For{k =\rm{ 1, 2,...}}{
$x^{k+1/2} = x^k - P_k \nabla f(x^k)$  \hspace{4pt}$\triangleright$ \textsf{Scheduler suggests the next step based on $P_k$}\\
$x^{k+1} = \Mcal(x^k, x^{k+1/2})$ \hspace{4.9pt}\quad $\triangleright$ \textsf{Landscape  evaluates $P_k$, chooses  $x^{k+1}$, and gives feedback $\ell_{x^k}(P_k)$}\\
$P_{k+1} = \mathcal{A}(P_k, \{\ell_j\}_{j \leq k})$ ~\quad$\triangleright$  \textsf{Scheduler learns a better stepsize from feedback} \\ 
}
\caption{Online scaled gradient methods ({\osgm})\label{alg:osgm}}
\end{algorithm} 

\begin{center}
\includegraphics[scale=0.36]{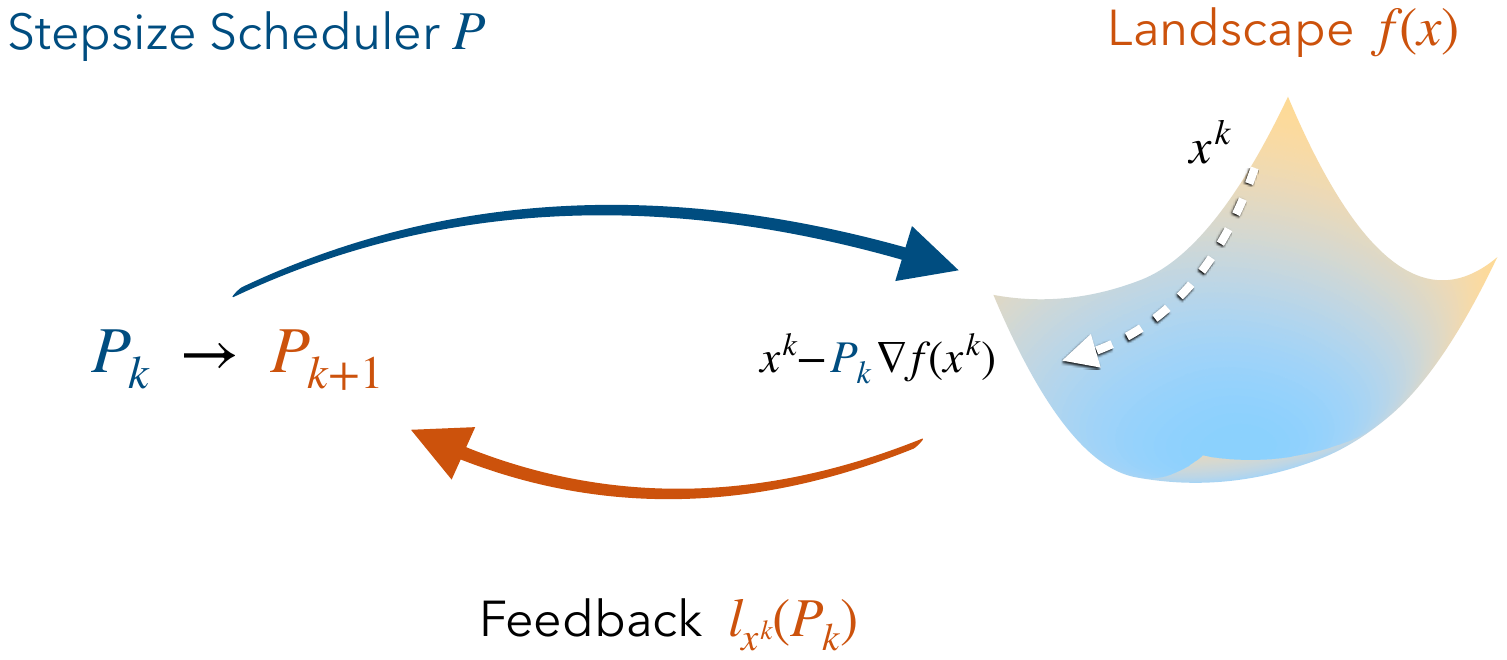}\qquad\qquad
\end{center}

\paragraph{Comparison with online learning.} Our framework builds on online
learning but slightly differs in the spirit: the environment (Landscape) in online learning is often assumed to be adversarial,
while the landscape in {\osgm} shares the same objective of minimizing $f$ and collaborates with the scheduler in two ways:
1) Landscape provides faithful feedback to the scheduler to improve the stepsize. 
2) Landscape can reject bad decisions $P_k$ made by the scheduler and suppress regret in online learning through the landscape action.

\paragraph{Structure of the paper.}
The paper is organized as follows to introduce different components of {\osgm}.
\begin{itemize}[leftmargin=10pt,itemsep=0pt,topsep=5pt]
  \item {\tmem{Feedback $\ell_x (P)$}}. 
A function that locally measures the quality
  of stepsize $P$ at $x$. In our context, {\tmem{smaller}} feedback means
  {\tmem{better}} quality. \Cref{sec:feedback} introduces two important feedback functions.
  
    \item {\tmem{Landscape action $\mathcal{M}$}}. Landscape selectively accepts stepsizes $P_k$ with small feedback and suppresses the
  regret of the learning algorithm $\mathcal{A}$. \Cref{sec:action} discusses different landscape actions and their consequences. 
  
  \item {\tmem{Learning algorithm $\mathcal{A}$}}. An online learning algorithm that
  guarantees performance concerning feedback $\ell_x (P)$. In our
  context, the guarantee refers to the (sublinear) regret with respect to any stepsize $P \in \mathcal{P}$:
  \[ \textstyle \sum_{k = 1}^K \ell_{x^k} (P_k) \leq \sum_{k = 1}^K \ell_{x^k} (P) + o
     (K) . \]
   \Cref{sec:oco} discusses the regret guarantees of online gradient descent on our feedback functions.
\end{itemize}

After establishing these components, we introduce several important instantiations of {\osgm} in \Cref{sec:algodemos} and discuss the practical aspects of {\osgm} in \Cref{sec:practical}.

%% file: sec_landscape.tex
\section{Feedback design and analysis}\label{sec:feedback}

We introduce two feedback functions $\ell_x(P)$: \tmtextit{ratio} feedback and \tmtextit{hypergradient} feedback\footnote{We consider two feedback functions in this paper for simplicity. More feedback functions can be found in \cite{gao2024gradient, chu2025provable}}. 
The ratio feedback requires knowledge of the optimal value $f^{\star}$ and has strong theoretical guarantees. The hypergradient feedback does not require $f^{\star}$ and is more practical. We analyze the properties of each feedback function and associate each feedback function with a minimax optimal stepsize.

\subsection{Ratio feedback}

Given any $x \nin \mathcal{X}^{\star}$, the \tmtextit{ratio feedback} $r_x
(P)$ measures the quality of stepsize $P$ by the contraction ratio of the suboptimality
as the iterate moves from $x$ to $x - P \nabla f (x)$:
\[ r_x (P) \assign \tfrac{f (x - P \nabla f (x)) - f^{\star}}{f (x) - f^{\star}} . \]
The ratio feedback inherits both convexity and smoothness from  $f$. 
Since $r_x (P)$ simply translates and scales the function value at the proposed iterate $u_x (P) \assign f (x - P \nabla f
(x))$, we begin by analyzing $u_x(P)$.

\begin{prop}[Properties of $u_x$]
  \label{prop:property-ux}For any $x \in \Rbb^n$, the function $u_x (P) = f (x - P \nabla f (x))$ has
  gradient
  \[ \nabla u_x (P) = - \nabla f (x - P \nabla f (x)) \nabla f (x)^{\top} . \]
  Moreover, it satisfies the following properties:
  \begin{enumerate}[leftmargin=12pt,itemsep=0pt]
    \item If $f$ is convex and $L$-smooth, then $u_x (P)$ is convex and $L \| \nabla f (x)
    \|^2$-smooth in $P$.
    
    \item If $f$ is $L$-smooth and $\diam(\mathcal P) \leq D$, 
    then $u_x (P)$ is $(L D + 1) \|
    \nabla f (x) \|^2$-Lipschitz in $P$.
  \end{enumerate}
\end{prop}

\input{proofs/proof-property-ux.tex}

With the properties of $u_x$, the analytic properties of $r_x (P)$ follow immediately.

\begin{lem}[Properties of $r_x$]
  \label{lem:property-rx}For any $x \nin \mathcal{X}^{\star}$, the ratio
  feedback $r_x (P)$ has gradient
  \[ \nabla r_x (P) = - \tfrac{\nabla f (x - P \nabla f (x)) \nabla f
     (x)^{\top}}{f (x) - f^{\star}} . \]
  Moreover, it satisfies the following properties:
  \begin{enumerate}[leftmargin=12pt,itemsep=0pt]
    \item If $f$ is convex and $L$-smooth, then $r_x (P)$ is convex, non-negative, and $2
    L^2$-smooth in $P$.
    
    \item If $f$ is $L$-smooth and $\diam (\mathcal P) \leq D$, then $r_x (P)$ is $2L(LD+1)$-Lipschitz in $P$.
  \end{enumerate}
\end{lem}

\input{proofs/proof-property-rx.tex}

\subsection{Hypergradient feedback}

Given any $x \nin \mathcal{X}^{\star}$, the \tmtextit{hypergradient feedback} $h_x (P)$ measures the quality of $P$ by the function value progress relative to the size of the gradient as the iterate moves from $x$ to $x - P \nabla f (x)$:
\[ h_x (P) \assign \tfrac{f (x - P \nabla f (x)) - f (x)}{\| \nabla f (x)
   \|^2}. \]
The hypergradient feedback is motivated by the descent lemma for $L$-smooth
functions:
\[ f ( x - \tfrac{1}{L} \nabla f (x) ) - f (x) \leq - \tfrac{1}{2
   L} \| \nabla f (x) \|^2, \]
which states that the improvement in the function value of a gradient step with
stepsize $\tfrac{1}{L} I$ is proportional to the size of the gradient $\| \nabla f (x)
\|^2$ with ratio $\tfrac{1}{2 L}$. In a similar vein, hypergradient feedback
$h_x (P)$ quantifies the progress in terms of this ratio when stepsize $P$ is used.
The properties of hypergradient feedback $h_x(P)$ also follow from \Cref{prop:property-ux}.

\begin{lem}[Properties of $h_x$]
  \label{lem:property-hx}For any $x \nin \mathcal{X}^{\star}$, the
  hypergradient feedback $h_x (P)$ has gradient
  \[ \nabla h_x (P) = - \tfrac{\nabla f (x - P \nabla f (x)) \nabla f (x)^{\top}}{\| \nabla f (x) \|^2} . \]
  Moreover, it satisfies the following properties:
  \begin{enumerate}[leftmargin=12pt,itemsep=0pt]
    \item If $f$ is convex and $L$-smooth, then $h_x (P)$ is convex and $L$-smooth in $P$.
    
    \item If $f$ is $L$-smooth and $\diam (\mathcal P) \leq D$, then $h_x (P)$ is $(L D +
    1)$-Lipschitz in $P$.
  \end{enumerate}
\end{lem}

\input{proofs/proof-property-hx.tex} 

\subsection{Minimax optimal stepsizes} \label{sec:minimax-stepsize}

Given a feedback function $\ell_x (P)$, we can define its corresponding minimax optimal stepsize $P^{\star}_{\ell}$ by
\begin{align}
  P^{\star}_{\ell} \in & \argmin_{P \in \mathcal{P}} \max_{x \nin \Xcal^\star} ~ \ell_x (P).
  \nonumber
\end{align}
Minimax optimal stepsize $P^{\star}_{\ell}$ represents the best decision the scheduler can make when $x$ is selected adversarially. 
We are interested in the feedback achieved by the minimax optimal stepsize $\ell_x(P^\star_\ell)$ since it represents the performance of the scheduler in the worst case. The minimax optimal feedback $\ell_x(P^\star_\ell)$ can be bounded by quantities studied in the literature. For example, when $\tfrac{1}{L} I \in \mathcal{P}$, the minimax optimal hypergradient feedback $h_x(P^\star_h)$ is no worse than $h_x(\tfrac{1}{L} I) \leq -\frac{1}{2L}$ guaranteed by the descent lemma. For ratio feedback, the performance of $P^\star_r$ can be quantified using the globally optimal preconditioner defined below.

\begin{definition}[Globally optimal preconditioner \cite{kunstner2024searching}] \label{def:opt-preconditioner}
Suppose $f$ is $L$-smooth, $\mu$-strongly convex, and twice-differentiable. Given a candidate set of stepsizes $\mathcal{P}$ such that $\tfrac{1}{L} I \in \mathcal{P}$, the globally optimal preconditioner $P^{\star}_+ \in \mathbb{S}_+^n$ 
  and the optimal condition number $\kappa^\star$ are defined as the optimal solution and the optimal value of
  \[ \min_{\hat{\kappa} \geq 0 , P \in
     \mathcal{P} \cap \mathbb{S}_+^n}  \hat{\kappa} \quad \mathrm{subject~to}
     \quad \tfrac{1}{\hat{\kappa}} P^{- 1} \preceq \nabla^2 f (x) \preceq P^{-
     1} \quad \mathrm{for~ all~} x. \]
\end{definition}

The globally optimal condition number satisfies $\kappa^{\star} \leq \kappa = L / \mu$ since $\frac{1}{L}I \in \Pcal$. 
In practice, recent works {\cite{gao2023scalable,qu2024optimal,kunstner2024searching}} observe that $\kappa^{\star} \ll \kappa$ is often the case. 
Gradient descent with stepsize $P^{\star}_+$ guarantees a contraction ratio of $1-\tfrac{1}{\kappa^\star}$ (see \Cref{sec:proof-opt-preconditioner} for the proof):
\begin{equation} \label{eq:opt-preconditioner}
  f (x - P_+^{\star} \nabla f (x)) - f^{\star} \leq ( 1 - \tfrac{1}{\kappa^{\star}} ) [f (x) - f^{\star}].
\end{equation}
Since $P^{\star}_+ \in \mathcal{P}$, the ratio feedback of the globally optimal preconditioner $r_x (P^{\star}_+)$ upper bounds that of the minimax optimal stepsize $r_x (P^{\star}_r)$. The bound is tight on strongly convex quadratics \cite{gao2024gradient}.

\input{proofs/proof-opt-preconditioner.tex}

\begin{prop}[Feedback of the minimax stepsize]
  \label{prop:minimax-feedback}For any $x \nin \mathcal{X}^{\star}$, the minimax optimal stepsizes satisfy:
  \begin{itemize}[leftmargin=15pt]
    \item If $f$ is $L$-smooth and $\mu$-strongly convex, then $r_x
    (P^{\star}_r) \leq 1 - \tfrac{1}{\kappa^{\star}}$.

    \item If $f$ is $L$-smooth and $\tfrac{1}{L} I \in \mathcal{P}$, then $h_x (P^{\star}_h) \leq -
    \tfrac{1}{2 L}$.
  \end{itemize}
  In addition, if $f$ is a strongly convex quadratic, then the minimax optimal stepsizes are $P^{\star}_r = P^{\star}_+$ and $P^{\star}_h = \tfrac{1}{L} I$.
\end{prop}

\input{proofs/proof-minimax-feedback.tex}

\paragraph{Limitations of minimax optimal stepsize.} Minimax stepsizes are useful for establishing global convergence results. However, they are not affine invariant and may not reflect the properties of the local landscape. For example, the globally optimal condition number $\kappa^{\star}$ depends on the choice of coordinate system when $\mathcal{P} \neq \mathbb{R}^{n \times n}$.

\begin{exple}[Limitations of globally optimal preconditioner] \label{exple:orientation}
Consider two strongly convex quadratics
\[ f_1 (x) = \tfrac{1}{2}\langle x, \Lambda x\rangle, \quad f_2 (x) = \tfrac{1}{2}\langle x, T_n x\rangle, 
\qquad \tiny T_n = \left(\begin{array}{cccc}
  2 & -1 & & \\
  -1 & \ddots & \ddots & \\
  & \ddots & \ddots & -1 \\
  & & -1 & 2
\end{array}\right), \]
and $\Lambda$ has eigenvalues of $T_n$ on the diagonal $\lambda_k  = 4 \sin^2 \big[ \tfrac{\pi k}{2 (n + 1)} \big]$ for $k=1, \ldots, n$. 
The functions $f_1$ and $f_2$ are related by a rotation: they share the same spectrum, and thus the same smoothness constant, strong convexity constant, and condition number $\kappa (T_n) = \Theta(n^2)$.
However, their optimal condition numbers with respect to diagonal stepsizes are different:
\begin{align}
  \kappa^{\star}_{f_1} = 1 & \quad \text{with optimal diagonal stepsize }
  P_1^{\star} = \Lambda^{-1} ; \nonumber\\
  \kappa^{\star}_{f_2} = \kappa (T_n) & \quad \text{with optimal diagonal stepsize } P_2^{\star} = I_n. \nonumber
\end{align}
The optimal diagonal preconditioner for $f_1$ gives perfect conditioning $\kappa_{f_1}^{\star} = 1$, yet no diagonal stepsize improves the conditioning of $f_2$. The orientation of eigenvectors affects the optimal condition number $\kappa^{\star}$
when $\mathcal{P} \neq \mathbb{R}^{n \times n}$.
\end{exple}

%% file: proofs/proof-property-ux.tex
By the chain rule, the gradient of $u_x (P)$ with respect to $P$ takes the form
\begin{equation} \label{eqn:property-ux-definition}
	\nabla u_x (P) = - \nabla f (x - P \nabla f (x)) \nabla f (x)^{\top}.
\end{equation}
Since $u_x (P) = f (x - P \nabla f (x))$ is the composition
between affine function $x - P \nabla f (x)$ and convex function $f$, it is convex by the composition rule.
Next, we show the smoothness of $u_x(P)$. For any $P_1, P_2$, we deduce that
\begin{align}
  \| \nabla u_x (P_1) - \nabla u_x (P_2) \|_F & = \| [\nabla f (x - P_1 \nabla
  f (x)) - \nabla f (x - P_2 \nabla f (x))] \nabla f (x)^{\top} \|_F
  \tag{by definition \eqref{eqn:property-ux-definition}} \\
  & = \| \nabla f (x - P_1 \nabla f (x)) - \nabla f (x - P_2 \nabla f (x)) \|
  \| \nabla f (x) \| \tag{by $\|a b^\top\|_F = \|a\| \|b\|$} \\
  & \leq L \| (P_1 - P_2) \nabla f (x) \| \| \nabla f (x) \| \tag{by $L$-smoothness}\\
  & \leq L \| \nabla f (x) \|^2  \| P_1 - P_2 \| \tag{by submultiplicativity of $\|\cdot\|$}\\
  & \leq L \| \nabla f (x) \|^2  \| P_1 - P_2 \|_F. \tag{by $\|\cdot\| \leq \|\cdot\|_F$}
\end{align}
Now, suppose $\diam(\mathcal{P}) \leq D$. We show the Lipschitz continuity of $u_x (P)$. 
For any $P \in \mathcal{P}$, we deduce that
\begin{align}
  \| \nabla u_x (P) \|_F & = \| \nabla f (x - P \nabla f (x)) \nabla f
  (x)^{\top} \|_F \tag{by definition \eqref{eqn:property-ux-definition}} \\
  & = \| \nabla f (x - P \nabla f (x)) \| \| \nabla f (x) \| \tag{by $\|a b^\top\|_F = \|a\| \|b\|$} \\
  & \leq (\| \nabla f (x - P \nabla f (x)) - \nabla f (x) \| + \| \nabla f
  (x) \|) \| \nabla f (x) \| \tag{by $\|a\| \leq \|a - b\| + \|b\|$}\\
  & \leq (L \| P \nabla f (x) \| + \| \nabla f (x) \|) \| \nabla f (x) \|
  \tag{by $L$-smoothness} \\
  & \leq (L \| P \| + 1) \| \nabla f (x) \|^2 \tag{by submultiplicativity of $\|\cdot\|$} \\
  & \leq (L D + 1) \| \nabla f (x) \|^2. \tag{by $\|P\|\leq \diam(\mathcal{P}) \leq D$}
\end{align}

%% file: proofs/proof-property-rx.tex
When $x \nin \mathcal{X}^{\star}$, the denominator $f (x) - f^{\star} > 0$
and the ratio feedback is well-defined.
Note that $r_x(P)$ simply translates and scales $u_x(P)$ by a positive factor $f (x) - f^{\star} > 0$ (since $x \nin \mathcal{X}^{\star}$):
\begin{equation*}
  r_x (P) = \tfrac{u_x (P) - f^{\star}}{f (x) - f^{\star}}.
\end{equation*}
Hence the gradient of $r_x(P)$ with respect to $P$ is $\nabla r_x (P) = \tfrac{\nabla u_x (P)}{f (x) - f^{\star}} =
-  \tfrac{\nabla f (x - P \nabla f (x)) \nabla f (x)^{\top}}{f (x) - f^{\star}}$. Non-negativity of $r_x (P)$ follows from the fact that $f(x) > f^{\star}$ for all $x \nin \Xcal^\star$ and the numerator is non-negative.
The convexity of $r_x (P)$ follows from the convexity of $u_x (P)$ in \tmtextbf{Proposition \ref{prop:property-ux}}. 
Since $u_x (P)$ is $L \| \nabla f (x) \|^2$-smooth, the ratio feedback $r_x (P)$ is also smooth with smoothness constant $\tfrac{L \| \nabla f (x) \|^2}{f (x) - f^{\star}} \leq 2 L^2$.
Suppose further $\diam(\mathcal{P}) \leq D$. 
Since $u_x (P)$ is $(LD + 1)\| \nabla f(x)\|^2$-Lipschitz, the ratio feedback has Lipschitz constant $\tfrac{(LD + 1)\| \nabla f(x)\|^2}{f (x) - f^{\star}} \leq 2L(LD + 1)$.

%% file: proofs/proof-property-hx.tex
Hypergradient feedback $h_x(P) = \tfrac{u_x(P) - f(x)}{\| \nabla f(x) \|^2}$ simply translates $u_x(P)$ and divides by a factor of $\|\nabla f(x)\|^2$. 
The convexity, smoothness, and Lipschitz continuity of $h_x(P)$ follow immediately from the same properties of $u_x(P)$ in \Cref{prop:property-ux}.
The expression of gradient $\nabla h_x(P)$, smoothness constant, and Lipschitz constant are obtained by dividing those of $u_x(P)$ by the scaling factor $\| \nabla f(x) \|^2$.

%% file: proofs/proof-opt-preconditioner.tex
For any $x, y \in \Rbb^n$ and $t \in [0, 1]$, define $x_t \assign x + t
(y - x)$. Since $(P^{\star}_+)^{1 / 2} \nabla^2 f (z) (P^{\star}_+)^{1 / 2}
\succeq \tfrac{1}{\kappa^\star} I$ for any $z \in \mathbb{R}^n$,
\begin{align}
  f (y) & = f (x) + \langle \nabla f (x), y - x \rangle + \textstyle \int^1_0 (1 - t)
  \langle y - x, \nabla^2 f (x_t) (y - x) \rangle \mathd t \nonumber\\
  & \geq f (x) + \langle \nabla f (x), y - x \rangle +
  \tfrac{1}{\kappa^{\star}} \| x - y \|^2_{(P^{\star}_+)^{- 1}} \textstyle \int^1_0 (1 -
  t) \mathd t \tag{by $(P^{\star}_+)^{1 / 2} \nabla^2 f (z) (P^{\star}_+)^{1 / 2}
\succeq \tfrac{1}{\kappa^\star} I$} \\
  & = f (x) + \langle \nabla f (x), y - x \rangle + \tfrac{1}{2
  \kappa^{\star}} \| x - y \|^2_{(P^{\star}_+)^{- 1}} . \nonumber
\end{align}
Minimize both sides over $y \in \mathbb{R}^n$ to obtain
\[ f^{\star} \geq f (x) + \min_{d \in \mathbb{R}^n} \langle \nabla f (x),
   d \rangle + \tfrac{1}{2 \kappa^{\star}} \| d \|^2_{(P^{\star}_+)^{- 1}} = f
   (x) - \tfrac{\kappa^{\star}}{2} \| \nabla f (x) \|_{P^{\star}_+}^2, \]
which rearranges to 
\begin{equation} \label{eqn:pf-opt-preconditioner-bound-1}
  f (x) - f^{\star} \leq \tfrac{\kappa^{\star}}{2} \| \nabla f (x) \|_{P^{\star}_+}^2.
\end{equation}
Using $(P^{\star}_+)^{1 / 2} \nabla^2 f (z) (P^{\star}_+)^{1 / 2} \preceq I$, we have
\begin{align}
  f (x - P^{\star}_+ \nabla f (x)) ={} & f (x) - \langle \nabla f (x),
  P^{\star}_+ \nabla f (x) \rangle + \textstyle \int^1_0 (1 - t) \langle \nabla f (x),
  P^{\star}_+ \nabla^2 f (x_t) P^{\star}_+ \nabla f (x) \rangle \mathd t
  \nonumber\\
  \leq{} & f (x) - \langle \nabla f (x), P^{\star}_+ \nabla f (x) \rangle +
  \textstyle \int^1_0 (1 - t) \| \nabla f (x) \|_{P^{\star}_+}^2 \mathd t \tag{by $(P^{\star}_+)^{1 / 2} \nabla^2 f (z) (P^{\star}_+)^{1 / 2} \preceq I$}\\
  ={} & f (x) - \| \nabla f (x) \|_{P^{\star}_+}^2 + \tfrac{1}{2} \| \nabla f
  (x) \|_{P^{\star}_+}^2 \nonumber\\
  ={} & f (x) - \tfrac{1}{2} \| \nabla f (x) \|_{P^{\star}_+}^2 \nonumber\\
  \leq{} & f (x) - \tfrac{1}{\kappa^{\star}} [f (x) - f^{\star}]. \tag{by \eqref{eqn:pf-opt-preconditioner-bound-1}}
\end{align}
Rearrange to conclude $f (x - P^{\star}_+ \nabla f (x)) - f^{\star} \leq ( 1 - \tfrac{1}{\kappa^{\star}} ) [f (x) - f^{\star}]$.

%% file: proofs/proof-minimax-feedback.tex
\begin{itemize}[leftmargin=10pt]
\item Ratio feedback: Definition of $P_r^{\star}$ and $P_+^{\star} \in \mathcal{P}$ together imply $r_x (P_r^{\star}) \leq \max_{x \nin \Xcal^\star } r_x (P^{\star}_+) \overset{\eqref{eq:opt-preconditioner}}{\leq} 1 - \tfrac{1}{\kappa^{\star}}$.

  If $f$ is a strongly convex quadratic with Hessian $A \succ 0$, then both global optimal preconditioner $P^{\star}_+$ and the minimax stepsize $P_r^{\star}$ are Hessian inverse with optimal condition number $\kappa^{\star} = 1$ and minimax feedback $r_x (A^{-1}) = 0$.

\item Hypergradient feedback: Since $\tfrac{1}{L} I \in \mathcal{P}$, the definition of $P_h^{\star}$ and descent lemma together imply
  \begin{equation*}
    h_x (P^{\star}_h) \leq \max_{x \nin \Xcal^\star } h_x (\tfrac{1}{L} I) \leq - \tfrac{1}{2 L}.
  \end{equation*}
Now we show that $P_h^\star = \tfrac{1}{L} I$ if $f$ is a strongly convex quadratic. Without loss of generality, assume $f(x) = \tfrac{1}{2} \langle x, A x \rangle$ is homogeneous. Then $x^{\star} = 0$ is the unique optimal solution and 
  \begin{align}
    \min_{P \in \mathcal{P}} \max_{x \neq 0} ~ h_x (P)
    ={} & \min_{P \in \mathcal{P}} \max_{x \neq 0} \tfrac{f(x - P \nabla f(x)) - f(x)}{\| \nabla f(x) \|^2} \nonumber\\
    \geq{} &\min_{P \in \Rbb^{n\times n}}  \max_{x \neq 0} \tfrac{f(x - P \nabla f(x)) - f(x)}{\| \nabla f(x) \|^2} \tag{by $\Pcal \subseteq \Rbb^{n \times n}$} \\
    \geq{} & \max_{x \neq 0} \min_{P \in \Rbb^{n\times n}} \tfrac{f(x - P \nabla f(x)) - f(x)}{\| \nabla f(x) \|^2} \tag{by weak duality} \\
    ={} & \max_{x \neq 0} \tfrac{f^{\star} - f(x)}{\| \nabla f(x) \|^2} \label{eqn:minimax-feedback-1} \\
    ={} & \max_{x \neq 0} - \tfrac{\langle x, A x \rangle}{2 \| A x
    \|^2} \tag{by $f^{\star} = 0$ and definition of $f$}\\
    ={} & \max_{x \neq 0} \tfrac{\langle A x, (-A^{-1}) A x \rangle}{2 \| A x
    \|^2}
    ={} \tfrac{1}{2} \lambda_{\max} (- A^{- 1}) = - \tfrac{1}{2 L}, \nonumber
  \end{align}
where \eqref{eqn:minimax-feedback-1} plugs in the minimizer $P = A^{-1}$ that drives $x$ to optimal solution $x^\star$ in one step. 
On the other hand, descent lemma guarantees $\max_{x \neq 0} h_x ( \tfrac{1}{L} I ) \leq - \tfrac{1}{2 L}$. Hence, $\tfrac{1}{L} I$ achieves the minimax feedback and is the minimax optimal stepsize with respect to hypergradient feedback.
\end{itemize}

%% file: sec_action.tex
\section{Progress reduction and landscape action} \label{sec:action}

The example algorithm in \Cref{sec:intro} is based on a simple observation that the suboptimality after $K$ iterations can be bounded by the sum of contraction ratios at each step (see \eqref{eqn:intro-reduction}).
This observation, which we call \emph{reduction}, reduces the optimization problem $\min_{x \in \mathbb{R}^n} f(x)$ to minimizing the sum of \emph{per iteration progress}. 
The per iteration progress with respect to ratio feedback and hypergradient feedback is measured by
\begin{equation*}
  r_k \assign \tfrac{f (x^{k + 1}) - f^{\star}}{f (x^k) - f^{\star}}, \qquad h_k \assign \tfrac{f (x^{k + 1}) - f (x^k)}{\| \nabla f (x^k) \|^2}.
\end{equation*}
Two reductions below guarantee that smaller cumulative progress $\sum_{k = 1}^K r_k$ or $\sum_{k = 1}^K h_k$ yields faster convergence.

\begin{thm}[Reduction for $r_k$]
  \label{thm:reduction-ratio}Let $\{x^k\}$ be any sequence such that $x^k \nin \Xcal^\star$. Then
  \[ \textstyle f (x^{K + 1}) - f^{\star} \leq [f (x^1) - f^{\star}] (
     \tfrac{1}{K} \sum_{k = 1}^K r_k )^K . \]
\end{thm}

\input{\input{proofs/proof-reduction-ratio.tex}}

\begin{thm}[Reductions for $h_k$]
  \label{thm:reduction-hypergrad}Let $\{x^k\}$ be any sequence such that $x^k \nin \Xcal^\star$ and $f(x^{k+1}) \leq f(x^k)$. 
  \begin{itemize}[leftmargin=10pt,itemsep=0pt,topsep=5pt]
    \item If $f$ is convex, then
    \[ \textstyle f (x^{K + 1}) - f^{\star} \leq \min \{
       \tfrac{\Delta^2}{K} \frac{1}{\frac{1}{K}\sum^K_{k = 1} - h_k}, f (x^1) - f^{\star}
       \}, \]
    where $\Delta \assign \max_{x \in \{x: f(x) \leq f(x^1)\}} \min_{x^{\star}
    \in \mathcal{X}^{\star}} \| x - x^{\star} \|$.
    
    \item If $f$ is $\mu$-strongly convex, then
    \[ \textstyle f (x^{K + 1}) - f^{\star} \leq [f (x^1) - f^{\star}] ( 1 -
       \tfrac{2 \mu}{K} \sum^K_{k = 1} -h_k )^K . \]
  \end{itemize}
\end{thm}

\input{\input{proofs/proof-reduction-hypergrad.tex}}

These two reductions are blackbox: they are independent of the mechanism generating the iterates $\{x^k\}$ and universally apply to (monotone) algorithms other than {\osgm}. 
\Cref{thm:reduction-ratio} even does not require $f$ to be convex.
\Cref{thm:reduction-hypergrad} comes with a restriction of non-increasing function values ${f(x^k)}$.
In {\osgm}, when the landscape always accepts the scheduler's proposal by taking $x^{k+1} = x^{k+1/2}$, the algorithm may not be monotone, and the reductions for $h_k$ in \Cref{thm:reduction-hypergrad} may not apply. However, we can design the landscape agent so that it takes action that filters out bad stepsizes and ensures monotonicity of the iterates.  

\subsection{Landscape actions $\mathcal{M}$ and progress}
The landscape action in {\osgm} is denoted by $x^{k + 1} =\mathcal{M} (x^{k + 1 / 2}, x^k)$.
We consider four landscape actions below.
\begin{itemize}[leftmargin=10pt]
  \item {\textit{Vanilla.}} $x^{k + 1} = x^{k + 1 / 2}$
  
  Landscape accepts all decisions from the scheduler.
  
  \item {\textit{Monotone.}} $x^{k+1}$ satisfies $f(x^{k + 1}) \leq \min \{ f (x^{k+1/2}), f(x^k) \}$
  
  Landscape moves to $x^{k+1}$ that is no worse than the current iterate $x^k$ and the suggested iterate $x^{k+1/2}$ in terms of objective value. This action can be implemented by line search or a null step: $x^{k+1} = \argmin_{x \in \{x^{k+1/2}, x^k\}} f(x)$.
  
  \item {\textit{Lookahead.}} $x^{k + 1} = x^{k + 1 / 2} -
  \tfrac{1}{L} \nabla f (x^{k + 1 / 2})$
  
  Landscape takes an additional gradient descent step on top of the suggested iterate $x^{k + 1 / 2}$ but does not enforce monotonicity.
  
  \item {\textit{Monotone Lookahead.}} $x^{k+1}$ satisfies $f(x^{k + 1}) \leq \min \{ f (x^{k+1/2} - \frac{1}{L}\nabla f(x^{k + 1/2})), f(x^k) \}$  

Landscape moves to $x^{k+1}$ that is no worse than the current iterate $x^k$ and the lookahead iterate $x^{k+1/2} - \frac{1}{L}\nabla f(x^{k + 1/2})$ in terms of objective value. 

\end{itemize}

All these actions, except vanilla, improve the iteration progress $r_k$ and $h_k$ upon the feedback $r_{x^k}(P_k)$ and $h_{x^k}(P_k)$ achieved by the scheduler, which we illustrate in \Cref{lem:feedback-progress} below.

\begin{lem}[Feedback and progress]
  \label{lem:feedback-progress}
  Let $f$ be convex and $L$-smooth. Each of the above four landscape actions guarantees the following relation between the feedback and per iteration progress:
  \begin{itemize}[leftmargin=10pt]
    \item Vanilla. $r_k = r_{x^k} (P_k)$ and $h_k = h_{x^k}
    (P_k)$.
    
    \item Monotone. $r_k \leq \min \{ r_{x^k} (P_k), 1 \}$ and $h_k \leq \min \{ h_{x^k} (P_k), 0 \}$.
    
    \item Lookahead. $r_k \leq r_{x^k} (P_k) - \tfrac{1}{4 L^2}
    \| \nabla r_{x^k} (P_k) \|_F^2$ and $h_k \leq h_{x^k} (P_k) - \tfrac{1}{2
    L} \| \nabla h_{x^k} (P_k) \|_F^2$.
    
    \item Monotone Lookahead. $r_k \leq \min \{ r_{x^k}
    (P_k) - \tfrac{1}{4 L^2} \| \nabla r_{x^k} (P_k) \|_F^2, 1 \}$ and
    $h_k \leq \min \{ h_{x^k} (P_k) - \tfrac{1}{2 L} \| \nabla h_{x^k}
    (P_k) \|_F^2, 0 \}$.
  \end{itemize}
\end{lem}

\input{proofs/proof-feedback-progress.tex}

The per iteration progress $r_k$ and $h_k$ under these landscape actions are never worse than the feedback $r_{x^k}(P_k)$ and $h_{x^k}(P_k)$.
To accelerate convergence, 

it suffices for the scheduler to minimize the cumulative feedback $\sum_{k = 1}^K r_{x^k} (P_k)$ and $\sum_{k = 1}^K h_{x^k} (P_k)$, a task well-suited for online learning algorithms.

%% file: proofs/proof-reduction-ratio.tex
The result follows immediately from the definition of $r_k$ and the AM-GM inequality:
\[ \textstyle \tfrac{f (x^{K + 1}) - f^{\star}}{f (x^1) - f^{\star}} = \prod_{k
   = 1}^K \tfrac{f (x^{k + 1}) - f^{\star}}{f (x^k) - f^{\star}} =
   \prod_{k = 1}^K r_k \leq ( \tfrac{1}{K}  \sum_{k = 1}^K r_k )^K.
\]

%% file: proofs/proof-reduction-hypergrad.tex
\paragraph{Convex $f$.} Observe that
\begin{align}
  \tfrac{1}{f (x^{K + 1}) - f^{\star}} & = \textstyle \sum_{k = 1}^K \Big[
  \tfrac{1}{f (x^{k + 1}) - f^{\star}} - \tfrac{1}{f (x^k) - f^{\star}} \Big] + \tfrac{1}{f (x^1) - f^{\star}} \nonumber\\
  & = \textstyle \sum_{k = 1}^K \Big[ \tfrac{f (x^k) - f (x^{k + 1})}{[f (x^{k + 1}) -
  f^{\star}] [f (x^k) - f^{\star}]} \Big] + \tfrac{1}{f (x^1) - f^{\star}} \nonumber\\
  & = \textstyle \sum_{k = 1}^K \Big[ \tfrac{- h_k \| \nabla f (x^k) \|^2}{[f (x^{k
  + 1}) - f^{\star}] [f (x^k) - f^{\star}]} \Big] + \tfrac{1}{f
  (x^1) - f^{\star}} . \label{eqn:pf-reduction-hypergrad-1}
\end{align}
Since $h_k = \tfrac{f(x^{k+1}) - f(x^{k})}{\|\nabla f(x^{k})\|^2} \leq 0$ by monotonicity $f(x^{k+1}) \leq f(x^k)$ and $f (x) - f^{\star} \leq \| \nabla f (x) \| \cdot \| x - x^{\star}\|$ by convexity of $f$, we have
\[ \tfrac{- h_k \| \nabla f (x^k) \|^2}{[f (x^{k + 1}) - f^{\star}] [f
   (x^k) - f^{\star}]} \geq \tfrac{- h_k \| \nabla f (x^k) \|^2}{[f (x^k)
   - f^{\star}]^2} \geq \tfrac{- h_k}{\tmop{dist} (x^k,
   \mathcal{X^{\star}})^2} \geq \tfrac{- h_k}{\Delta^2} . \]
Hence, \eqref{eqn:pf-reduction-hypergrad-1} can be lower bounded by
\begin{equation*}
  \tfrac{1}{f (x^{K + 1}) - f^{\star}} \geq - \tfrac{1}{\Delta^2}
  \textstyle \sum_{k = 1}^K h_k + \tfrac{1}{f (x^1) - f^{\star}} \geq \max \{ - \tfrac{1}{\Delta^2} \textstyle \sum_{k = 1}^K h_k, \tfrac{1}{f
  (x^1) - f^{\star}} \}.
\end{equation*}

The desired result follows by taking the reciprocal on both sides of the
inequality.

\paragraph{$\mu$-strongly convex $f$.}
The desired inequality follows from the following steps:
\begin{align}
  \tfrac{f (x^{K + 1}) - f^{\star}}{f (x^1) - f^{\star}}
  \leq{} & ( \tfrac{1}{K} \textstyle \sum_{k=1}^K \tfrac{f(x^{k+1}) - f^{\star}}{f(x^{k}) - f^{\star}} )^K \tag{by AM-GM inequality} \\
  ={} & ( 1 + \tfrac{1}{K} \textstyle\sum_{k = 1}^K \tfrac{f(x^{k+1})-f(x^k)}{f (x^k) - f^{\star}} )^K
  \nonumber\\
  ={} & ( 1 + \tfrac{1}{K} \textstyle\sum_{k = 1}^K h_k \tfrac{\|\nabla f(x^k)\|^2}{f (x^k) - f^{\star}} )^K
  \tag{by definition of $h_k$}\\
  \leq{} & ( 1 - \tfrac{2 \mu}{K} \textstyle\sum_{k = 1}^K -h_k )^K. \tag{by $h_k \leq 0$ and $\tfrac{\|\nabla f(x^k)\|^2}{f (x^k) - f^{\star}} \geq 2 \mu$}
\end{align}

%% file: proofs/proof-feedback-progress.tex
Recall that $x^{k+1} = \mathcal{M}(x^{k+1/2}, x^k)$ and $x^{k+1/2} = x^k - P_k \nabla f(x^k)$. The per iteration progress and feedback are
\begin{equation*}
  r_k = \tfrac{f (x^{k + 1}) - f^{\star}}{f (x^k) - f^{\star}}, \quad
  r_{x^k} (P_k) = \tfrac{f (x^{k + 1 / 2}) - f^{\star}}{f (x^k) - f^{\star}}, \quad
  h_k = \tfrac{f (x^{k + 1}) - f (x^k)}{\| \nabla f (x^k) \|^2}, \quad
  h_{x^k} (P_k) = \tfrac{f (x^{k + 1 / 2}) - f (x^k)}{\| \nabla f (x^k) \|^2}.
\end{equation*}

\begin{itemize}[leftmargin=10pt]
\item \tmtextit{Vanilla} landscape satisfies $x^{k + 1} = x^{k + 1 / 2}$. Then  $r_k = r_{x^k} (P_k)$ and $h_k = h_{x^k}(P_k)$ follow immediately.

\item \tmtextit{Monotone} landscape satisfies $f (x^{k + 1}) \leq \min \{ f (x^{k + 1 / 2}), f(x^k) \}$ and thus
\begin{align}
r_k & = \tfrac{f (x^{k + 1}) - f^{\star}}{f (x^k) - f^{\star}} \leq
\tfrac{\min \{ f (x^{k + 1 / 2}), f (x^k) \} - f^{\star}}{f (x^k) - f^{\star}} = \min \{ r_{x^k} (P_k), 1 \} ; \nonumber\\
h_k & = \tfrac{f (x^{k + 1}) - f (x^k)}{\| \nabla f (x^k) \|^2} \leq
\tfrac{\min \{ f (x^{k + 1 / 2}), f (x^k) \} - f (x^k)}{\| \nabla f (x^k)
\|^2} = \min \{ h_{x^k} (P_k), 0 \} . \nonumber
\end{align}
\item \tmtextit{Lookahead} landscape satisfies $x^{k + 1} = x^{k + 1 / 2} -  \tfrac{1}{L} \nabla f (x^{k + 1 / 2})$ and descent lemma implies
\[ f (x^{k + 1}) = f ( x^{k + 1 / 2} - \tfrac{1}{L} \nabla f (x^{k + 1 /
  2}) ) \leq f (x^{k + 1 / 2}) - \tfrac{1}{2 L} \| \nabla f (x^{k + 1 /  2}) \|^2. \]
Then it follows that
\begin{align}
r_k = \tfrac{f (x^{k + 1}) - f^{\star}}{f (x^k) - f^{\star}} \leq{} &
\tfrac{f (x^{k + 1 / 2}) - f^{\star}}{f (x^k) - f^{\star}} -
\tfrac{1}{2 L} \tfrac{\| \nabla f (x^{k + 1 / 2}) \|^2}{f (x^k) - f^{\star}} \nonumber\\
={} & \tfrac{f (x^{k + 1 / 2}) - f^{\star}}{f (x^k) - f^{\star}} -
\tfrac{1}{2 L} \tfrac{\| \nabla f (x^{k + 1 / 2}) \|^2 \| \nabla f (x^k)
\|^2}{[f (x^k) - f^{\star}]^2} \tfrac{f (x^k) - f^{\star}}{\| \nabla
f (x^k) \|^2} \nonumber\\
={} & r_{x^k} (P_k) - \tfrac{1}{2 L}  \| \nabla r_{x^k} (P_k) \|_F^2  \tfrac{f
(x^k) - f^{\star}}{\| \nabla f (x^k) \|^2} \label{eqn:pf-feedback-progress-1}\\
\leq{} & r_{x^k} (P_k) - \tfrac{1}{4 L^2} \| \nabla r_{x^k} (P_k) \|_F^2, \label{eqn:pf-feedback-progress-2}
\end{align}
where \eqref{eqn:pf-feedback-progress-1} uses the relation
\[ \| \nabla r_{x^k} (P_k) \|^2_F = \big\| \tfrac{\nabla f (x^{k + 1 / 2})
  \nabla f (x^k)^{\top}}{f (x^k) - f^{\star}} \big\|^2_F = \tfrac{\|
  \nabla f (x^{k + 1 / 2}) \|^2 \| \nabla f (x^k) \|^2}{[f (x^k) - f^{\star}]^2} \]
and the inequality in \eqref{eqn:pf-feedback-progress-2} uses the fact $f(x) - f^{\star} \geq \tfrac{1}{2L} \|\nabla f(x)\|^2$ for $L$-smooth $f$.
Similarly,
\begin{align}
h_k = \tfrac{f (x^{k + 1}) - f (x^k)}{\| \nabla f (x^k) \|^2} \leq{} & \tfrac{f
(x^{k + 1 / 2}) - f^{\star}}{\| \nabla f (x^k) \|^2} - \tfrac{1}{2 L}
\tfrac{\| \nabla f (x^{k + 1 / 2}) \|^2}{\| \nabla f (x^k) \|^2} \nonumber\\
={} & h_{x^k} (P_k) - \tfrac{1}{2 L} \| \nabla h_{x^k} (P_k) \|_F^2, \label{eqn:pf-feedback-progress-3}
\end{align}
where \eqref{eqn:pf-feedback-progress-3} uses the relation
\[ \| \nabla h_{x^k} (P_k) \|^2_F = \Big \| \tfrac{\nabla f (x^{k + 1 / 2})
  \nabla f (x^k)^{\top}}{\| \nabla f (x^k) \|^2} \Big\|^2_F = \tfrac{\|
  \nabla f (x^{k + 1 / 2}) \|^2}{\| \nabla f (x^k) \|^2} . \]
\item \tmtextit{Monotone Lookahead} landscape satisfies
\[ f (x^{k + 1}) \leq \min \{ f ( x^{k + 1 / 2} - \tfrac{1}{L} \nabla
  f (x^{k + 1 / 2}) ), f (x^k) \}. \]
Combining the analysis of lookahead and monotone landscape completes the proof.

\end{itemize}

%% file: sec_oco.tex
\section{Stepsize update and online learning} \label{sec:oco}

The scheduler aims to generate a sequence of stepsizes $\{ P_k \}$ that reduces the cumulative feedback $\sum_{k = 1}^K \ell_{x^k} (P_k)$ as much as possible. This can be done by the existing online learning algorithms with sublinear regret guarantees:

\begin{equation} \label{eqn:regret-single}
   \textstyle  \tfrac{1}{K}\sum_{k = 1}^K \ell_{x^k} (P_k) \leq  \tfrac{1}{K} \sum_{k = 1}^K \ell_{x^k} (\hat{P}) + o (1) \quad \text{ for any } \hat{P} \in \mathcal{P}.
\end{equation}
This guarantee \eqref{eqn:regret-single} says that the average feedback $\tfrac{1}{K} \sum_{k = 1}^K \ell_{x^k} (P_k)$ achieved by the scheduler is asymptotically no worse than that achieved by any fixed stepsize $\hat{P} \in \Pcal$. 
We consider the online learning algorithm $\Acal$ to be \emph{online gradient descent} in this paper for concreteness, but our arguments for the convergence of {\osgm} apply to other online learning algorithms as well. We will discuss the practical choice of online algorithm in \Cref{sec:practical}.

\subsection{Regret guarantees of online gradient descent}
Online gradient descent ({\ogd}) updates the stepsize by
\begin{equation} \label{eqn:online-gradient-descent}
   P_{k + 1} = \Pi_{\mathcal{P}} [P_k - \eta_k \nabla \ell_{x^k} (P_k)],
\end{equation}
where $\{\eta_k\}$ is a sequence of non-negative and non-increasing online learning stepsizes which will be specified later and $\Pi_{\mathcal{P}}$ is the projection operator onto the candidate set $\mathcal{P}$. {\ogd} on convex feedback $\ell_{x^k}$ satisfies

\begin{equation} \label{eqn:ogd-Pk-ineq}
   \| P_{k + 1} - \hat{P} \|_F^2 \leq \| P_k - \hat{P} \|_F^2 - 2 \eta_k [\ell_{x^k} (P_k) - \ell_{x^k} (\hat{P})] + \eta_k^2 \| \nabla \ell_{x^k} (P_k) \|_F^2 \quad \text{ for any } \hat{P} \in \mathcal{P}.
\end{equation}
According to \eqref{eqn:ogd-Pk-ineq}, whenever the scheduler makes a
stepsize decision $P_k$ that underperforms $\hat{P}$ (i.e., $\ell_{x^k} (P_k) -
\ell_{x^k} (\hat{P}) > 0$), the refined stepsize $P_{k + 1}$ approaches $\hat{P}$ up to some
error $\eta_k^2 \| \nabla \ell_{x^k} (P_k) \|_F^2$. Rearranging the inequality \eqref{eqn:ogd-Pk-ineq} and telescoping from $k=1$ to $K$, we obtain the following regret guarantee:

\begin{lem}[Static regret]
\label{lem:ogd-sqrtK-eta} Suppose $\{ \ell_{x^k} (P) \}$ are convex in
  $P$. Then {\ogd} with constant stepsize $\eta_k \equiv \eta > 0$
  generates a sequence $\{ P_k \}$ such that
\begin{equation} \label{eqn:ogd-regret-1}
\textstyle \sum_{k = 1}^K \ell_{x^k} (P_k) \leq \sum_{k = 1}^K \ell_{x^k} (\hat{P}) + \tfrac{1}{2 \eta} \| P_1 - \hat{P} \|_F^2 + \frac{\eta}{2} \sum_{k =
   1}^K \| \nabla \ell_{x^k} (P_k) \|_F^2 \text{~~ for any $\hat{P} \in
   \mathcal{P}$} .
\end{equation}
Besides, if $\{ \ell_{x^k} (P) \}$ are $\sigma$-Lipschitz w.r.t. $\| \cdot \|_F$ and $\diam (\mathcal{P}) \leq D$, then {\ogd} with $\eta_k \equiv \frac{c}{\sqrt{K}}$ or $\eta_k = \frac{c}{\sqrt{k}}$ satisfies
\begin{equation} \label{eqn:ogd-regret-2}
\textstyle  \sum_{k = 1}^K \ell_{x^k} (P_k) \leq \sum_{k = 1}^K \ell_{x^k} (\hat{P}) + ( \tfrac{D^2}{2 c} + c \sigma^2 ) \sqrt{K} \text{~~ for any $\hat{P} \in \mathcal{P}$}.
\end{equation}
\end{lem}

\Cref{lem:ogd-sqrtK-eta} competes the adaptive stepsizes $\{ P_k \}$ against a static stepsize $\hat{P} \in \mathcal{P}$. A more refined analysis of online gradient descent, referred to as \emph{dynamic regret} \cite{hazan2016introduction, orabona2019modern}, enables the scheduler to compete with a sequence of stepsizes $\{ \hat{P}_k \}, \hat{P}_k \in \mathcal{P}$, at the cost of the \emph{path length} defined by 
\begin{equation}\label{eqn:pathlength}
\textstyle \mathsf{PL} (\{ \hat{P}_k \}) \assign \sum_{k = 1}^{K-1} \| \hat{P}_k - \hat{P}_{k + 1} \|_F.
\end{equation}

\begin{lem}[Dynamic regret]
\label{lem:dynamic-sqrtK}Suppose $\{ \ell_{x^k} (P) \}$ are convex in
$P$. For any benchmark sequence of stepsizes $\{ \hat{P}_k \},  P_k \in \mathcal{P}$, {\ogd} with constant stepsize $\eta > 0$ generates $\{ P_k \}$ such that
\begin{equation}\label{eqn:dynamic-regret-1}
\textstyle \sum_{k = 1}^K \ell_{x^k} (P_k) \leq \textstyle \sum_{k = 1}^K \ell_{x^k}
   (\hat{P}_k) + \frac{\eta}{2} \textstyle \sum_{k = 1}^K \| \nabla \ell_{x^k} (P_k) \|_F^2 + \tfrac{\| \hat{P}_{K} - P_1 \|^2_F}{2 \eta} + \tfrac{\max_{k \leq K} \| P_{k} - P_1 \|_F}{\eta} \mathsf{PL} (\{
   \hat{P}_k \}).
\end{equation}
Besides, if $\{ \ell_{x^k} (P) \}$ are $\sigma$-Lipschitz w.r.t. $\| \cdot \|_F$ and  $\diam (\mathcal{P}) \leq D$, then {\ogd} with $\eta_k \equiv \frac{c}{\sqrt{K}}$ satisfies
\begin{equation}\label{eqn:dynamic-regret-2}
\textstyle \sum_{k = 1}^K \ell_{x^k} (P_k) \leq \textstyle \sum_{k = 1}^K \ell_{x^k}
   (\hat{P}_k) + [ \tfrac{c \sigma^2}{2} + \tfrac{D^2}{2 c} +
   \tfrac{D}{c} \mathsf{PL} (\{ \hat{P}_k \}) ] \sqrt{K}
   \text{ \ for any $\hat{P}_k \in \mathcal{P}$} .
\end{equation}
\end{lem}

When $\hat{P_k} \equiv \hat{P}$, the path length vanishes $\mathsf{PL} (\{ \hat{P}_k \}) = 0$ and dynamic regret (\Cref{lem:dynamic-sqrtK}) reduce to static regret (\Cref{lem:ogd-sqrtK-eta}).
The relation \eqref{eqn:dynamic-regret-1} holds for \emph{any} sequence $\{ \hat{P}_k \},  \hat{P}_k \in \mathcal{P}$, indicating that {\ogd} is asymptotically competitive with the optimal sequence $\{ \hat{P}_k \}$ that minimizes the right-hand side of the bounds.
The optimal sequence balances the trade-off between the cumulative feedback $\sum_{k = 1}^K \ell_{x^k} (\hat{P}_k)$ and the path length $\mathsf{PL} (\{ \hat{P}_k \})$. \Cref{exple:local} illustrates the effect of dynamic regret guarantees.

\begin{exple} \label{exple:local}
   Consider diagonal stepsizes $\{P_k\}, P_k \in \mathcal{P} \assign \{P = \diag(d): d \in \mathbb{R}^n \}$ for a 2-dimensional smooth convex objective defined as
   \[ f (x_1, x_2) = 
   \begin{cases}
   \tfrac{1}{4} x_1^2 + \tfrac{1}{2} x_2^2, & \text{ if } (x_1, x_2) \in R_1 \assign \{ (x_1, x_2) : x_1 \geq 0 \} ;\\
   \tfrac{3}{4} x_1^2 + \tfrac{1}{2} x_2^2, & \text{ if } (x_1, x_2) \in R_2 \assign \{ (x_1, x_2)
   : x_1 < 0 \}.
   \end{cases}
   \]
   \begin{minipage}{0.58\textwidth}
   The function $f$ is strongly convex and piecewise quadratic with different curvatures in the regions $R_1$ and $R_2$.
   For $x \in R_1$, the Hessian inverse $\diag(0.5, 1)$ achieves zero ratio feedback $r_x(\diag(0.5, 1)) = 0$ since it sends any $x \in R_1$ to the optimal solution $x^\star = (0, 0)$ in one preconditioned gradient step, resulting in $f(x - \diag(0.5, 1) \nabla f(x)) - f^{\star} = 0$. 
   Similarly, for $x \in R_2$, the Hessian inverse $\diag(1.5, 1)$ also achieves zero ratio feedback.
   The trajectory $\{ x^k \}$ of gradient descent (with appropriate stepsize) remains in the same region as the initial point $x^1$, ensuring that the common Hessian inverse on that region minimizes the cumulative feedback over the trajectory with an optimal value of zero: $\min_{P \in \Pcal} \sum_{k = 1}^K r_{x^k} (P) = 0$.\\
   \end{minipage}
   \hfill
   \begin{minipage}{0.4\textwidth}
      \centering
      \includegraphics[width=0.75\linewidth]{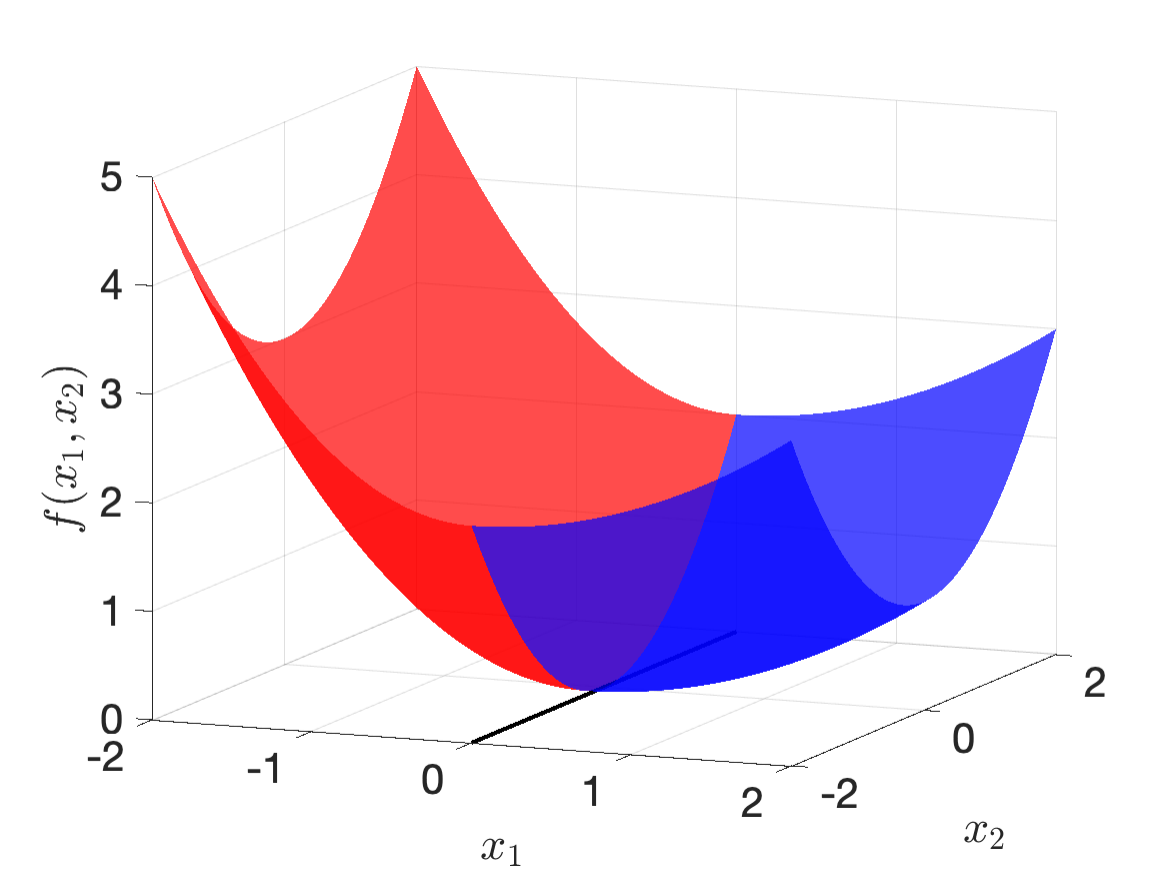}
      \captionof{figure}{Illustration of $f(x_1, x_2)$.}
      \label{fig:ex_pquad}
      \vspace{2em}
   \end{minipage}

   However, no diagonal stepsize can achieve $r_x(P) = 0$ for all $x$ simultaneously. Thus, the minimax optimal stepsize $P^\star_r$ must incur positive cumulative feedback $\sum_{k = 1}^K r_{x^k} (P_r^\star) > 0$, which is strictly worse than the trajectory-based bound.
   Even if the trajectory $\{ x^k \}$ does not always stay in the same region, the dynamic regret bound in \eqref{eqn:dynamic-regret-2} allows us to compare with a benchmark sequence $\{ \hat{P}_k \}$.
   For example, let $\hat{P}_k = \diag(0.5,1)$ if $x^k \in R_1$; and $\hat{P}_k = \diag(1.5, 1)$ if $x^k \in R_2$. Then $\textstyle \sum_{k=1}^{K} r_{x^k}(\hat{P}_k) = 0$ and the path length $\mathsf{PL} (\{ \hat{P}_k \}) = \sum_{k = 1}^{K-1} \| \hat{P}_k - \hat{P}_{k + 1} \|_F$ is proportional to the number of times the trajectory $\{ x^k \}$ switches from one region to the other.
\end{exple}

%% file: sec_algodemo.tex
\section{Algorithm design and analysis} \label{sec:algodemos}

This section demonstrates the convergence guarantees of {\osgm} through two variants, one with ratio feedback and the other with hypergradient feedback.
We denote each variant of {\osgm} as \texttt{\{Action\}} {\osgm\texttt{-\{Feedback\}}}, where \texttt{Action} refers to the landscape action and $\texttt{Feedback} \in \{\texttt{R}, \texttt{H}\}$ represents the initial letter of the type of feedback. Without loss of generality, we assume $x^k \nin \Xcal^\star$; otherwise the algorithm can be stopped.

\subsection{\osgmhandsonrx}\label{sec:handson-rx}

In this section, we assume $f$ is $L$-smooth and $\mu$-strongly convex  and instantiate {\osgm} with
\[ \ell_x (P) \assign r_x (P) \qquad \text{Lookahead landscape:~} x^{k + 1} = x^{k + 1 / 2}
   - \tfrac{1}{L} \nabla f (x^{k + 1 / 2}) \qquad \mathcal{A} \assign
   \text{Online gradient descent} . \]
The algorithm is called {\osgmhandsonrx} (\Cref{alg:demo}).
The candidate stepsize set is $\Pcal = \Rbb^{n\times n}$.

\begin{algorithm}[H]
{\textbf{input:} Initial point $x^1$, initial stepsize $P_1 \in \Pcal = \Rbb^{n\times n}$, online gradient stepsize $\eta_k \equiv \eta > 0$}\\
\For{$k = 1, 2, \dots$}{
      $x^{k + 1/2} = x^{k} - P_k \nabla f(x^k)$\\
      $x^{k + 1} = x^{k + 1/2} - \frac{1}{L} \nabla f(x^{k + 1/2})$\\
      $P_{k+1} = P_k - \eta \nabla r_{x^k}(P_k)$ \nllabel{alg:demo-1}
}
\caption{{\osgmhandsonrx}\label{alg:demo}}
\end{algorithm}
Recall that the gradient of ratio feedback in \Cref{alg:demo-1} of \Cref{alg:demo} takes the form
\begin{equation*}
   \nabla r_{x^k}(P_k) = - \tfrac{\nabla f(x^{k + 1/2})\nabla f(x^k)^\top}{f(x^k) - f^{\star}}.
\end{equation*}
The rest of the section analyzes the convergence behavior of {\osgmhandsonrx}, including global convergence and local superlinear convergence.

\subsubsection{Global convergence}

We analyze the convergence of {\osgmhandsonrx} in three steps.
First, according to \Cref{thm:reduction-ratio}, the convergence of {\osgmhandsonrx} follows from bounding the suboptimality in terms of the cumulative progress
\begin{equation} \label{eqn:algodemo-reduction}
	\textstyle f (x^{K + 1}) - f^{\star} \leq [f (x^1) - f^{\star}] (
   \tfrac{1}{K} \sum_{k = 1}^K r_k )^K.
\end{equation}
Second, the cumulative progress $\sum_{k=1}^K r_k$
under lookahead action is bounded by 
(\Cref{lem:feedback-progress}):
\begin{equation} \label{eqn:algodemo-progress}
   \textstyle \sum_{k=1}^K r_k \leq \textstyle \sum_{k=1}^K r_{x^k} (P_k) - \tfrac{1}{4 L^2} \textstyle \sum_{k=1}^K \| \nabla r_{x^k} (P_k) \|_F^2.
\end{equation}
Third, the regret guarantee of online gradient descent bounds the cumulative feedback $\sum_{k=1}^K r_{x^k} (P_k)$ by (\Cref{lem:ogd-sqrtK-eta}):
\begin{equation} \label{eqn:algodemo-regret}
\textstyle \sum_{k=1}^K r_{x^k} (P_k) \leq \sum_{k=1}^K r_{x^k} (\hat{P}) + \tfrac{1}{2 \eta} \| P_1 - \hat{P}
   \|_F^2 + \tfrac{\eta}{2} \sum_{k=1}^K \| \nabla r_{x^k}
   (P_k) \|_F^2.
\end{equation}
Putting \eqref{eqn:algodemo-progress} and \eqref{eqn:algodemo-regret} together, the cumulative progress is bounded by
\begin{equation}\label{eqn:algodemo-together}
   \textstyle \sum_{k=1}^K r_k
   \leq  \sum_{k=1}^K r_{x^k} (\hat{P}) + \tfrac{1}{2 \eta} \| P_1 - \hat{P} \|_F^2 + (\tfrac{\eta}{2}
   - \tfrac{1}{4 L^2}) \textstyle \sum_{k=1}^K \| \nabla r_{x^k} (P_k) \|_F^2. 
\end{equation}
When $\eta \leq \tfrac{1}{2 L^2}$, lookahead action helps
the scheduler suppress the $\tfrac{\eta}{2} \| \nabla r_{x^k} (P_k) \|_F^2$ error term from online gradient descent. 
In particular, stepsize $\eta = \tfrac{1}{2 L^2}$ simplifies \eqref{eqn:algodemo-together} to
\begin{equation} \label{eqn:demo-ratio-constant-regret}
	\textstyle \sum_{k = 1}^K r_k \leq \sum_{k = 1}^K r_{x^k} (\hat{P}) + L^2 \| P_1 -
   \hat{P} \|_F^2.
\end{equation}
The global convergence guarantee follows immediately by plugging \eqref{eqn:demo-ratio-constant-regret} into the reduction \eqref{eqn:algodemo-reduction}, summarized in \Cref{thm:algodemo-global-conv} below.
The collaboration between landscape and scheduler allows {\osgmhandsonrx} to converge as if the online algorithm incurred only constant regret. 

\begin{figure}[h]
\centering
\includegraphics[scale=0.5]{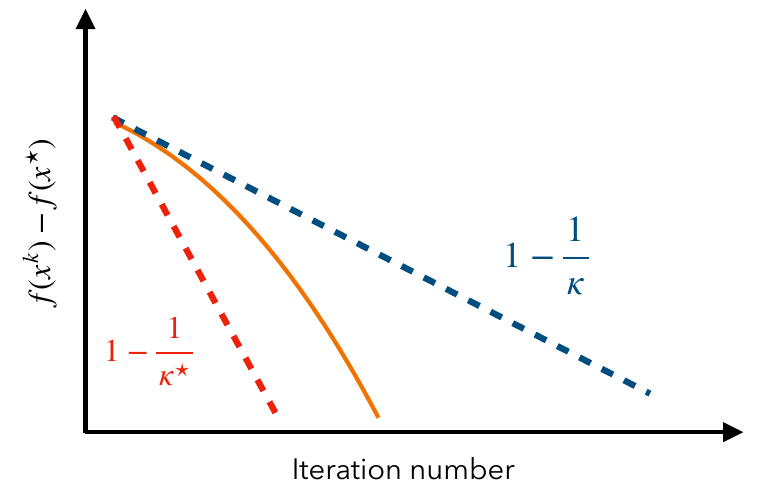}
\caption{Theoretical performance of {\osgmhandsonrx}. The linear convergence rate initially matches the $1 - \frac{1}{\kappa}$ rate of vanilla gradient descent and accelerates to (at least) $1 - \frac{1}{\kappa^\star}$. \label{fig:demo}}
\end{figure}

\begin{thm}[Global convergence] \label{thm:algodemo-global-conv} 
Let $f$ be $L$-smooth and $\mu$-strongly convex.
For any benchmark stepsize $\hat{P} \in \mathbb{R}^{n \times n}$, {\osgmhandsonrx} (\Cref{alg:demo}) with $\eta = 1 / (2 L^2)$ satisfies
\begin{equation} \label{eqn:handsonrx-tjry-conv}
   \textstyle f (x^{K + 1}) - f^{\star} \leq [f (x^1) - f^{\star}] (
   \tfrac{1}{K} \sum_{k = 1}^K r_{x^k} (\hat{P}) + \tfrac{L^2 \| P_1 -
   \hat{P} \|_F^2}{K} )^K.
\end{equation}
In particular, if {\osgmhandsonrx} is initialized with $P_1 = \frac{1}{L} I$, then 
\begin{equation} \label{eqn:algodemo-global-conv-rate}
	  f (x^{K + 1}) - f^{\star} \leq [f (x^1) - f^{\star}] \min\{ ( 1 -
   \tfrac{1}{\kappa} )^K,  ( 1 -
   \tfrac{1}{\kappa^{\star}} + \tfrac{L^2}{K} \| \tfrac{1}{L} I -
   P_r^{\star} \|_F^2 )^K\},
\end{equation}
where $\kappa$ and $\kappa^\star$ are the condition number and optimal condition number (\Cref{def:opt-preconditioner}) of $f$.
\end{thm}
The convergence guarantee in \eqref{eqn:handsonrx-tjry-conv} is powerful since it holds for any benchmark stepsize $\hat{P}$. In particular, one can choose the benchmark stepsize $\hat{P}$ that achieves the best average feedback along the algorithm trajectory. 
Moreover, if {\osgmhandsonrx} is initialized with $P_1 = \frac{1}{L} I$, according to \eqref{eqn:algodemo-global-conv-rate}, {\osgmhandsonrx} can automatically adapt to the best convergence rate among $1 - \tfrac{1}{\kappa}$ and $1 - \tfrac{1}{\kappa^{\star}} + \tfrac{L^2}{K} \| \tfrac{1}{L} I - P_r^{\star} \|_F^2$. The former is the rate of vanilla gradient descent, and the latter is asymptotically the rate of gradient descent with the optimal stepsize. 
When $K$ is small, the first rate $1 - \tfrac{1}{\kappa}$ is smaller; when $K$ is large, the second is smaller. 
Relation \eqref{eqn:algodemo-global-conv-rate} guarantees the best-of-both-worlds: {\osgmhandsonrx} at least matches vanilla gradient descent, and asymptotically converges at least as fast as gradient descent with the best possible stepsize.
\Cref{fig:demo} illustrates the expected convergence behavior.
An asymptotic $\mathcal{O} (\kappa^{\star} \log (1 / \varepsilon))$ complexity
can be obtained from \Cref{thm:algodemo-global-conv}. Alternatively, we can obtain an explicit complexity bound
through the lens of a potential function.
\paragraph{Potential reduction.}
Let $\hat{P}$ be a benchmark stepsize such that $r_x (\hat{P}) \leq 1 - \tfrac{1}{\kappa_{\hat{P}}} <
1$ for all $x$. For example, it suffices to take the optimal stepsize $\hat{P} = P_r^\star$ with its associated optimal condition number $\kappa_{\hat{P}} = \kappa^\star$. In each iteration of {\osgmhandsonrx}, either 
\begin{enumerate}[leftmargin=15pt, itemsep=0pt, label=\roman*)]
  \item the feedback $r_{x^k} (P_k)$ is smaller (better) than $r_{x^k} (\hat{P}) < 1$ and the suboptimality contracts:
  \[ f (x^{k + 1}) - f^{\star} = r_k [f (x^k) - f^{\star}] < f (x^k) - f^{\star}, \]
  \item or the feedback $r_{x^k} (P_k)$ is larger (worse) than $r_{x^k} (\hat{P})$ but $\| P_{k + 1} - \hat{P} \|_F^2$ shrinks:
  \[ \| P_{k + 1} - \hat{P} \|_F^2 \leq \| P_k - \hat{P} \|_F^2 - 2 \eta
     [r_{x^k} (P_k) - r_{x^k} (\hat{P})] + \eta^2 \| \nabla r_{x^k} (P_k)
     \|_F^2 . \]
\end{enumerate}
This observation motivates the definition of a joint potential in $x$ and $P$,
parametrized by the benchmark $\hat P$:
\[ \varphi (x, P) \assign \rho \log (f (x) - f^{\star}) + \| P - \hat{P}
   \|^2_F , \quad \rho > 0. \]
This potential combines the primal suboptimality $f (x) - f^{\star}$ and the distance to the benchmark stepsize $\| P - \hat{P} \|^2_F$. {\osgmhandsonrx} decreases this potential by at least a constant at every iteration.
\begin{thm}[Potential reduction]\label{thm:potential-handson-rx}
Let $f$ be $L$-smooth and $\mu$-strongly convex. 
At every iteration, {\osgmhandsonrx} with $\eta = 1 / (2 L^2)$ decreases the potential $\varphi (x, P)$ with $\rho = 1 / L^2$ by
\[ \varphi (x^{k + 1}, P_{k + 1}) \leq \varphi (x^k, P_k) - \tfrac{1}{\kappa_{\hat{P}} L^2}. \]
In particular, {\osgmhandsonrx} finds an $\varepsilon$-optimal solution within complexity
\[ K_{\varepsilon} \assign \Big\lceil \min_{\hat{P}} \big\{ \kappa_{\hat{P}} L^2
\| P_1 - \hat{P} \|^2_F + \kappa_{\hat{P}} \log \big( \tfrac{f(x^1) - f^{\star}}{\varepsilon}
\big) \big\} \Big\rceil. \]
\end{thm}

\input{proofs/proof-potential-handson-rx.tex}

The simple algorithm {\osgmhandsonrx} (\Cref{alg:demo}) delivers strong worst-case convergence guarantees without additional assumptions beyond smoothness and strong convexity. In the next subsection, we show that {\osgmhandsonrx} also achieves strong problem-dependent convergence as well as local convergence guarantees.

\subsubsection{Local convergence}

Online gradient descent replaces stale information with new updates, which allows {\osgm} to adapt to the local landscape. 
This subsection highlights three consequences of this adaptivity:
 1) trajectory-level adaptivity to the local landscape 2) local superlinear convergence 3) asymptotic optimality compared to any fixed stepsize.

\paragraph{Trajectory-level adaptivity.}
We combine \Cref{thm:reduction-ratio} and the dynamic regret guarantee in \Cref{lem:dynamic-sqrtK} to show local adaptivity to the trajectory.

\begin{thm}[Local adaptivity] \label{thm:local-trajectory-conv-handson-rx}
Let $f$ be $L$-smooth and $\mu$-strongly convex. 
For any benchmark sequence of stepsizes $\{ \hat{P}_k \}$, {\osgmhandsonrx} with $\eta = 1/(2L^2)$ satisfies
\[\textstyle f (x^{K + 1}) - f^{\star} \leq [f (x^1) - f^{\star}] (
   \tfrac{1}{K} \sum_{k = 1}^K r_{x^k} (\hat{P}_k) + \tfrac{L^2}{K} [\| \hat{P}_{K} - P_1 \|^2_F + 2 \displaystyle \max_{k \leq K} \{\| P_{k} - P_1 \|_F\} \mathsf{PL} (\{
      \hat{P}_k \})])^K, \]
where $\mathsf{PL} (\{ \hat{P}_k \}) \assign \sum_{k = 1}^{K-1} \| \hat{P}_k - \hat{P}_{k + 1} \|_F$ is the path length defined in \eqref{eqn:pathlength}.
\end{thm}
The convergence rate of {\osgmhandsonrx} competes with the average contraction ratio  $\tfrac{1}{K} \sum_{k = 1}^K r_{x^k} (\hat{P}_k)$ achieved by any benchmark sequence of stepsizes $\{ \hat{P}_k \}$, incurring an additional cost from the path length.

\input{proofs/proof-local-trajectory-conv.tex}

\paragraph{Local superlinear convergence.}
When $f$ is twice continuously differentiable, $f$ behaves like a quadratic function around $x^{\star}$: $f (x) = f^{\star} + \tfrac{1}{2} \langle x - x^{\star}, \nabla^2 f (x^{\star}) (x - x^{\star}) \rangle +\mathcal{O} (\| x - x^{\star} \|^3)$, for which the fixed stepsize $[\nabla^2 f (x^{\star})]^{- 1}$ enjoys (local) quadratic convergence.
Since {\osgmhandsonrx} can compete with any fixed stepsize, including $[\nabla^2 f (x^{\star})]^{- 1}$, {\osgmhandsonrx} converges superlinearly.  First, we formalize the intuition that $[\nabla^2 f (x^{\star})]^{- 1}$ is locally a good fixed stepsize by bounding its ratio feedback.

\begin{lem} \label{lem:local-opthessian-rx}
Suppose $f$ is $L$-smooth $\mu$-strongly convex and has $H$-Lipschitz Hessian. Then the ratio feedback of Hessian inverse at $x^\star$ is bounded by
$r_x ([\nabla^2 f (x^{\star})]^{- 1}) \leq \tfrac{H^2 \kappa}{4 \mu^2} \| x - x^{\star} \|^2$ for all $x \nin \Xcal^\star$.
\end{lem}

\input{proofs/proof-local-opthessian-rx.tex}

When $f$ is a quadratic, $H=0$ and the Hessian inverse drives any point $x$ to the optimal solution $x^{\star}$ in one step and thus achieves zero ratio feedback.
The superlinear convergence of {\osgmhandsonrx} follows from \Cref{lem:local-opthessian-rx} and \Cref{thm:algodemo-global-conv}, along with the observation that $\|x^k - x^\star\|\rightarrow 0$ geometrically.

\begin{thm}[Superlinear convergence] \label{thm:local-superlinear-conv-rx}
Suppose $f$ is $L$-smooth and $\mu$-strongly convex and has $H$-Lipschitz Hessian. Then {\osgmhandsonrx} (\Cref{alg:demo}) with $\eta = 1/(2L^2)$ and $P_1 = \frac{1}{L} I $ satisfies
\[ f (x^{K + 1}) - f^{\star} \leq [f (x^1) - f^{\star}] (
   \tfrac{C}{K} )^K, \]
where $C \assign \tfrac{H^2 \kappa^2}{2 \mu^3} [f (x^1) - f^{\star}] +
L^2 \| \frac{1}{L} I - [\nabla^2 f (x^{\star})]^{- 1} \|_F^2$.
\end{thm}

This result shows superlinear convergence of {\osgmhandsonrx}. The convergence appears linear (by \Cref{thm:algodemo-global-conv}) until $K \sim C$ and superlinear convergence becomes apparent when $K$ is large. 
The convergence behavior of {\osgmhandsonrx} in practice is the best of \Cref{thm:algodemo-global-conv} and \Cref{thm:local-superlinear-conv-rx}.

\input{proofs/proof-local-superlinear-conv-rx.tex}

\paragraph{Negative regret.} Relation \eqref{eqn:demo-ratio-constant-regret} suggests that the progress of {\osgmhandsonrx} is no worse than any benchmark stepsize $\hat{P}$ up to a fixed additive constant. Actually, we can further show that progress of {\osgmhandsonrx} is strictly better than any benchmark stepsize $\hat{P}$
unless there is a perfect stepsize that drives $x$ to $x^{\star}$ in one step (e.g., $\hat P = [\nabla^2 f (x^{\star})]^{- 1}$). 
\begin{thm}[Negative regret]  \label{thm:local-negative-regret-rx}
Fix an arbitrary benchmark stepsize $\hat{P}$. For each $K \geq 1$, {\osgmhandsonrx} with $\eta \in ( 0, \tfrac{1}{4 L^2} ]$ satisfies exactly one of the cases below:
\begin{itemize}[leftmargin=10pt,itemsep=0pt,]
   \item the average progress is smaller than that achieved by  $\hat{P}$: $\tfrac{1}{K} \textstyle \sum_{k = 1}^K r_k \leq \frac{1}{K} \textstyle \sum_{k = 1}^K r_{x^k} (\hat{P})$, 
   
   \item the suboptimality satisfies a superlinear convergence bound: $f (x^{K + 1}) - f (x^1) \leq \tfrac{1}{2 \mu} \| \nabla f (x^1) \|^2 ( \tfrac{\kappa^2 \| P_1 - \hat{P} \|_F^2}{\eta K} )^K$.
\end{itemize}
\end{thm}

In other words, {\osgmhandsonrx} outperforms the linear convergence rate achievable by any fixed stepsize. The only exception is when a fixed stepsize achieves a convergence rate of $0$, the case of superlinear convergence.

\paragraph{Knowledge of $f^{\star}$.} One downside of {\osgmrx} is the requirement of optimal value $f^{\star}$. We can relax this requirement by running an outer loop to search for $f^{\star}$ and obtain an $\Ocal(\kappa^\star \log^2 (1/\varepsilon))$ complexity result. See the conference version \cite{gao2024gradient} for more details.

\subsection{\osgmhandsonhx} \label{sec:monotone-hx}

In this section, we assume $f$ is $L$-smooth, optionally $\mu$-strongly convex and instantiate {\osgm} with
\[ \ell_x (P) \assign h_x (P), ~ \text{Monotone Lookahead landscape:~} f(x^{k + 1}) \leq \min \{ f(x^{k + 1 / 2}
   - \tfrac{1}{L} \nabla f (x^{k + 1 / 2})), f(x^k) \}, ~ \mathcal{A} \assign {\ogd} . \]
The algorithm is called {\osgmhandsonhx} (\Cref{alg:osgmhandsonhx}). The candidate stepsize set is $\Pcal = \Rbb^{n \times n}$. 
\begin{algorithm}[H]
{\textbf{input:}  Initial point $x^1$, initial stepsize $P_1 \in \Pcal = \Rbb^{n\times n}$, online gradient stepsize $\eta_k \equiv \eta > 0$}\\
\For{$k = 1, 2, \dots$}{
      $x^{k + 1/2} = x^{k} - P_k \nabla f(x^k)$\\
      Choose $x^{k+1}$ that satisfies $f(x^{k+1}) \leq \min \{ f(x^{k + 1 / 2}
   - \tfrac{1}{L} \nabla f (x^{k + 1 / 2})), f(x^k) \}$ \\
      $P_{k+1} = P_k - \eta \nabla h_{x^k}(P_k)$ \label{alg:osgmhandsonhx-line}
}
\caption{{\osgmhandsonhx}\label{alg:osgmhandsonhx}}
\end{algorithm}
Recall that the gradient of hypergradient feedback in \Cref{alg:osgmhandsonhx-line} of \Cref{alg:osgmhandsonhx} takes the form
\begin{equation*}
   \nabla h_{x^k}(P_k) = - \tfrac{\nabla f(x^{k + 1/2})\nabla f(x^k)^\top}{\|\nabla f(x^k) \|^2}.
\end{equation*}
The rest of the section analyzes the convergence behavior of {\osgmhandsonhx}, including global convergence and local superlinear convergence.

\subsubsection{Global convergence.} 
{\osgmhandsonhx} enjoys similar convergence guarantees to {\osgmhandsonrx}.

\begin{thm}[Global convergence] \label{thm:global-conv-handson-hx} 
Let $f$ be $L$-smooth and ($\mu$-strongly) convex. For any benchmark stepsize $\hat{P} \in \mathbb{R}^{n \times n}$, {\osgmhandsonhx} with $\eta = 1 / L$ satisfies
\begin{align}
 f (x^{K + 1}) - f^{\star} \leq{} & \min \{ \tfrac{\Delta^2}{K \max
   \{ \frac{1}{K} \sum^K_{k = 1} - h_{x^k} (\hat{P} ) - \frac{L}{2K}
   \| P_1 - \hat{P} \|_F^2, 0 \}}, f (x^1) - f^{\star} \},
   \tag{convex}\\
 f (x^{K + 1}) - f^{\star} \leq{} & [f (x^1) - f^{\star}] ( 1 -
   2 \mu \max \{ \tfrac{1}{K} \textstyle\sum^K_{k = 1} - h_{x^k} (\hat{P} ) -
   \tfrac{L}{2K} \| P_1 - \hat{P} \|_F^2 , 0 \} )^K,\tag{$\mu$-strongly convex}
\end{align}
where the constant $\Delta$ is defined in \Cref{thm:reduction-hypergrad}. 
In particular, if {\osgmhandsonhx} is initialized with $P_1 = \tfrac{1}{L} I$, then
\begin{align}
   f (x^{K + 1}) - f^{\star} \leq{} & \min \{ \tfrac{2 L \Delta^2}{K}, f (x^1) - f^{\star} \}, \tag{convex}\\
   f (x^{K + 1}) - f^{\star} \leq{} & [f (x^1) - f^{\star}] ( 1 - \tfrac{1}{\kappa} )^K . \tag{$\mu$-strongly convex}
\end{align}
\end{thm}

\input{proofs/proof-global-conv-handson-hx.tex}

According to \Cref{thm:global-conv-handson-hx}, {\osgmhandsonhx} converges no slower than vanilla gradient descent.
More importantly, whenever there exists a
$\hat{P}$ such that the average hypergradient feedback along the trajectory satisfies $\tfrac{1}{K} \sum^K_{k = 1} - h_{x^k} (\hat{P} ) \gg \tfrac{1}{2 L}$, {\osgmhandsonhx} converges faster. 

\paragraph{Potential reduction.} {\osgmhandsonhx} also admits a potential function analysis. Define\begin{align}
   \omega (x, P) &\coloneqq{}  \tfrac{\rho_{\omega}}{f (x) - f^{\star}} +
  \| P - \tfrac{1}{L} I \|_F^2, \tag{convex}\\
  \varphi (x, P) &\coloneqq{}  \rho_{\varphi} \log (f (x) - f^{\star}) +  \| P - \tfrac{1}{L} I \|_F^2. \tag{$\mu$-strongly convex}
\end{align}
{\osgmhandsonhx} decreases these potentials by at least a constant at every iteration.

\begin{thm}[Potential reduction] \label{thm:potential-handson-hx}
Let $f$ be $L$-smooth and ($\mu$-strongly) convex. 
At every iteration, {\osgmhandsonhx} with $\eta = 1 / L$ decreases the potential $\varphi (x, P)$ with $\rho_{\varphi} = 1 / (L \mu)$ and potential $\omega (x, P)$ with $\rho_{\omega} = 2
\Delta^2 / L$  by
\begin{align}
   &\omega (x^{k + 1}, P_{k + 1}) - \omega (x^k, P_k) \leq - \tfrac{1}{L^2}, \tag{convex}\\
   &\varphi (x^{k + 1}, P_{k + 1}) - \varphi (x^k, P_k) \leq - \tfrac{1}{L^2}. \tag{$\mu$-strongly convex}
\end{align}
In particular, {\osgmhandsonhx} with $P_1 = \tfrac{1}{L} I$ finds an $\varepsilon$-optimal solution within
 complexity
\begin{equation*}
   K_{\varepsilon} \assign \big\lceil \tfrac{2 L
   \Delta^2}{\varepsilon} \big\rceil \quad \text{for convex $f$}; \quad \text{and} \quad 
   K_{\varepsilon} \assign \big\lceil \kappa \log ( \tfrac{1}{\varepsilon} ) \big\rceil \quad \text{for $\mu$-strongly convex $f$}.
\end{equation*}

\end{thm}

\subsubsection{Local convergence}
The local convergence guarantees of {\osgmhandsonhx} are similar to {\osgmhandsonrx}.
\begin{thm}[Local adaptivity] \label{thm:local-trajectory-conv-handsonhx}
Let $f$ be $L$-smooth and ($\mu$-strongly) convex. 
For any benchmark sequence of stepsizes $\{ \hat{P}_k \}$, 
{\osgmhandsonhx} with stepsize $\eta = 1 / L$ satisfies
  \begin{align}
    f (x^{K + 1}) - f^{\star} \leq{} & \min \{ \tfrac{\Delta^2}{K \max
    \{ \frac{1}{K} \sum_{k = 1}^K -h_{x^k} (\hat{P}_k) - \rho_K(\{\hat{P}_k\}), 0 \}}, f (x^1) - f^{\star} \},
    \tag{convex}\\
    f (x^{K + 1}) - f^{\star} \leq{} & \textstyle [f (x^1) - f^{\star}] ( 1 -
    2 \mu \max \{ \frac{1}{K} \sum_{k = 1}^K -h_{x^k} (\hat{P}_k) - \rho_K(\{\hat{P}_k\}),0 \} )^K . \tag{$\mu$-strongly convex}
  \end{align}
where $\rho_K(\{\hat{P}_k\}) \assign \tfrac{L}{2K} [ \|\hat{P}_{K} - P_1\|_F^2 + 2 \max_{k \leq K} \{\|P_{k} - P_1\|_F\} \mathsf{PL} (\{ \hat{P}_k \})]$ and $\mathsf{PL} (\{ \hat{P}_k \})$ is defined in \eqref{eqn:pathlength}.
\end{thm}

\begin{thm}[Superlinear convergence] \label{thm:local-superlinear-conv-handson-hx}
Suppose $f$ is $L$-smooth, $\mu$-strongly convex and has $H$-Lipschitz Hessian. Then {\osgmhandsonhx} with stepsize
$\eta = 1 / L$ and $P_1 = \frac{1}{L} I$ satisfies
\[ f (x^{K + 1}) - f^{\star} \leq [f(x^1) - f^{\star}]( \tfrac{C}{K} )^K, \]
where $C \assign \tfrac{H^2 \kappa^3}{2 \mu^3} [f(x^1) - f^{\star}] + L^2 \| \tfrac{1}{L} I - [\nabla^2 f(x^{\star})]^{- 1} \|_F^2$.
\end{thm}

\begin{thm}[Negative regret]  \label{thm:local-negative-regret-handson-hx}
Fix an arbitrary benchmark stepsize $\hat{P}$. For each $K \geq 1$, {\osgmhandsonhx} with $\eta \in ( 0, \tfrac{1}{2 L} ]$ satisfies exactly one of the cases below:
\begin{itemize}[leftmargin=10pt,itemsep=0pt]
   \item the average progress is smaller than achieved by $\hat{P}$: $\tfrac{1}{K} \textstyle \sum_{k = 1}^K h_k \leq \frac{1}{K} \textstyle \sum_{k = 1}^K h_{x^k} (\hat{P})$, 
   
   \item the suboptimality satisfies a superlinear convergence bound: $f (x^{K + 1}) - f (x^1) \leq \tfrac{1}{2 \mu} \| \nabla f (x^1) \|^2
   ( \tfrac{2 L \| P_1 - \hat{P} \|_F^2}{\eta K} )^K$.
\end{itemize}
\end{thm}

\subsection{Other instances of {\osgm}}
Variants of {\osgm} can be obtained by enumerating combinations of feedback and landscape actions, but noting that hypergradient feedback should be used with monotone (lookahead) action due to the reduction for $h_k$.
{\osgmhandsoffrx} introduced in \Cref{sec:intro} is the only variant that requires no environment action and is still guaranteed to converge on smooth strongly convex problems. 
{\osgmhandsoffhx} often demonstrates the strongest empirical performance among all variants.
The analyses for other variants resemble {\osgmhandsonrx} and {\osgmhandsonhx} presented in the previous two subsections, except for the choice of $\eta_k$ in online gradient descent. We defer the highly repetitive details to \Cref{sec:other-instances} and summarize the convergence rates of each variant in \Cref{table:global-conv}.
Notably, {\osgm} is the second family of first-order methods with superlinear convergence guarantees on smooth convex optimization problems, following the renowned quasi-Newton methods. The online learning argument offers a simple superlinear convergence analysis of {\osgm}.

\begin{table}[!h]
   \centering
   \caption{Global convergence rates of {\osgm} \label{table:global-conv}}
   \renewcommand{\arraystretch}{1}
   \resizebox{0.85\textwidth}{!}{
   \begin{tabular}{cccccc}
      \toprule
      Feedback & Convexity & Landscape action $\Mcal$ & Worst-case global convergence & Superlinear convergence \\
      \midrule
      \multirow{4}{*}{$r_x$} 
      & \multirow{4}{*}{Strongly convex} & Lookahead          & \multirow{2}{*}{$\min \{ ( 1 -
      \frac{1}{\kappa^{\star}} +\mathcal{O} ( \frac{1}{K} )
      )^K, ( 1 - \frac{1}{\kappa} )^K \}$} & \multirow{2}{*}{$\mathcal{O}
      (e^{- K \log K})$} \\
      & & Monotone Lookahead & & \\
      \cmidrule{3-5}
      & & Vanilla         & \multirow{2}{*}{$( 1 - \frac{1}{\kappa^{\star}}
      +\mathcal{O} ( \frac{1}{\sqrt{K}} ) )^K$} & \multirow{2}{*}{$\mathcal{O}
      ( e^{- \frac{1}{2} K \log K} )$} \\
      & & Monotone & & \\
      \midrule
      \multirow{4}{*}{$h_x$}
      & \multirow{2}{*}{Strongly convex} & Monotone Lookahead  & $( 1 - \frac{1}{\kappa} )^K$ & {$\mathcal{O} (e^{- K \log K})$} \\
      & & Monotone & $( 1 - \frac{1}{\kappa} +\mathcal{O} ( \frac{1}{\sqrt{K}} ) )^K$ & {$\mathcal{O} (e^{- \frac{1}{2} K \log K})$}\\
      \cmidrule{2-5}
      & \multirow{2}{*}{Convex} & Monotone Lookahead & $\frac{2 L \Delta^2}{K}$ & \multirow{2}{*}{---} \\
      & & Monotone & $\frac{2L\Delta^2}{K} ( \frac{1}{1 -\mathcal{O} ( \frac{1}{\sqrt{K}} )} )$ & \\
      \bottomrule
   \end{tabular}}
\end{table}

%% file: proofs/proof-potential-handson-rx.tex
At every iteration, the function value part of the potential
changes by
\begin{equation} \label{eqn:pf-potential-handson-rx-0}
  \rho \log (f (x^{k + 1}) - f^{\star}) - \rho \log (f (x^k) - f^{\star}) = \rho \log \tfrac{f (x^{k + 1}) - f^{\star}}{f (x^{k}) - f^{\star}} = \rho \log r_k.
\end{equation}
On the other hand, the distance part of the potential changes by
\begin{align}
  \| P_{k + 1} - \hat{P} \|_F^2 & \overleq{\eqref{eqn:ogd-Pk-ineq}}  \| P_k - \hat{P} \|_F^2 - 2 \eta
  [r_{x^k} (P_k) - r_{x^k} (\hat{P})] + \eta^2 \| \nabla r_{x^k} (P_k) \|_F^2
  \nonumber\\
 & =  \| P_k - \hat{P} \|_F^2 - 2 \eta [ r_{x^k} (P_k) - \tfrac{\eta}{2}
  \| \nabla r_{x^k} (P_k) \|_F^2 - r_{x^k} (\hat{P}) ] \nonumber\\
 & ={}  \| P_k - \hat{P} \|_F^2 - 2 \eta [ r_{x^k} (P_k) - \tfrac{1}{4 L^2}
  \| \nabla r_{x^k} (P_k) \|_F^2 - r_{x^k} (\hat{P}) ] \tag{by $\eta = \tfrac{1}{2 L^2}$}\\
 & \leq{}  \| P_k - \hat{P} \|_F^2 - 2 \eta [r_k - r_{x^k} (\hat{P})], \label{eqn:pf-potential-handson-rx-2}
\end{align}
where \eqref{eqn:pf-potential-handson-rx-2} uses $r_k \leq r_{x^k} (P_k) - \tfrac{1}{4 L^2} \| \nabla r_{x^k} (P_k) \|_F^2$ by \Cref{lem:feedback-progress} with lookahead landscape action. Combining \eqref{eqn:pf-potential-handson-rx-0} and \eqref{eqn:pf-potential-handson-rx-2} gives
\begin{align}
  & \varphi (x^{k + 1}, P_{k + 1}) - \varphi (x^k, P_k) \nonumber\\
  ={} & \rho \log (f (x^{k + 1}) - f^{\star}) - \rho \log (f (x^k) - f^{\star}) + \| P_{k + 1} - \hat{P} \|_F^2 - \| P_k - \hat{P} \|_F^2
  \nonumber\\
  \leq{} & \rho \log r_k - 2 \eta [r_k - r_{x^k} (\hat{P})] \tag{by \eqref{eqn:pf-potential-handson-rx-2}}\\
  ={} & \rho \log r_k - 2 \eta r_k + 2 \eta r_{x^k} (\hat{P}) \nonumber\\
  \leq{} & \rho \log r_k - 2 \eta r_k + 2 \eta ( 1 -
  \tfrac{1}{\kappa_{\hat{P}}} ). \tag{by $r_{x^k} (\hat{P}) \leq 1 - \tfrac{1}{\kappa_{\hat{P}}}$} \\
 ={} & \alpha(r_k) + 2 \eta ( 1 - \tfrac{1}{\kappa_{\hat{P}}} ), \nonumber
\end{align}
where the function $\alpha (x) \assign \rho \log x - 2 \eta x$ is maximized at $x = \frac{\rho}{2 \eta}$ and $\rho \log r_k - 2 \eta r_k \leq \rho \log \tfrac{\rho}{2 \eta} - \rho$, which implies
\[ \varphi (x^{k + 1}, P_{k + 1}) - \varphi (x^k, P_k) \leq \rho \log
   \tfrac{\rho}{2 \eta} - \rho + 2 \eta ( 1 - \tfrac{1}{\kappa_{\hat{P}}}
   ) \]
Taking $\rho = 1 / L^2 = 2 \eta$, we get the desired reduction in the potential function:
\[ \varphi (x^{k + 1}, P_{k + 1}) - \varphi (x^k, P_k) \leq -
   \tfrac{1}{\kappa_{\hat{P}} L^2}. \]
To obtain the iteration complexity, note that
\begin{align}
  \tfrac{1}{L^2} \log (f (x^{K + 1}) - f^{\star}) \leq{} & \varphi (x^{K + 1}, P_{K +
  1}) \nonumber\\
  \leq{} & \varphi (x^1, P_1) - \tfrac{1}{\kappa_{\hat{P}} L^2} K \nonumber\\
  ={} & \tfrac{1}{L^2} \log [f (x^1) - f^{\star}] + \| P_1 - P^{\star} \|_F^2 -
  \tfrac{1}{\kappa_{\hat{P}} L^2} K. \nonumber
\end{align}
Hence, {\osgmhandsonrx} achieves $f (x^{K + 1}) - f^{\star} \leq \varepsilon$ for $K \geq \kappa_{\hat{P}} L^2 \| P_1 - P^{\star} \|_F^2 + \kappa_{\hat{P}} \log (\tfrac{f (x^1) - f^{\star}}{\varepsilon} )$. 
This holds for arbitrary benchmark stepsize $\hat{P}$. Taking the minimum (infimum) over $\hat{P}$ completes the proof.

%% file: proofs/proof-local-opthessian-rx.tex
We start by bounding the numerator $f (x - [\nabla^2 f (x^{\star})]^{- 1} \nabla f (x)) - f^{\star}$ of ratio feedback $r_x([\nabla^2 f (x^{\star})]^{- 1})$. 
By $f (x) - f^{\star}
\leq \tfrac{L}{2} \| x - x^{\star} \|^2$ from $L$-smoothness and $[\nabla^2 f (x^{\star})]^{- 1}
\preceq \tfrac{1}{\mu} I$ from $\mu$-strong convexity, we upper bound the numerator as follows:
\begin{align}
  f (x - [\nabla^2 f (x^{\star})]^{- 1} \nabla f (x)) - f^{\star}
  \leq{} & \tfrac{L}{2} \| x - [\nabla^2 f (x^{\star})]^{- 1} \nabla f (x) -
  x^{\star} \|^2 \nonumber\\
  ={} & \tfrac{L}{2} \| [\nabla^2 f (x^{\star})]^{- 1} [\nabla^2 f (x^{\star})
  (x - x^{\star}) - (\nabla f (x) - \nabla f (x^{\star}))] \|^2 \nonumber\\
  \leq{} & \tfrac{L}{2 \mu^2} \| \nabla^2 f (x^{\star}) (x - x^{\star}) -
  (\nabla f (x) - \nabla f (x^{\star})) \|^2. \label{eqn:local-opthessian-rx-0}
\end{align}

Using $\nabla f (x) - \nabla f (x^{\star}) = \textstyle\int_0^1 \nabla^2 f (x^{\star} +
t (x - x^{\star})) (x - x^{\star}) \mathd t$, the norm in \eqref{eqn:local-opthessian-rx-0} can be further bounded by
\begin{align}
  & \| \nabla^2 f (x^{\star}) (x - x^{\star}) - (\nabla f (x) - \nabla f
  (x^{\star})) \| \nonumber\\
  ={} & \| \nabla^2 f (x^{\star}) (x - x^{\star}) - \textstyle\int_0^1 \nabla^2 f
  (x^{\star} + t (x - x^{\star})) (x - x^{\star}) \mathd t \|
  \nonumber\\
  ={} & \| \textstyle\int_0^1 [\nabla^2 f (x^{\star}) - \nabla^2 f (x^{\star} + t (x
  - x^{\star}))] (x - x^{\star}) \mathd t \| \nonumber\\
  \leq{} & \textstyle\int_0^1 t H \| x - x^{\star} \|^2 \mathd t ={} \tfrac{H}{2} \| x -
  x^{\star} \|^2, \tag{by $H$-Lipschitz Hessian}
\end{align}
and thus \eqref{eqn:local-opthessian-rx-0} becomes
\begin{equation} \label{eqn:local-opthessian-rx-1}
	f (x - [\nabla^2 f (x^{\star})]^{- 1} \nabla f (x)) - f^{\star} \leq
\tfrac{L}{2 \mu^2} ( \tfrac{H}{2} \| x - x^{\star} \|^2 )^2 =
\tfrac{L H^2}{8 \mu^2} \| x - x^{\star} \|^4.
\end{equation}

Dividing both sides by $f (x) -
f^{\star}$ and using $f (x) - f^{\star} \geq \tfrac{\mu}{2} \| x
- x^{\star} \|^2$, we conclude that
\begin{align}
  r_x ([\nabla^2 f (x^{\star})]^{- 1}) ={} & \tfrac{f (x - [\nabla^2 f
  (x^{\star})]^{- 1} \nabla f (x)) - f^{\star}}{f (x) - f^{\star}}
  \nonumber\\
  \leq{} & \tfrac{L H^2}{8 \mu^2} \tfrac{\| x - x^{\star} \|^4}{f (x) -  f^{\star}} \leq{} \tfrac{L H^2}{4 \mu^3} \| x - x^{\star} \|^2 = \tfrac{H^2
  \kappa}{4 \mu^2} \| x - x^{\star} \|^2. \nonumber
\end{align}

%% file: proofs/proof-local-superlinear-conv-rx.tex
Plugging $r_x ([\nabla^2 f (x^{\star})]^{- 1}) \leq \tfrac{H^2 \kappa}{4 \mu^2}
\| x - x^{\star} \|^2$ from \Cref{lem:local-opthessian-rx}
into \Cref{thm:algodemo-global-conv}, we get
\begin{equation} \label{eqn:local-regret-rx}
\textstyle f (x^{K + 1}) - f^{\star} \leq [f (x^1) - f^{\star}] (
     \tfrac{H^2 \kappa}{4 \mu^2} \tfrac{1}{K} \sum_{k = 1}^K \| x^k -
     x^{\star} \|^2 + \tfrac{L^2}{K} \| P_1 - [\nabla^2 f (x^{\star})]^{- 1}
     \|_F^2 )^K .
\end{equation}
Recall from \Cref{thm:algodemo-global-conv} that {\osgmhandsonrx} satisfies $f (x^k) - f^{\star} \leq ( 1 - \tfrac{1}{\kappa} )^{k - 1} [f (x^1) - f^{\star}]$ and thus together with $\tfrac{\mu}{2} \| x^k - x^{\star} \|^2 \leq f (x^k) - f^{\star}$, we have
\begin{align}
  \textstyle \sum_{k = 1}^K \| x^k - x^{\star} \|^2 
  \leq{} & \textstyle \sum_{k = 1}^K \tfrac{2}{\mu} [f (x^k) - f^{\star}] \nonumber\\
  \leq{} & \tfrac{2}{\mu} [f (x^1) - f^{\star}] \textstyle \sum_{k = 1}^K ( 1 - \tfrac{1}{\kappa} )^{k - 1}
    \nonumber\\
    \leq{} & \tfrac{2}{\mu} [f (x^1) - f^{\star}] \kappa, \label{eqn:pf-superlinear-rx-1}
\end{align}
where the last inequality uses the relation $\sum_{k = 1}^K \gamma^{k-1} \leq \frac{1}{1 - \gamma}$ for $\gamma \in [0, 1)$.
Plugging \eqref{eqn:pf-superlinear-rx-1} back into \eqref{eqn:local-regret-rx} gives
\begin{align}
  f (x^{K + 1}) - f^{\star} \leq{} & \textstyle[f (x^1) - f^{\star}] (
  \tfrac{H^2 \kappa}{4 \mu^2} \tfrac{1}{K} \sum_{k = 1}^K \| x^k - x^{\star}
  \|^2 + \tfrac{L^2}{K} \| P_1 - [\nabla^2 f (x^{\star})]^{- 1} \|_F^2
  )^K \nonumber\\
  \leq{} & [f (x^1) - f^{\star}] ( \tfrac{H^2 \kappa^2}{2 \mu^3}
  \tfrac{f (x^1) - f^{\star}}{K} + \tfrac{L^2}{K} \| P_1 - [\nabla^2 f
  (x^{\star})]^{- 1} \|_F^2 )^K \nonumber\\
  ={} & [f (x^1) - f^{\star}] ( \tfrac{C}{K} )^K,
\end{align}
where $C = \tfrac{H^2 \kappa^2}{2 \mu^3} [f (x^1) - f^{\star}] + L^2 \| \frac{1}{L} I - [\nabla^2 f (x^{\star})]^{- 1} \|_F^2$.

%% file: proofs/proof-global-conv-handson-hx.tex
Using \Cref{lem:feedback-progress} with lookahead action, we bound the cumulative progress by
\begin{equation} \label{eqn:pf-global-conv-handson-hx-0}
  \textstyle \sum_{k = 1}^K h_k \leq \textstyle \sum_{k = 1}^K \min \{ h_{x^k} (P_k) - \tfrac{1}{2 L} \| \nabla h_{x^k} (P_k) \|_F^2, 0 \} 
  \leq \textstyle \sum_{k = 1}^K h_{x^k} (P_k) - \tfrac{1}{2 L} \textstyle \sum_{k = 1}^K \| \nabla h_{x^k} (P_k) \|_F^2.
\end{equation}
By \Cref{lem:ogd-sqrtK-eta} and the choice of stepsize $\eta = \tfrac{1}{L}$, we further bound the right-hand side of \eqref{eqn:pf-global-conv-handson-hx-0}:
\begin{align}
  \textstyle \sum_{k = 1}^K h_k
  \leq{} & \textstyle \sum_{k = 1}^K h_{x^k} (P_k) - \tfrac{1}{2 L} \textstyle \sum_{k = 1}^K \|
  \nabla h_{x^k} (P_k) \|_F^2 \nonumber\\
  \leq{} & \textstyle \sum_{k = 1}^K h_{x^k} (\hat{P}) + \tfrac{1}{2 \eta} \| P_1 - \hat{P}
  \|_F^2 + ( \tfrac{\eta}{2} - \tfrac{1}{2 L} ) \textstyle \sum_{k = 1}^K \|
  \nabla h_{x^k} (P_k) \|_F^2 \label{eqn:pf-global-conv-handson-hx-0.5} \\
  \leq{} & \textstyle \sum_{k = 1}^K h_{x^k} (\hat{P}) + \tfrac{L}{2} \| P_1 - \hat{P} \|_F^2 \label{eqn:pf-global-conv-handson-hx-1}.
\end{align}
Together with the hypergradient reduction from \Cref{thm:reduction-hypergrad} and $h_k \leq 0$, we conclude 
\paragraph{Convex $f$.}
    \begin{align}
      f (x^{K + 1}) - f^{\star} \leq{} & \min \{ \tfrac{\Delta^2}{\sum_{k = 1}^K - h_k}, f (x^1) - f^{\star}
      \} \nonumber\\
      \leq{} & \min \big\{ \tfrac{\Delta^2}{\max \{ \sum_{k = 1}^K - h_{x^k}
      (\hat{P}) - \frac{L}{2} \| P_1 - \hat{P} \|_F^2, 0 \}}, f (x^1) - f^{\star} \big\}; \nonumber
    \end{align}
\paragraph{$\mu$-strongly convex $f$.}
    \begin{align}
      f (x^{K + 1}) - f^{\star} 
      \leq{} & [f (x^1) - f^{\star}] ( 1 +
      \tfrac{2 \mu}{K} \textstyle \sum^K_{k = 1} h_k )^K \nonumber\\
      \leq{} & [f (x^1) - f^{\star}] ( 1 - 2
      \mu \max \{ \tfrac{1}{K} \textstyle \sum_{k = 1}^K - h_{x^k} (\hat{P}) -
      \tfrac{L}{2 K} \| P_1 - \hat{P} \|_F^2, 0 \} )^K . \nonumber
    \end{align}

%% file: sec_practical.tex
\section{Practical algorithm design with {\osgm}} \label{sec:practical}

This section discusses the implementation of {\osgm} for practical performance.

\subsection{Choice of candidate stepsize}

The candidate set $\mathcal{P}$ affects the convergence of {\osgm} through its expressive power and regret and determines the computational efficiency of {\osgm}.

\paragraph{Expressiveness of $\mathcal{P}$ and superlinear convergence.}
Richer candidate sets $\mathcal{P}$ leads to smaller minimum average feedback $\min_{\hat{P} \in \mathcal{P}} \tfrac{1}{K} \sum_{k = 1}^K \ell_{x^k} (\hat{P})$, indicating faster convergence rate {\osgm} can achieve asymptotically. 
In particular, inverse Hessian at optimal solution $[\nabla^2 f(x^{\star})]^{- 1}$ guarantees local superlinear convergence. If $[\nabla^2 f(x^{\star})]^{- 1}
\in \mathcal{P}$, then it is a legitimate benchmark $\hat{P}$ for {\osgm} to compete against, indicating local superlinear convergence of {\osgm} too.

\paragraph{Online regret.}
Richer candidate sets however may lead to larger dimension-dependent constant $\|P_1 - \hat{P}\|^2_F$ and online gradient norm $\| \nabla \ell_x (P)
\|_F^2$ in the regret and increase the kick-in period.
A scalar or diagonal stepsize often outperforms a full matrix if we expect to perform only a few iterations $K$.\\

Three common choices for $\mathcal P$ in practice are listed in \Cref{table:pattern}.

\begin{table}[h]
  \centering
  \caption{Examples of stepsize patterns and corresponding online gradients.}
  \label{table:pattern}
  \resizebox{0.85\textwidth}{!}{
  \begin{tabular}{cccc}
  \toprule
  \multirow{2}{*}{Pattern} & \multirow{2}{*}{Candidate set $\mathcal{P}$} & \multicolumn{2}{c}{Online gradient} \\ 
  \cmidrule{3-4}
  & & Ratio $r_x$ & Hypergradient $h_x$ \\
  \midrule
  Scalar &
  $\{ \alpha I : \alpha \in \mathbb{R} \}$ &
  $r'_x (\alpha) = - \frac{\langle \nabla f (x - \alpha \nabla f (x)), \nabla f (x) \rangle}{f(x) - f^{\star}}$ &
  $h'_x (\alpha) = - \frac{\langle \nabla f (x - \alpha \nabla f (x)), \nabla f (x) \rangle}{\|\nabla f(x)\|^2}$ \\[5pt]
  Diagonal &
  $\{ \operatorname{Diag}(d) : d \in \mathbb{R}^n \}$ &
  $\nabla r_x (d) = - \frac{\nabla f (x - d \circ \nabla f (x)) \circ \nabla f (x)}{f(x) - f^{\star}}$ &
  $\nabla h_x (d) = - \frac{\nabla f (x - d \circ \nabla f (x)) \circ \nabla f (x)}{\|\nabla f(x)\|^2}$ \\[5pt]
  Full matrix &
  $\mathbb{R}^{n \times n}$ &
  $\nabla r_x (P) = - \frac{\nabla f (x - P \nabla f (x)) \nabla f (x)^{\top}}{f(x) - f^{\star}}$ &
  $\nabla h_x (P) = - \frac{\nabla f (x - P \nabla f (x)) \nabla f (x)^{\top}}{\|\nabla f(x)\|^2}$ \\[5pt]
  \bottomrule
  \end{tabular}
}
\end{table}

Additional constraints, such as positive-semidefiniteness and boundedness, can
be enforced in {\ogd} update $P_{k + 1} = \Pi_{\mathcal{P}} [P_k - \eta \nabla \ell_{x^k} (P_k)]$ by an explicit projection but at different computational costs.
For example, 
entrywise projection onto an interval simply truncates the entries to values in the interval,
while projection on the semidefinite cone requires a full eigendecomposition. 
Our theory for {\osgm} applies even when {\ogd} involves the projector $\Pi_{\mathcal{P}}[\cdot]$
due to the non-expansiveness of the projection operator:
\[ \| P_{k + 1} - \hat{P} \|_F^2 = \| \Pi_{\mathcal{P}} [P_k - \eta \nabla
   \ell_{x^k} (P_k)] - \hat{P} \|_F^2 \leq \| P_k - \eta \nabla \ell_{x^k} (P_k) - \hat{P}
   \|_F^2, \text{~~for all~~} \hat{P} \in \Pcal.\]

\paragraph{Memory and iteration cost.}
Storing and updating a full matrix costs $\mathcal{O}(n^{2})$ memory and time, whereas a diagonal or scalar parameterization is essentially free, with memory and computational costs that match the cost of a gradient step. 
Projecting onto the set of positive semidefinite matrices $\mathcal P = \mathbb S^n_+$ requires a full eigendecomposition, at a cost of $\Ocal(n^3)$, which is very rarely worthwhile.
Indeed, one of the advantages of {\osgm} over traditional preconditioning is that we directly parametrize the scaling of the gradient, and so need not ensure the learned scaling is positive semidefinite.

\paragraph{Recommended choice.}If a good guess of $[\nabla^2 f
(x^{\star})]^{- 1}$ is available, {\osgm} locates $[\nabla^2 f(x^{\star})]^{-
1}$ efficiently in the first few iterations, and a full matrix stepsize is
preferred. Otherwise, a diagonal stepsize is generally the most practical option.
We observe that projection onto positive scalings is sometimes helpful for diagonal stepsizes.

\subsection{Choice of feedback}
Hypergradient feedback $h_x$ generally has stronger empirical performance than ratio feedback $r_x$, even though the latter has desirable theoretical properties.
There are two possible reasons:
\begin{itemize}[leftmargin=10pt]
  \item Hypergradient feedback does not require the knowledge of optimal value $f^{\star}$.
  
  \item Hypergradient feedback has a better dependence on the smoothness constant: $h_{x}$ is $L$-smooth but $r_{x}$ is $2L^{2}$-smooth.
      Consequently, $r_{x}$ may incur larger regret and its theoretical advantage can only be observed for very large $K$.
\end{itemize}
An exception where ratio feedback is preferred occurs when
1) $\kappa^{\star} \leq 100$ and 
2) a good initial stepsize $P_1$ is available. 
In this case, {\osgmrx} sometimes outperforms {\osgmhx}. 

\subsection{Choice of landscape action}
Vanilla landscape action is less stable and not recommended. In most cases, a monotone landscape action suffices to yield satisfying performance.
A null step $x^{k+1} = \argmin_{x \in \{x^{k+1/2}, x^k\}} f(x)$ is particularly recommended since the number of gradient oracle calls at each iteration equals that of vanilla gradient descent. 

\subsection{Choice of online learning algorithm}
Other advanced online learning algorithms can be used in place of {\ogd}, and our convergence analysis still applies.
Advanced online learning algorithms often outperform {\ogd}, and we recommend using {\adagrad} \cite{duchi2011adaptive} and other advanced online learning algorithms from \cite{orabona2016coin,jacobsen2023unconstrained} in {\osgm}.
Most convergence guarantees in this paper hold for {\adagrad} variant of {\osgm}, though possibly in a weaker statement.

\paragraph{Online learning parameters.} 
Online learning algorithms often come with their own hyperparameters, such as the stepsize $\eta$ in {\ogd}. 
Our analysis assumes the knowledge of $L$ or at least an upper bound. 
One option to relax this assumption is to use parameter-free online algorithms such as \cite{orabona2016coin}, which typically require $\diam (\Pcal) < \infty$.
Another option is to apply back-tracking line-search to estimate $L$, relying on the fact that both $f$ and feedback $r_x (P)$/$h_x (P)$ are smooth. 
For example, given an estimate $L'$ of $L$, we can monitor two conditions below:
\begin{align}
  h_{x^k} (P_k - \tfrac{1}{L'} \nabla h_{x^k} (P_k)) - h_{x^k} (P_k) & \leq -
  \tfrac{1}{2 L'} \| \nabla h_{x^k} (P_k) \|_F^2, \nonumber\\
  f ( x^{k + 1 / 2} - \tfrac{1}{L'} \nabla f (x^{k + 1 / 2}) ) - f
  (x^{k + 1 / 2}) & \leq -\tfrac{1}{2 L'} \| \nabla f (x^{k + 1 / 2}) \|^2 .
  \nonumber
\end{align}
Whenever either condition fails, one can backtrack $L'$ by a certain fraction to find new estimates.

\subsection{Momentum-based optimization step}

The most practical variant of {\osgm} adds the heavy-ball momentum to the step of gradient descent:
\[x^{k+1} = x^{k} - P_k \nabla f(x^k) + \beta_k (x^k - x^{k-1})\]
and jointly learns the stepsize $P_k$ and momentum parameter $\beta_k$ by online learning.
The hypergradient feedback is modified based on the potential function for heavy-ball momentum \cite{danilova2020non}:
\[ \varphi(x,x^-) = f(x) - f^{\star} + \tfrac{\omega}{2} \|x - x^-\|^2.  \]
Heavy-ball variant of {\osgmhx} demonstrates promising performance in \Cref{sec:exp} and \cite{chu2025provable}.

%% file: sec_exp.tex
\section{Experiments} \label{sec:exp}

We have developed an efficient variant of {\osgmhx} for large-scale machine learning tasks. See Part II of this paper for more details.

%% file: sec_literature.tex
\section{Related work} \label{sec:literature}

{\osgm} is closely related to several lines of research, which we detail below.

\paragraph{Preconditioned gradient descent.}The update $x^{k + 1} = x^k - P_k
\nabla f (x^k)$ that scales the gradient by $P_K$ is known as preconditioned gradient descent in literature. Preconditioning has become a standard tool in
both optimization and numerical linear algebra and is widely used in modern
optimization algorithms
{\cite{lin2025pdcs,xiong2024role,lu2023cupdlp,lu2023cupdlpc,goulart2024clarabel,lin2021admm,deng2024enhanced,o2016conic,applegate2021practical,huang2024restarted}}.
Some recent research has tried to empirically quantify the effect of
preconditioning on linear systems {\cite{qu2024optimal,gao2023scalable}}.
However, identifying a good preconditioner with theoretical guarantees is
often challenging in practice. {\osgm} allows gradient-based methods
to converge as if it is preconditioned by the best possible preconditioner,
with the cost amortized along the iterations.

\paragraph{Adaptive first-order methods and stepsize
scheduling.}{\osgm} is closely related to the literature on stepsize
scheduling and the widely used adaptive gradient methods. Examples in this family include stepsize schedulers such as
{\cite{li2021second,wang2023convergence}} and adaptive methods such as
\texttt{AdaGrad} {\cite{duchi2011adaptive}}, \texttt{Adam}
{\cite{zhang2024adam,nesterov1983method}}, \texttt{RMSProp}
{\cite{hinton2012neural}} among many other popular variants
{\cite{gupta2018shampoo,reddi2019convergence,zhuang2020adabelief,disentangling2020}}.
Many of the adaptive gradient methods originate from online learning and have
provable regret guarantees, and they also relate to parameter-free online
learning algorithms (see {\cite{orabona2019modern}} for more details). Aside
from stepsize selection strategies motivated by online learning, there are
well-known stepsize selection strategies such as the Barzilai-Borwein (BB) step
{\cite{barzilai1988two}}, Polyak stepsize
{\cite{polyak1987introduction,hazan2019revisiting}} and notable recently
developed adaptive stepsizes {\cite{malitsky2020adaptive,malitsky2024adaptive}}. Another relevant line of research focuses on uniformly optimal first-order methods
{\cite{li2023simple,deng2024uniformly}}, which choose stepsize to adapt to
problem parameters, such as the smoothness constant. The idea of analyzing the multi-step convergence relates to recent work on stepsize hedging \cite{altschuler2024acceleration, altschuler2025acceleration}.

\paragraph{Quasi-Newton methods.}One notable feature of {\osgm} is its non-asymptotic superlinear convergence. See \Cref{table:superlin} for a comparison with recent results showcasing superlinear convergence rates of quasi-Newton methods.

\begin{table}[h]
  \centering
  \caption{Superlinear convergence rates of various methods \label{table:superlin}}
\resizebox{0.65\textwidth}{!}{
  \begin{tabular}{ccc}
    \toprule
    Algorithm & Reference & Superlinear convergence rate \\ \midrule
    Greedy quasi-Newton& \cite{rodomanov2021greedy}  & $\mathcal{O}(e^{-\frac{1}{2} K^{2}})$ \\
    Broyden family& \cite{rodomanov2022rates,rodomanov2021new}  & $\mathcal{O}(e^{-\frac{1}{2} K \log K})$ \\
   Online-learning guided quasi-Newton & \cite{jiang2023online}  & $\mathcal{O}(e^{-\frac{1}{2} K \log K})$ \\
    BFGS with line-search& \cite{jin2024non,jin2024nonArmijo}  & $\mathcal{O}(e^{-K \log K})$ \\
   {Monotone Lookahead  \osgm} & This paper  & $\mathcal{O}(e^{-K \log K})$ \\
   {Vanilla/Monotone \osgm} & This paper  & $\mathcal{O}(e^{-\frac{1}{2}K \log K})$ \\
    \bottomrule
  \end{tabular}
}
\end{table}

Our results identify a similarity between {\osgm} and quasi-Newton methods.
Both {\osgm} and {\hdm} learn the inverse Hessian operator $g \mapsto [\nabla^2 f(x^{\star})]^{-1} g$ as the algorithm progresses, but through different properties of the operator.
The quasi-Newton methods use the secant equation $x - y \approx [\nabla^2 f(x^{\star})]^{-1}(\nabla f(x)- \nabla f(y))$ for $x, y$ close to $x^\star$ and enforce this equation, replacing the inverse Hessian by $P_k$,
to guide learning \cite{jiang2023online, jiang2024online}.
In contrast, {\osgm} learns an optimal stepsize for the function. 
Since the function is locally quadratic, this optimal stepsize is the inverse Hessian.
{\osgm} uses the ratio/hypergradient feedback to measure the quality of the stepsize directly and can search for an optimal stepsize in a given closed convex set $\mathcal{P}$. Both approaches require a safeguard to prevent divergence in the warm-up phase, which is achieved by line-search in quasi-Newton and landscape actions in {\osgm}. In a word, both {\osgm} and quasi-Newton leverage complementary perspectives on $g \mapsto [\nabla^2 f(x^{\star})]^{-1} g$, so it is natural that they achieve similar convergence guarantees.

\paragraph{Hypergradient descent heuristic.}{\osgm} is closely
relevant to the hypergradient descent heuristic. Hypergradient descent method ({\hdm}) dates back to 1999 \cite{almeida1999parameter}, which was first proposed as a heuristic to accelerate stochastic gradient descent.  Similar updates were also explored in \cite{sutton1992adapting,schraudolph1999local,jacobs1988increased,mahmood2012tuning}, while those works employed slightly different algorithmic updates. 
Later, \cite{gunes2018online} rediscovered the {\hdm} and named it ``hypergradient descent''; \cite{gunes2018online} also extended {\hdm} to other first-order methods with extensive experimental validation of its practical efficacy.  Recent studies \cite{jie2022adaptive,chandra2022gradient,ozkara2024mada} further empirically enhanced {\hdm} for broader applicability, reporting promising numerical results. Despite these empirical successes, a rigorous theoretical understanding of {\hdm} has emerged only recently. \cite{rubio2017convergence} showed that a variant of {\hdm} converges on convex quadratic functions and established several analytic properties. \cite{powerball2023} showed the convergence of hypergradient descent with a ``powerball'' technique and variance-reduction on stochastic finite-sum optimization. However, these results do not explain the empirical performance of {\hdm}. Subsequently, \cite{kunstner2024searching} demonstrated that when using a diagonal preconditioner, hypergradient can be employed to generate cutting planes in the preconditioner space, achieving an $\Ocal(\sqrt{n}\kappa^\star \log (1/\varepsilon))$ complexity result on smooth strongly convex functions. Prior to \cite{kunstner2024searching}, a similar idea was exploited in \cite{MonteiroONealNemirovski2004}, where the authors used the ellipsoid method to update the preconditioner in the conjugate gradient method for solving ill-conditioned linear systems. The preliminary versions of this work \cite{gao2024gradient, chu2025provable} showed that {\hdm} can be viewed as online gradient descent applied to some surrogate loss function and that {\hdm} has strong trajectory-based convergence guarantees.

\paragraph{Connection between {\osgm} and {\hdm} in literature.}
The most typical version of hypergradient descent \cite{gunes2018online, rubio2017convergence} swaps the order of primal update and stepsize update in {\osgmhx}:
\begin{equation*}
  \renewcommand{\arraystretch}{1.3}
  \begin{array}{lcl}
    {\hdm} & \qquad 
    & \texttt{Lookahead } {\osgmhx}
    \\
    P_{k + 1} = P_k - \eta \nabla h_{x^k} (P_k) & 
    & x^{k+1/2} = x^k - P_k \nabla f (x^k) 
    \\
    x_{\hdm}^{k+1} = x^k - P_{k + 1} \nabla f (x^k) &
    & x_{{\osgm}}^{k+1} = x^{k+1/2} - \tfrac{1}{L} \nabla f(x^{k+1/2}) \\
    & & P_{k + 1} = P_k - \eta \nabla h_{x^k} (P_k) \\
  \end{array}
\end{equation*}

{\hdm} first updates the matrix stepsize $P_k$ using the feedback $h_{x^k}(P)$ and then makes a gradient step on $x^k$ using $P_{k+1}$ to arrive at the next iterate $x^{k + 1}$. The per-iteration progress of {\hdm} is therefore
\[ 
h_k = \tfrac{f(x^{k+1}_{\hdm}) - f(x^k)}{\| \nabla f(x^k)\|^2} 
= \tfrac{f(x^{k} - P_{k+1} \nabla f(x^k)) - f(x^k)}{\| \nabla f(x^k)\|^2}
= h_{x^k}(P_{k+1}),
\]
which is the hypergradient feedback evaluated at $P_{k+1}$ whose update already uses the feedback function $h_{x^k}(P)$.
This setting of {\hdm} can be modeled as {\textit{prescient}} online learning {\cite{wang2024no}}, in which $h_{x^k}(P)$ is revealed to the scheduler before making a decision. The received feedback $h_{x^k}(P_{k+1})$ is evaluated at the newly chosen stepsize.
The knowledge of the future allows the stepsize scheduler in {\hdm} to achieve a constant regret {\cite{wang2024no,orabona2019modern,rakhlin2009lecture}}, which is similar to the guarantee achieved by \texttt{Lookahead} {\osgm}. 
This is not a coincidence: when $\Pcal = \Rbb^{n \times n}$, {\hdm} reduces to \texttt{Lookahead} {\osgmhx} except the stepsize of additional gradient step is not $\tfrac{1}{L}$ but $\eta$:
\begin{align}
  x_{\hdm}^{k + 1} ={} & x^k - P_{k + 1} \nabla f (x^k) \nonumber\\
  ={} & x^k - (P_k - \eta \nabla h_{x^k} (P_k)) \nabla f (x^k) \tag{By update of $P_{k+1}$} \\
  ={} & x^k - ( P_k + \eta \tfrac{\nabla f (x^k - P_k \nabla f (x^k)) \nabla f
  (x^k)^{\top}}{\| \nabla f (x^k) \|^2} ) \nabla f (x^k) \tag{By definition of $\nabla h_x(P)$}\\
  ={} & x^k - P_k \nabla f (x^k) - \eta \nabla f ( x^k - P_k \nabla f (x^k)) \nonumber\\
  ={} & x^{k + 1 / 2} - \eta \nabla f (x^{k + 1 / 2}) . \tag{$x^{k+1/2} = x^k - P_k \nabla f(x^k)$}
\end{align}
In other words, {\hdm} uses the same stepsize $\eta$ for both additional gradient step and online gradient descent, while \texttt{Lookahead} {\osgmhx} decouples the two stepsizes. 
{\hdm} and \texttt{Lookahead} {\osgm} are equivalent if the stepsize of online gradient descent is set to $\eta = \tfrac{1}{L}$ and $\mathcal{P}=\mathbb{R}^{n \times n}$ (no projection in online gradient descent).
When $\eta < 1/L$, the lookahead landscape action makes more progress than the swapped update $x_{\hdm}^{k + 1} = x^k - P_{k + 1} \nabla f (x^k)$.

\paragraph{Reconcile with lower bound.} {\osgm} provides an problem-dependent $\Ocal{(\kappa^\star\log(1/\varepsilon))}$ acceleration result that can outperform the accelerated rate $\Ocal{(\sqrt{\kappa}\log(1/\varepsilon	))}$ when $\kappa^\star < \sqrt{\kappa}$ and $\varepsilon \rightarrow 0$. However, our result does not conflict with the known lower bound for smooth strongly convex optimization \cite{nesterov2013introductory}. Instead, the acceleration effect should be considered as the effect of implicit preconditioning.

%% file: conclusion.tex
\section{Conclusions}

Online scaled gradient methods ({\osgm}) offer a new mechanism for accelerating gradient-based methods and constitute a new family of first-order methods with superlinear convergence. {\osgm} has global, local, and trajectory-dependent convergence guarantees, laying the theoretical foundation for hypergradient-descent-type first-order methods.
We hope {\osgm} paves the way for new directions in the design of adaptive optimization methods.

%% file: app_landscape.tex
\section{Proof of results in \Cref{sec:feedback}} 

\subsection{Proof of Proposition \ref{prop:property-ux}}
\input{proofs/proof-property-ux.tex}   

\subsection{Proof of Lemma \ref{lem:property-rx}}
\input{proofs/proof-property-rx.tex}

\subsection{Proof of Lemma \ref{lem:property-hx}}
\input{proofs/proof-property-hx.tex} 

\subsection{Proof of Equation \ref{eq:opt-preconditioner}} \label{sec:proof-opt-preconditioner}
\input{proofs/proof-opt-preconditioner.tex}

\subsection{Proof of Proposition \ref{prop:minimax-feedback}}
\input{proofs/proof-minimax-feedback.tex}

%% file: app_actions.tex
\section{Proof of results in Section \ref{sec:action}}

\subsection{Proof of Theorem \ref{thm:reduction-ratio}}

\input{proofs/proof-reduction-ratio.tex}
\subsection{Proof of Theorem \ref{thm:reduction-hypergrad}}

\input{proofs/proof-reduction-hypergrad.tex}\subsection{Proof of Lemma \ref{lem:feedback-progress}}
\input{proofs/proof-feedback-progress.tex}

%% file: app_oco.tex
\section{Proof of results in \Cref{sec:oco}}

\subsection{Proof of Lemma \ref{lem:ogd-sqrtK-eta}}
\input{proofs/proof-ogd-static-regret.tex}

\subsection{Proof of Lemma \ref{lem:dynamic-sqrtK}}
\input{proofs/proof-ogd-dynamic-regret.tex}

%% file: proofs/proof-ogd-static-regret.tex
For any $\hat{P} \in \Pcal$, online gradient descent update $P_{k + 1} = \Pi_{\mathcal{P}} [P_k - \eta_k \nabla \ell_{x^k} (P_k)]$ satisfies
\begin{align}
  \| P_{k + 1} - \hat{P} \|_F^2 
  ={} & \| \Pi_{\mathcal{P}} [P_k - \eta_k \nabla \ell_{x^k} (P_k)] -
  \hat{P} \|_F^2 \nonumber\\
  \leq{} & \| P_k - \eta_k \nabla \ell_{x^k} (P_k) -
  \hat{P} \|_F^2 \tag{by non-expansiveness of projection}\\
  ={} & \| P_k - \hat{P} \|_F^2 - 2 \eta_k \langle \nabla \ell_{x^k} (P_k), P_k -
  \hat{P} \rangle + \eta_k^2 \| \nabla \ell_{x^k} (P_k) \|^2_F \nonumber\\
  \leq{} & \| P_k - \hat{P} \|_F^2 - 2 \eta_k [\ell_{x^k} (P_k) - \ell_{x^k}
  (\hat{P})] + \eta_k^2 \| \nabla \ell_{x^k} (P_k) \|^2_F, \label{eqn:pf-ogd-static-2}
\end{align}
where \eqref{eqn:pf-ogd-static-2} uses convexity $\ell_{x^k} (\hat{P}) \geq \ell_{x^k} (P_k) + \langle \nabla \ell_{x^k} (P_k), \hat{P} - P_k \rangle$. Dividing both sides by $2 \eta_k$, we get
\[ \ell_{x^k} (P_k) \leq{} \ell_{x^k} (\hat{P}) + \tfrac{1}{2 \eta_k} [\| P_k -
   \hat{P} \|_F^2 - \| P_{k + 1} - \hat{P} \|_F^2] + \tfrac{\eta_k}{2} \| \nabla
   \ell_{x^k} (P_k) \|^2_F . \]
Telescoping the relation from $k = 1$ to $K$ and dropping the negative term $-\tfrac{1}{2 \eta_K} \| P_{K + 1} - \hat{P} \|_F^2$, we have
\begin{align}
   \textstyle \sum_{k = 1}^K \ell_{x^k} (P_k) &\leq \textstyle\sum_{k = 1}^K \ell_{x^k} (\hat{P}) +
   \tfrac{1}{2 \eta_1} \| P_1 - \hat{P} \|^2_F + \sum_{k = 1}^{K - 1} (
   \tfrac{1}{2 \eta_{k + 1}} - \tfrac{1}{2 \eta_k} ) \| P_{k + 1} -
   \hat{P} \|^2_F + \sum_{k = 1}^K \tfrac{\eta_k}{2} \| \nabla \ell_k (P_k)
   \|_F^2 \nonumber\\
   &\leq \textstyle \sum_{k = 1}^K \ell_{x^k} (\hat{P}) + \tfrac{1}{2 \eta_1} \| P_1 -
   \hat{P} \|^2_F + \sum_{k = 1}^K \tfrac{\eta_k}{2} \| \nabla \ell_k (P_k) \|_F^2, \label{eqn:pf-ogd-static-3}
\end{align}
where \eqref{eqn:pf-ogd-static-3} uses the fact $\eta_{k+1} \leq \eta_k$. 
Equation \eqref{eqn:pf-ogd-static-3} reduces to \eqref{eqn:ogd-regret-1} when $\eta_k \equiv \eta$. Now, suppose that $\{ \ell_{x^k} (P) \}$ are $\sigma$-Lipschitz and
$\tmop{diam} (\mathcal{P}) \leq D$. Then $\| \nabla \ell_k (P_k) \|_F^2
\leq \sigma^2$ and $\| P_1 - \hat{P} \|^2_F \leq D^2$. And \eqref{eqn:pf-ogd-static-3} implies
\[  \textstyle \sum_{k = 1}^K \ell_{x^k} (P_k) \leq \sum_{k = 1}^K \ell_{x^k} (\hat{P}) + \tfrac{D^2}{2 \eta_1} + \tfrac{\sigma^2}{2} \sum_{k = 1}^K \eta_k . \]
\begin{itemize}[leftmargin=10pt]
   \item The choice of the stepsize $\eta_k \equiv \tfrac{c}{\sqrt{K}}$ gives
   \begin{equation*}
      \textstyle \sum_{k = 1}^K \ell_{x^k} (P_k) \leq \sum_{k = 1}^K \ell_{x^k} (\hat{P}) + ( \tfrac{D^2}{2 c} + \tfrac{c \sigma^2}{2} ) \sqrt{K}.
   \end{equation*}
   \item The choice of the stepsize $\eta_k = \tfrac{c}{\sqrt{k}}$ and the fact $\sum_{k = 1}^K \tfrac{1}{\sqrt{k}} \leq 2 \sqrt{K}$ together give
   \begin{equation*}
     \textstyle \sum_{k = 1}^K \ell_{x^k} (P_k) \leq \sum_{k = 1}^K \ell_{x^k} (\hat{P}) +
      \tfrac{D^2}{2 c} + \tfrac{c \sigma^2}{2} \sum_{k = 1}^K \tfrac{1}{\sqrt{k}}
      \leq \sum_{k = 1}^K \ell_{x^k} (\hat{P}) + \tfrac{D^2}{2 c} + c \sigma^2 \sqrt{K}.
   \end{equation*}
\end{itemize}
Either case can be further bounded by the right-hand side of \eqref{eqn:ogd-regret-2}.

%% file: proofs/proof-ogd-dynamic-regret.tex
Recall from \eqref{eqn:pf-ogd-static-2} that for any $\hat{P} \in \mathcal{P}$ we have
\[ \| P_{k + 1} - \hat{P} \|_F^2 \leq \| P_k - \hat{P} \|_F^2 - 2 \eta [\ell_{x^k} (P_k) -
   \ell_{x^k} (\hat{P})] + \eta^2 \| \nabla \ell_{x^k} (P_k) \|_F^2 . \]
Plug in $\hat{P} = \hat{P}_k$ to obtain $\| P_{k + 1} - \hat{P}_k \|_F^2 \leq \|
P_k - \hat{P}_k \|_F^2 - 2 \eta [\ell_{x^k} (P_k) - \ell_{x^k} (\hat{P}_k)] +
\eta^2 \| \nabla \ell_{x^k} (P_k) \|_F^2$ and telescoping from $k = 1$ to $K$ yields
\begin{align}
  & \textstyle \sum_{k = 1}^K \ell_{x^k} (P_k) - \textstyle \sum_{k = 1}^K \ell_{x^k} (\hat{P}_k) -
  \tfrac{\eta}{2} \textstyle \sum_{k = 1}^K \| \nabla \ell_k (P_k) \|_F^2 \nonumber\\
  \leq{} & \tfrac{\| \hat{P}_1 - P_1 \|^2_F - \| \hat{P}_{K} - P_{K + 1}
  \|^2_F}{2 \eta} + \tfrac{1}{2 \eta} \textstyle \sum_{k = 1}^{K-1} [\| \hat{P}_{k + 1} -
  P_{k + 1} \|_F^2 - \| \hat{P}_k - P_{k + 1} \|_F^2] . \label{eqn:pf-ogd-dynamic-1}
\end{align}
Observe that the sum in the last term can be simplified to
\begin{align}
  & \textstyle \sum_{k = 1}^{K-1} [\| \hat{P}_{k + 1} - P_{k + 1} \|_F^2 - \| \hat{P}_k -
  P_{k + 1} \|_F^2] \nonumber\\
  ={} & \textstyle \sum_{k = 1}^{K-1} [\| \hat{P}_{k + 1} \|_F^2 - \| \hat{P}_k \|_F^2 + 2
  \langle P_{k + 1}, \hat{P}_k - \hat{P}_{k + 1} \rangle] \nonumber\\
  ={} & \textstyle \sum_{k = 1}^{K-1} [\| \hat{P}_{k + 1} \|_F^2 - \| \hat{P}_k \|_F^2 + 2
  \langle P_{k + 1} - P_1, \hat{P}_k - \hat{P}_{k + 1} \rangle + 2 \langle
  P_1, \hat{P}_k - \hat{P}_{k + 1} \rangle] \nonumber\\
  ={} & \| \hat{P}_{K} \|_F^2 - \| \hat{P}_1 \|_F^2 + 2 \textstyle \sum_{k = 1}^{K-1}
  \langle P_{k + 1} - P_1, \hat{P}_k - \hat{P}_{k + 1} \rangle + 2 \langle
  P_1, \hat{P}_1 - \hat{P}_{K} \rangle \nonumber\\
  \leq{} & \| \hat{P}_{K} \|_F^2 - \| \hat{P}_1 \|_F^2 + 2 \max_{k \leq K} \| P_{k} - P_1 \|_F \textstyle \sum_{k = 1}^{K-1} \| \hat{P}_k - \hat{P}_{k + 1} \|_F + 2 \langle
  P_1, \hat{P}_1 - \hat{P}_{K} \rangle \nonumber\\
  ={} & \| \hat{P}_{K} - P_1 \|_F^2 - \| \hat{P}_1 - P_1 \|_F^2 + 2 \max_{k \leq K}  \| P_{k} - P_1 \|_F \textstyle \sum_{k = 1}^{K-1} \| \hat{P}_k - \hat{P}_{k + 1} \|_F . \label{eqn:pf-ogd-dynamic-2}
\end{align}
Combining \eqref{eqn:pf-ogd-dynamic-1} and \eqref{eqn:pf-ogd-dynamic-2}, we get
\begin{align}
  & \textstyle \sum_{k = 1}^K \ell_{x^k} (P_k) - \textstyle \sum_{k = 1}^K \ell_{x^k} (\hat{P}_k) -
  \tfrac{\eta}{2} \textstyle \sum_{k = 1}^K \| \nabla \ell_k (P_k) \|_F^2 \nonumber\\
  \leq{} & \tfrac{\| \hat{P}_{K} - P_1 \|^2_F - \| \hat{P}_{K} - P_{K +
  1} \|^2_F}{2 \eta} + \tfrac{\max_{k \leq K} \| P_{k} - P_1 \|_F}{\eta} \textstyle \sum_{k = 1}^{K-1} \| \hat{P}_k - \hat{P}_{k + 1} \|_F . \nonumber
\end{align}
Moving $\tfrac{\eta}{2} \textstyle \sum_{k = 1}^K \| \nabla \ell_{x^k} (P_k) \|_F^2$ to the right and dropping the negative term $-\tfrac{\| \hat{P}_{K} - P_{K +1} \|^2_F}{2 \eta}$ prove \eqref{eqn:dynamic-regret-1}. 
When $\{ \ell_{x^k}
(P_k) \}$ are $\sigma$-Lipschitz, we can bound $\| \nabla \ell_{x^k} (P_k) \|_F^2 \leq \sigma^2$ and plug in $\eta = \tfrac{c}{\sqrt{K}}$ to arrive at \eqref{eqn:dynamic-regret-2}.

%% file: app_algodemo.tex
\section{Proof of results in Section \ref{sec:algodemos}}

\subsection{Proof of Theorem \ref{thm:algodemo-global-conv}}
\input{proofs/proof-global-conv-handson-rx.tex}

\subsection{Proof of Theorem \ref{thm:potential-handson-rx}}
\input{proofs/proof-potential-handson-rx.tex}

\subsection{Proof of Theorem \ref{thm:local-trajectory-conv-handson-rx}}
\input{proofs/proof-local-trajectory-conv-rx.tex}

\subsection{Proof of Lemma \ref{lem:local-opthessian-rx}}
\input{proofs/proof-local-opthessian-rx.tex}

\subsection{Proof of Theorem \ref{thm:local-superlinear-conv-rx}}
\input{proofs/proof-local-superlinear-conv-rx.tex}

\subsection{Proof of Theorem \ref{thm:local-negative-regret-rx}}
\input{proofs/proof-local-negative-regret-rx.tex}

\subsection{Proof of Theorem \ref{thm:global-conv-handson-hx}}
\input{proofs/proof-global-conv-handson-hx.tex}

\subsection{Proof of Theorem \ref{thm:potential-handson-hx}}
\input{proofs/proof-potential-handson-hx.tex}

\subsection{Proof of Theorem \ref{thm:local-trajectory-conv-handsonhx}}
\input{proofs/proof-local-trajectory-conv-handson-hx.tex}

\subsection{Proof of Theorem \ref{thm:local-superlinear-conv-handson-hx}}
\input{proofs/proof-local-superlinear-conv-handson-hx.tex}

\subsection{Proof of Theorem \ref{thm:local-negative-regret-handson-hx}}
\input{proofs/proof-local-negative-regret-handson-hx.tex}

%% file: proofs/proof-global-conv-handson-rx.tex
\Cref{eqn:handsonrx-tjry-conv} follows by plugging the bound \eqref{eqn:demo-ratio-constant-regret} into \Cref{thm:reduction-ratio}.
Now, suppose {\osgmhandsonrx} is initialized with $P_1 = \tfrac{1}{L} I$. 
To show \eqref{eqn:algodemo-global-conv-rate}, we plug in $\hat{P} = \tfrac{1}{L} I$ into \eqref{eqn:handsonrx-tjry-conv} and use $r_x(\tfrac{1}{L} I) \leq 1 - \tfrac{1}{\kappa}$ to obtain
\begin{equation} \label{eqn:pf-global-conv-handsonrx-1}
    \textstyle f (x^{K + 1}) - f^{\star} 
    \leq [f (x^1) - f^{\star}] ( \tfrac{1}{K} \sum_{k = 1}^K r_{x^k} (\tfrac{1}{L} I) + \tfrac{L^2}{K}  \| \tfrac{1}{L} I - \tfrac{1}{L} I \|_F^2 )^K 
    \leq [f (x^1) - f^{\star}] (1 - \tfrac{1}{\kappa})^K.
\end{equation}
On the other hand, we can plug in $\hat{P} = P^{\star}_r$ into \eqref{eqn:handsonrx-tjry-conv} and use $r_x(P^{\star}_r) \leq 1 - \tfrac{1}{\kappa^\star}$ to obtain
\begin{align} 
    \textstyle f (x^{K + 1}) - f^{\star} 
    &\leq [f (x^1) - f^{\star}] ( \tfrac{1}{K} \textstyle \sum_{k = 1}^K r_{x^k} (P^{\star}_r) + \tfrac{L^2}{K}  \| \tfrac{1}{L} I - P^{\star}_r \|_F^2 )^K \nonumber\\
    &\leq [f (x^1) - f^{\star}] ( 1 - \tfrac{1}{\kappa^{\star}} + \tfrac{L^2}{K} \| \tfrac{1}{L} I - P_r^{\star} \|_F^2 )^K. \label{eqn:pf-global-conv-handsonrx-2}
\end{align}
Take the minimum of two bounds in \eqref{eqn:pf-global-conv-handsonrx-1} and \eqref{eqn:pf-global-conv-handsonrx-2} to conclude \eqref{eqn:algodemo-global-conv-rate}.

%% file: proofs/proof-local-trajectory-conv-rx.tex
The result follows from the following chain of inequalities:
\begin{align*}
  {} & f (x^{K + 1}) - f^{\star} \\
  \leq{} & [f (x^1) - f^{\star}] (\tfrac{1}{K} \textstyle \sum_{k = 1}^K r_k )^K \tag{by \Cref{thm:reduction-ratio}}\\
  \leq{} & [f (x^1) - f^{\star}] (\tfrac{1}{K} \textstyle \sum_{k = 1}^K r_{x^k}(P_k) - \tfrac{1}{4 L^2} \| \nabla r_{x^k} (P_k) \|_F^2)^K \tag{by \Cref{lem:feedback-progress} $+$ lookahead action}\\
  ={} & [f (x^1) - f^{\star}] (\tfrac{1}{K} \textstyle \sum_{k = 1}^K r_{x^k}(P_k) - \tfrac{\eta}{2} \| \nabla r_{x^k} (P_k) \|_F^2)^K \tag{by $\eta = \tfrac{1}{2L^2}$}\\
  \leq{} & [f (x^1) - f^{\star}] (\tfrac{1}{K} \textstyle \sum_{k = 1}^K r_{x^k}
  (\hat{P}_k) + \tfrac{L^2}{K} \rho_{K}(\{\hat{P}_k \}))^K, \tag{by \Cref{lem:dynamic-sqrtK}}
\end{align*}
where $\rho_{K}(\{\hat{P}_k \}) \assign \| \hat{P}_{K + 1} - P_1 \|^2_F + 2  \max_{k \leq K} \{\| P_{k} - P_1 \|_F\} \mathsf{PL}(\{ \hat{P}_k \})$.

%% file: proofs/proof-local-negative-regret-rx.tex
Fix a benchmark stepsize $\hat{P}$ and $\eta \leq \tfrac{1}{4 L^2}$. 
Recall from \eqref{eqn:algodemo-together} that
\begin{align}
  \textstyle \sum_{k = 1}^K r_k \leq{} & \textstyle \sum_{k = 1}^K r_{x^k} (\hat{P}) + \tfrac{1}{2
  \eta} \| P_1 - \hat{P} \|_F^2 + \tfrac{1}{2} ( \eta - \tfrac{1}{2 L^2}
  ) \textstyle \sum_{k = 1}^K \| \nabla r_{x^k} (P_k) \|_F^2 \tag{by \eqref{eqn:algodemo-together}}\\
  \leq{} & \textstyle \sum_{k = 1}^K r_{x^k} (\hat{P}) + \tfrac{1}{2 \eta} \| P_1 - \hat{P}
  \|_F^2 - \tfrac{1}{8 L^2} \textstyle \sum_{k = 1}^K \| \nabla r_{x^k} (P_k) \|_F^2.
  \tag{by $\eta \leq \tfrac{1}{4 L^2}$}
\end{align}
To lower bound the gradient norm of ratio feedback, observe that
\begin{align}
  \| \nabla r_{x^k} (P_k) \|_F^2 ={} \tfrac{\| \nabla f (x^k) \|^2 \| \nabla f
  (x^{k + 1 / 2}) \|^2}{[f (x^k) - f^{\star}]^2}
  ={} \tfrac{\| \nabla f (x^{k + 1}) \|^2}{\| \nabla f (x^k) \|^2} \tfrac{\|
  \nabla f (x^{k + 1 / 2}) \|^2}{\| \nabla f (x^{k + 1}) \|^2} \tfrac{\|
  \nabla f (x^k) \|^4}{[f (x^k) - f^{\star}]^2}.
  \nonumber
\end{align}
The middle fraction can be bounded using $\|\nabla f(x - \tfrac{1}{L} \nabla f(x)) \| \leq \|\nabla f(x)\|$; the last fraction can be lower bounded using the relation $\frac{1}{2 \mu} \| \nabla f (x^k) \|^2 \geq f (x^k) - f^{\star}$. 

Hence, we have $\| \nabla r_{x^k} (P_k) \|_F^2 \geq 4 \mu^2 \tfrac{\| \nabla f (x^{k + 1}) \|^2}{\| \nabla f (x^k) \|^2}$ and the bound on the cumulative progress becomes
\begin{align}
  \textstyle \sum_{k = 1}^K r_k
  \leq{} & \textstyle \sum_{k = 1}^K r_{x^k} (\hat{P}) + \tfrac{1}{2 \eta} \| P_1 - \hat{P}
  \|_F^2 - \tfrac{1}{8 L^2} \textstyle \sum_{k = 1}^K \| \nabla r_{x^k} (P_k) \|_F^2 \nonumber\\
  \leq{} & \textstyle \sum_{k = 1}^K r_{x^k} (\hat{P}) + \tfrac{1}{2 \eta} \| P_1 - \hat{P}\|_F^2 - \tfrac{1}{2 \kappa^2} \textstyle \sum_{k = 1}^K \tfrac{\| \nabla f (x^{k + 1}) \|^2}{\| \nabla f (x^k) \|^2}. \nonumber
\end{align}
A rearrangement gives
\begin{equation} \label{eqn:pf-local-negative-regret-rx-1}
  \textstyle \sum_{k = 1}^K \tfrac{\| \nabla f (x^{k + 1}) \|^2}{\| \nabla f (x^k) \|^2} \leq 2 \kappa^2 [ \textstyle \sum_{k = 1}^K r_{x^k}(\hat{P}) - \textstyle \sum_{k = 1}^K r_k ] + \tfrac{\kappa^2}{\eta} \| P_1 - \hat{P} \|_F^2.   
\end{equation}
Then we can bound the suboptimality by 
\begin{align}
  f (x^{K + 1}) - f^{\star} \leq{} & \tfrac{1}{2 \mu} \| \nabla f (x^{K + 1}) \|^2
  \nonumber\\
  \leq{} & \tfrac{1}{2 \mu} \| \nabla f (x^1) \|^2 ( \tfrac{1}{K} \textstyle \sum_{k =
  1}^K \tfrac{\| \nabla f (x^{k + 1}) \|^2}{\| \nabla f (x^k) \|^2} )^K
  \tag{by AM-GM inequality}\\
  \overleq{\eqref{eqn:pf-local-negative-regret-rx-1}} & \tfrac{1}{2 \mu} \| \nabla f (x^1) \|^2 ( \tfrac{2 \kappa^2}{K}
  [ \textstyle \sum_{k = 1}^K r_{x^k} (\hat{P}) - \textstyle \sum_{k = 1}^K r_k ] + \tfrac{\kappa^2}{\eta K} \| P_1 - \hat{P} \|_F^2 )^K. \label{eqn:pf-local-negative-regret-rx-2}
\end{align}

Finally, we do a case analysis on the sign of $\textstyle \sum_{k = 1}^K r_{x^k} (\hat{P}) - \textstyle \sum_{k = 1}^K r_k$.

\paragraph{Case 1.} $\textstyle \sum_{k = 1}^K r_k \leq \textstyle \sum_{k = 1}^K r_{x^k} (\hat{P})$: The progress of {\osgmhandsonrx} is better than stepsize $\hat{P}$.

\paragraph{Case 2.} $\textstyle \sum_{k = 1}^K r_k \geq \textstyle \sum_{k = 1}^K r_{x^k} (\hat{P})$: relation \eqref{eqn:pf-local-negative-regret-rx-2} can be further bounded by
\[ f (x^{K + 1}) - f (x^1) \leq \tfrac{1}{2 \mu} \| \nabla f (x^1) \|^2
  ( \tfrac{\kappa^2 \| P_1 - \hat{P} \|_F^2}{\eta K} )^K. \]

%% file: proofs/proof-potential-handson-hx.tex
For brevity, we drop the subscript $\varphi$ and $\omega$ for parameters $\rho_{\varphi}$ and $\rho_{\omega}$ in both potential functions.

\paragraph{Convex $f$.} We analyze the potential function
\[ \omega (x, P) = - \tfrac{\rho}{f (x) - f^{\star}} + \| P - \tfrac{1}{L} I \|^2_F . \]
At every iteration, the function value part of the potential changes by
\begin{align}
  \tfrac{1}{f (x^{k + 1}) - f^{\star}} & = \tfrac{1}{f (x^k) - f^{\star}} - \tfrac{\| \nabla f (x^k) \|^2 h_k}{[f (x^{k + 1}) - f^{\star}] [f (x^k) - f^{\star}]} \nonumber\\
  & \geq \tfrac{1}{f (x^k) - f^{\star}} - \tfrac{\| \nabla f (x^k) \|^2
  h_k}{[f (x^k) - f^{\star}]^2} \tag{by $f(x^{k+1}) \leq f(x^k)$ and $h_k \leq 0$}\\
  & \geq \tfrac{1}{f (x^k) - f^{\star}} - \tfrac{h_k}{\Delta^2} .
  \tag{by $\tfrac{f(x^k) - f^{\star}}{\| \nabla f(x^k)\|^2 } \geq \Delta^2$}
\end{align}
On the other hand, the distance part changes by
\begin{align}
  \| P_{k + 1} - \tfrac{1}{L} I \|^2_F & \leq{} \| P_k - \tfrac{1}{L} I \|^2_F - 2 \eta
  [h_{x^k} (P_k) - h_{x^k} (\tfrac{1}{L} I)] + \eta^2 \| \nabla h_{x^k} (P_k) \|^2_F
  \nonumber\\
  & \leq{} \| P_k - \tfrac{1}{L} I \|^2_F - 2 \eta [ h_{x^k} (P_k) - \tfrac{\eta}{2} \| \nabla h_{x^k} (P_k) \|^2_F ] - \tfrac{\eta}{L} \tag{by descent lemma $h_{x^k}(\tfrac{1}{L} I) \leq - \tfrac{1}{2 L}$}\\
  & \leq{} \| P_k - \tfrac{1}{L} I \|^2_F - 2 \eta [ h_{x^k} (P_k) - \tfrac{1}{2L} \| \nabla h_{x^k} (P_k) \|^2_F ] - \tfrac{\eta}{L} \tag{by $\eta \leq \tfrac{1}{L}$}\\
  & \leq{} \| P_k - \tfrac{1}{L} I \|^2_F - 2 \eta h_k - \tfrac{\eta}{L} . \tag{by \Cref{lem:feedback-progress} $+$ lookahead action}
\end{align}
Combining both parts to obtain the change of potential:
\begin{align}
  \omega (x^{k + 1}, P_{k + 1}) - \omega (x^k, P_k) 
  \leq \tfrac{\rho h_k}{\Delta^2} - 2 \eta h_k - \tfrac{\eta}{L}
  = ( \tfrac{\rho}{\Delta^2} - 2 \eta ) h_k - \tfrac{\eta}{L} .
  \nonumber
\end{align}
For $\rho \geq 2 \eta \Delta^2$, the potential function will strictly
decrease. Taking $\rho = 2 \eta \Delta^2$ and $\eta = \tfrac{1}{L}$, we have
\[ \omega (x^{k + 1}, P_{k + 1}) - \omega (x^k, P_k) \leq - \tfrac{1}{L^2}. \]
To obtain the iteration complexity, note that 
\begin{align}
  - \tfrac{1}{L} \tfrac{2 \Delta^2}{f (x^{K + 1}) - f^{\star}} + \| P_{K + 1} - \tfrac{1}{L} I \|_F^2 \leq & - \tfrac{1}{L} \tfrac{2 \Delta^2}{f (x^1) - f^{\star}} + \| P_1 - \tfrac{1}{L} I \|_F^2 - \tfrac{1}{L^2} K. \nonumber
\end{align}
Hence, {\osgmhandsonhx} with $P_1 = \tfrac{1}{L} I$ achieves $f (x^{K + 1}) - f^{\star} \leq \varepsilon$ for
\[K \geq L^2 \| P_1 - \tfrac{1}{L} I \|_F^2 + \tfrac{2 L \Delta^2}{\varepsilon} = \tfrac{2 L \Delta^2}{\varepsilon}.\]

\paragraph{$\mu$-strongly convex $f$.} We analyze the potential function
\[ \varphi (x, P) = \rho \log (f (x) - f^{\star}) + \| P - \tfrac{1}{L} I \|_F^2. \]
At every iteration, the function value part of the
potential changes by
\begin{align}
  f (x^{k + 1}) - f^{\star} ={} & f (x^{k + 1}) - f (x^k) + f (x^k) - f^{\star} \nonumber\\
  ={} & h_k \| \nabla f (x^k) \|^2 + f (x^k) - f^{\star} \nonumber\\
  ={} & ( 1 + h_k \tfrac{\| \nabla f (x^k) \|^2}{f (x^k) - f^{\star}}
  ) (f (x^k) - f^{\star}) \nonumber\\
  \leq{} & (1 + 2 \mu h_k) (f (x^k) - f^{\star}), \nonumber
\end{align}
where the last inequality uses $\tfrac{\| \nabla f (x^k) \|^2}{f (x^k) - f^{\star}} \geq 2
\mu$. 
The change of distance part $\| P - \tfrac{1}{L} I \|^2_F$ is the same as convex case.
Combine both parts to obtain the change of potential:
\begin{align}
  & \varphi (x^{k + 1}, P_{k + 1}) - \varphi (x^k, P_k) \nonumber\\
  ={} & \rho \log (f (x^{k + 1}) - f^{\star}) - \rho \log (f (x^k) - f^{\star}) + \| P_{k + 1} - \tfrac{1}{L} I \|^2_F - \| P_k - \tfrac{1}{L} I \|^2_F
  \nonumber\\
  \leq{} & \rho \log (1 + 2 \mu h_k) - 2 \eta h_k - \tfrac{\eta}{L} = \alpha(h_x) - \tfrac{\eta}{L},  \nonumber
\end{align}
where the function $\alpha (x) \assign \rho \log (1 + 2 \mu x) - 2 \eta x$ is
maximized at $x = \frac{\mu \rho - \eta}{2 \eta \mu}$. Taking $\rho = \eta /
\mu$, we get
\begin{equation*}
  \varphi (x^{k + 1}, P_{k + 1}) - \varphi (x^k, P_k) 
  \leq \rho \log (1 + \tfrac{\mu \rho - \eta}{\eta}) - \tfrac{\mu \rho - \eta}{\mu} - \tfrac{\eta}{L} = - \tfrac{\eta}{L}.
\end{equation*}
Take $\eta = 1 / L$ to get
\[ \varphi (x^{k + 1}, P_{k + 1}) - \varphi (x^k, P_k) \leq - \tfrac{1}{L^2}. \]
To obtain the iteration complexity, note that 
\begin{align}
   \tfrac{1}{\mu L} \log (f (x^{K + 1}) - f^{\star}) + \| P_{K + 1} -
  \tfrac{1}{L} I \|_F^2  \leq \tfrac{1}{\mu L} \log [f (x^1) - f^{\star}] + \| P_1 - \tfrac{1}{L} I
  \|_F^2 - \tfrac{1}{L^2} K. \nonumber
\end{align}
Hence, {\osgmhandsonhx} with $P_1 = \tfrac{1}{L} I$ achieves $f (x^{K + 1}) - f^{\star} \leq \varepsilon$ for 
\[ K \geq L^2 \| P_1 - \tfrac{1}{L} I \|_F^2 + \kappa \log ( \tfrac{f(x^1) - f (x)}{\varepsilon} ) = \kappa \log ( \tfrac{f(x^1) - f (x)}{\varepsilon} ). \]

%% file: proofs/proof-local-trajectory-conv-handson-hx.tex
By \eqref{eqn:pf-global-conv-handson-hx-0} and \Cref{lem:dynamic-sqrtK}, the cumulative progress of {\osgmhandsonhx} is bounded by
\begin{align}
\textstyle \sum_{k = 1}^K h_k
\leq{} & \textstyle \sum_{k = 1}^K h_{x^k} (P_k) - \tfrac{1}{2 L} \textstyle \sum_{k = 1}^K \| \nabla h_{x^k} (P_k) \|_F^2 \tag{by \eqref{eqn:pf-global-conv-handson-hx-0}}\\
\leq{} & \textstyle \sum_{k = 1}^K h_{x^k} (\hat{P}_k) + \tfrac{\| \hat{P}_{K} - P_1 \|^2_F}{2 \eta} + \tfrac{\max_{k \leq K} \| P_{k} - P_1 \|_F}{\eta} \mathsf{PL} (\{\hat{P}_k \}) \tag{by \Cref{lem:dynamic-sqrtK}} \\
\leq{} & \textstyle \sum_{k = 1}^K h_{x^k} (\hat{P}_k) + \tfrac{L}{2} [ \| \hat{P}_{K} - P_1 \|_F^2 + 2 \max_{k \leq K} \{\| P_{k} - P_1 \|_F\} \mathsf{PL} (\{\hat{P}_k \})] \tag{by $\eta=\tfrac{1}{L}$}.
\end{align}
Plugging the minimum of this bound and the bound $h_k \leq 0$ into hypergradient reduction from \Cref{thm:reduction-hypergrad} completes the proof.

%% file: proofs/proof-local-superlinear-conv-handson-hx.tex
First, we establish a bound on
$h_x ([\nabla^2 f (x^{\star})]^{- 1})$ below.

\begin{lem} \label{lem:local-opthessian-hx}
Suppose $f$ is $L$-smooth $\mu$-strongly convex has $H$-Lipschitz Hessian. 
Then the hypergradient feedback of Hessian inverse at $x^\star$ satisfies $h_x ([\nabla^2 f (x^{\star})]^{- 1}) + \tfrac{f (x) - f^{\star}}{\| \nabla f (x) \|^2} \leq \tfrac{H^2 \kappa}{8 \mu^3} \| x - x^{\star} \|^2$ for all $x \nin \Xcal^\star$.
\end{lem}	
\begin{proof}
Notice that
\[ h_x ([\nabla^2 f (x^{\star})]^{- 1}) + \tfrac{f (x) - f^{\star}}{\|
   \nabla f (x) \|^2} = \tfrac{f (x - [\nabla^2 f (x^{\star})]^{- 1} \nabla f
   (x)) - f^{\star}}{\| \nabla f (x) \|^2} . \]
Using \eqref{eqn:local-opthessian-rx-1} and $\tfrac{\| x
   - x^{\star} \|^2}{\| \nabla f (x) \|^2} \leq \tfrac{1}{\mu^2}$, the right-hand side can be further bounded by
\[ \tfrac{f (x - [\nabla^2 f (x^{\star})]^{- 1} \nabla f (x)) - f
   (x^{\star})}{\| \nabla f (x) \|^2} \overset{\eqref{eqn:local-opthessian-rx-1}}{\leq} \tfrac{L H^2}{8 \mu^2} \tfrac{\| x
   - x^{\star} \|^4}{\| \nabla f (x) \|^2} \leq \tfrac{L H^2}{8 \mu^4} \| x -
   x^{\star} \|^2 = \tfrac{H^2 \kappa}{8 \mu^3} \| x - x^{\star} \|^2. \]
\end{proof}

Now, we are ready to prove \Cref{thm:local-superlinear-conv-handson-hx}.
Note that
\begin{align}
 \textstyle f (x^{K + 1}) - f (x^k) \leq{} & [f (x^1) - f^{\star}] ( \tfrac{1}{K}
  \textstyle\sum_{k = 1}^K \tfrac{f (x^{k + 1}) - f^{\star}}{f (x^k) - f
  (x^{\star})} ) \tag{by AM-GM inequality}\\
  ={} & [f (x^1) - f^{\star}] ( \tfrac{1}{K} \textstyle \sum_{k = 1}^K \tfrac{f
  (x^{k + 1}) - f^{\star}}{\| \nabla f (x^k) \|^2} \tfrac{\| \nabla f
  (x^k) \|^2}{f (x^k) - f^{\star}} )^K \nonumber\\
  \leq{} & [f (x^1) - f^{\star}] ( \tfrac{2 L}{K} \textstyle \sum_{k = 1}^K
  \tfrac{f (x^{k + 1}) - f^{\star}}{\| \nabla f (x^k) \|^2} )^K
  \tag{by $L$-smoothness of $f$}\\
  ={} & [f (x^1) - f^{\star}] ( \tfrac{2 L}{K} \textstyle \sum_{k = 1}^K [
  h_k + \tfrac{f (x^k) - f^{\star}}{\| \nabla f (x^k) \|^2} ]
  )^K . \label{eqn:pf-local-superlinear-conv-handson-hx-1}
\end{align}
We now further bound the sum on the right by
\begin{align}
 \textstyle \sum_{k = 1}^K [ h_k + \tfrac{f (x^k) - f^{\star}}{\| \nabla f
  (x^k) \|^2} ] \leq{} & \textstyle \sum_{k = 1}^K [ h_{x^k} ([\nabla^2 f
  (x^{\star})]^{- 1}) + \tfrac{f (x^k) - f^{\star}}{\| \nabla f (x^k)
  \|^2} ] + \tfrac{L}{2} \| P_1 - [\nabla^2 f (x^{\star})]^{- 1} \|_F^2
  \tag{by \eqref{eqn:pf-global-conv-handson-hx-1}}\\
  \leq{} & \textstyle \tfrac{H^2 \kappa}{8 \mu^3} \sum_{k = 1}^K \| x^k - x^{\star} \|^2 +
  \tfrac{L}{2} \| P_1 - [\nabla^2 f (x^{\star})]^{- 1} \|_F^2 . \tag{by \Cref{lem:local-opthessian-hx}}
\end{align}

The linear convergence rate of {\osgmhandsonhx} in \Cref{thm:global-conv-handson-hx} implies
\begin{align*}
  \textstyle \sum_{k = 1}^K \| x^k - x^{\star} \|^2 \leq{} & \tfrac{2}{\mu} \textstyle \sum_{k = 1}^K [f(x^k) - f^{\star}] \tag{by $\mu$-strong convexity of $f$} \\
  \leq{} & \tfrac{2}{\mu} [f(x^1) - f^{\star}] \textstyle \sum_{k = 1}^K (1 - \tfrac{1}{\kappa})^k
  \leq \tfrac{2 \kappa}{\mu} [f(x^1) - f^{\star}]. \tag{by \Cref{thm:global-conv-handson-hx}}
\end{align*}
The above two inequalities, together with the choice $P_1 = \tfrac{1}{L} I$, imply
\begin{align}
  \textstyle \sum_{k = 1}^K [ h_k + \tfrac{f (x^k) - f^{\star}}{\| \nabla f(x^k) \|^2} ]
  \leq{} & \textstyle \tfrac{H^2 \kappa^2}{4 \mu^4} [f(x^1) - f^{\star}] + \tfrac{L}{2} \| \tfrac{1}{L} I - [\nabla^2 f (x^{\star})]^{- 1} \|_F^2. \label{eqn:pf-local-superlinear-conv-handson-hx-2}
\end{align}
Plugging \eqref{eqn:pf-local-superlinear-conv-handson-hx-2} back into \eqref{eqn:pf-local-superlinear-conv-handson-hx-1} completes the proof:
\begin{align}
  f (x^{K + 1}) - f (x^k) \leq{} & [f (x^1) - f^{\star}] ( \tfrac{2
  L}{K} [ \textstyle \tfrac{H^2 \kappa^2}{4 \mu^4} [f(x^1) - f^{\star}] + \tfrac{L}{2} \| \tfrac{1}{L} I - [\nabla^2 f (x^{\star})]^{- 1} \|_F^2 ] )^K \nonumber\\
  ={} & [f (x^1) - f^{\star}] ( \tfrac{1}{K} [ \tfrac{H^2
  \kappa^3}{2 \mu^3} [f(x^1) - f^{\star}] + L^2 \| \tfrac{1}{L} I - [\nabla^2 f
  (x^{\star})]^{- 1} \|_F^2 ] )^K \nonumber \\
  ={} & [f (x^1) - f^{\star}] (\tfrac{C}{K})^K,\nonumber
\end{align}
where $C \coloneqq \tfrac{H^2
  \kappa^3}{2 \mu^3} [f(x^1) - f^{\star}] + L^2 \| \tfrac{1}{L} I - [\nabla^2 f
  (x^{\star})]^{- 1} \|_F^2 ]$.

%% file: proofs/proof-local-negative-regret-handson-hx.tex
Fix a benchmark stepsize $\hat{P}$ and $\eta \leq \tfrac{1}{2 L}$. 
Recall from \eqref{eqn:pf-global-conv-handson-hx-0.5} that
\begin{align}
  \textstyle \sum_{k = 1}^K h_k \leq{} & \textstyle \sum_{k = 1}^K h_{x^k} (\hat{P}) + \tfrac{1}{2
  \eta} \| P_1 - \hat{P} \|_F^2 + \tfrac{1}{2} ( \eta - \tfrac{1}{L}
  ) \textstyle \sum_{k = 1}^K \| \nabla h_{x^k} (P_k) \|_F^2 \tag{by \eqref{eqn:pf-global-conv-handson-hx-0.5}}\\
  \leq{} & \textstyle \sum_{k = 1}^K h_{x^k} (\hat{P}) + \tfrac{1}{2 \eta} \| P_1 - \hat{P}
  \|_F^2 - \tfrac{1}{4 L} \textstyle \sum_{k = 1}^K \| \nabla h_{x^k} (P_k) \|_F^2.
  \tag{by $\eta \leq \tfrac{1}{2 L}$} \\
  \leq{} & \textstyle \sum_{k = 1}^K h_{x^k} (\hat{P}) + \tfrac{1}{2 \eta} \| P_1 - \hat{P}
  \|_F^2 - \tfrac{1}{4 L} \textstyle \sum_{k = 1}^K \tfrac{\| \nabla f (x^{k + 1})
  \|^2}{\| \nabla f (x^k) \|^2}, \label{eqn:pf-local-negative-regret-hx-0}
\end{align}
where \eqref{eqn:pf-local-negative-regret-hx-0} applies the definition of $h_x$ and that $\|\nabla f(x^{k+1})\| = \|\nabla f(x^{k+1/2} - \frac{1}{L} \nabla f(x^{k+1/2})) \| \leq \|\nabla f(x^{k+1/2})\|$:
\begin{equation*}
   \| \nabla h_{x^k} (P_k) \|_F^2 = \tfrac{\| \nabla f (x^k)
   \|^2 \| \nabla f (x^{k + 1 / 2}) \|^2}{\| \nabla f (x^k) \|^4} = \tfrac{\|
   \nabla f (x^{k + 1 / 2}) \|^2}{\| \nabla f (x^k) \|^2} \geq \tfrac{\| \nabla f (x^{k + 1}) \|^2}{\| \nabla f (x^k) \|^2}.
\end{equation*}
Rearranging \eqref{eqn:pf-local-negative-regret-hx-0}, we get
\begin{equation} \label{eqn:pf-local-negative-regret-hx-1}
  \textstyle \sum_{k = 1}^K \tfrac{\| \nabla f (x^{k + 1}) \|^2}{\| \nabla f (x^k) \|^2} \leq 4 L [ \textstyle \sum_{k = 1}^K h_{x^k}(\hat{P}) - \textstyle \sum_{k = 1}^K h_k ] + \tfrac{2 L}{\eta} \| P_1 - \hat{P} \|_F^2.   
\end{equation}
Then we can bound the suboptimality by 
\begin{align}
  f (x^{K + 1}) - f^{\star} \leq{} & \tfrac{1}{2 \mu} \| \nabla f (x^{K + 1}) \|^2
  \nonumber\\
  \leq{} & \tfrac{1}{2 \mu} \| \nabla f (x^1) \|^2 ( \tfrac{1}{K} \textstyle \sum_{k =
  1}^K \tfrac{\| \nabla f (x^{k + 1}) \|^2}{\| \nabla f (x^k) \|^2} )^K
  \tag{by AM-GM inequality}\\
  \overleq{\eqref{eqn:pf-local-negative-regret-hx-1}} & \tfrac{1}{2 \mu} \| \nabla f (x^1) \|^2 ( \tfrac{4 L}{K}
  [ \textstyle \sum_{k = 1}^K h_{x^k} (\hat{P}) - \textstyle \sum_{k = 1}^K h_k ] + \tfrac{2 L}{\eta K} \| P_1 - \hat{P} \|_F^2 )^K. \label{eqn:pf-local-negative-regret-hx-2}
\end{align}

Finally, we do a case analysis on the sign of $\textstyle \sum_{k = 1}^K h_{x^k} (\hat{P}) - \textstyle \sum_{k = 1}^K h_k$.

\paragraph{Case 1.} $\textstyle \sum_{k = 1}^K h_k \leq \textstyle \sum_{k = 1}^K h_{x^k} (\hat{P})$: The progress of {\osgmhandsonhx} is better than stepsize $\hat{P}$.

\paragraph{Case 2.} $\textstyle \sum_{k = 1}^K h_k \geq \textstyle \sum_{k = 1}^K h_{x^k} (\hat{P})$: relation \eqref{eqn:pf-local-negative-regret-hx-2} can be further bounded by
\[ f (x^{K + 1}) - f (x^1) \leq \tfrac{1}{2 \mu} \| \nabla f (x^1) \|^2
   ( \tfrac{2 L \| P_1 - \hat{P} \|_F^2}{\eta K} )^K. \]

%% file: app_moreexps.tex
\section{Other instances of {\osgm}} \label{sec:other-instances}

This section presents additional instantiations of {\osgm}. We will invoke the following assumptions.
\begin{enumerate}[leftmargin=25pt,label=\textbf{A\arabic*:},ref=\rm{\textbf{A\arabic*}},start=1]
\item $f$ is $L$-smooth and convex.  \label{A1}
\item $f$ is $\mu$-strongly convex.  \label{A2}
\item $\Pcal$ satisfies $0 \in \Pcal, \frac{1}{L} I \in \Pcal$ and $\diam (\Pcal) \leq D < \infty$.  \label{A3}
\item $f$ has $H$-Lipschitz Hessian.  \label{A4}
\end{enumerate}

\subsection{\osgmhandsoffrx} \label{sec:hands-off-rx}

In this section, we assume $f$ is $L$-smooth and $\mu$-strongly convex (\ref{A1}, \ref{A2}) and instantiate {\osgm} with
\[ \ell_x (P) \assign r_x (P), \qquad \text{Vanilla landscape:~} x^{k + 1} = x^{k + 1 / 2},\qquad \mathcal{A} \assign
   {\ogd} . \]
The algorithm is called {\osgmhandsoffrx} (\Cref{alg:osgmhandsoffrx}).

\begin{algorithm}[H]
{\textbf{input:}  $x^1, P_1 \in \Pcal , \eta_k = \frac{c}{\sqrt{k}} \text{ or } \eta_k \equiv \frac{c}{\sqrt{K}}, c > 0$}\\
\For{$k = 1, 2, \dots$}{
      $x^{k + 1} = x^{k} - P_k \nabla f(x^k)$\\
      $P_{k+1} = \Pi_{\Pcal}[ P_k - \eta_k \nabla r_{x^k}(P_k) ]$
}
\caption{{\osgmhandsoffrx}\label{alg:osgmhandsoffrx}}
\end{algorithm}

The convergence analysis of {\osgmhandsoffrx} is similar to that for {\osgmhandsonrx}. However, the convergence guarantees of {\osgmhandsoffrx} are weaker due to the vanilla landscape action.

\begin{thm}[Global convergence]
  \label{thm:global-conv-handsoff-rx}Under \ref{A1} to \ref{A3}, for any benchmark stepsize $\hat{P} \in \Pcal$, 
  {\osgmhandsoffrx}  (\Cref{alg:osgmhandsoffrx}) with $\eta_k \equiv \tfrac{D}{2 L (L D
  + 1) \sqrt{K}}$ satisfies
\begin{equation} \label{eqn:global-conv-handsoff-rx}
  f (x^{K + 1}) - f^{\star} \leq [f (x^1) - f^{\star}] (\tfrac{1}{K} \textstyle \sum_{k = 1}^K r_{x^k} (\hat{P}) + \tfrac{3 L D (L D +1)}{\sqrt{K}} )^K.
\end{equation}
Moreover, {\osgmhandsoffrx} with $\eta_k = \tfrac{D}{2 L (L D + 1) \sqrt{k}}$ satisfies \eqref{eqn:global-conv-handsoff-rx} for all $K \geq 1$.
\end{thm}

\Cref{thm:global-conv-handsoff-rx} suggests a divergence behavior of {\osgmhandsoffrx} in the earlier iterations. Indeed, when the landscape takes no action to filter out bad stepsizes, the algorithm will remain unstable until the scheduler learns a good stepsize.

\begin{thm}[Local adaptivity]
\label{thm:local-trajectory-conv-rx-handsoff}Under the same assumptions as
\Cref{thm:global-conv-handsoff-rx}, for any benchmark sequence of stepsizes $\{\hat{P}_k\}$, {\osgmhandsoffrx} with
$\eta_k \equiv \tfrac{D}{2 L (L D + 1) \sqrt{K}}$ satisfies
\[ f (x^{K + 1}) - f^{\star} \leq [f (x^1) - f^{\star}] (
   \tfrac{1}{K} \textstyle \sum_{k = 1}^K r_{x^k} (\hat{P_k}) + \tfrac{3 L (L D + 1)
   (2 D + \mathsf{PL} (\{ \hat{P}_k \}))}{\sqrt{K}} )^K \text{~~ for any $\hat{P}_k \in \mathcal{P}$}.
   \]
\end{thm}

\begin{thm}[Superlinear convergence]
  \label{thm:local-superlinear-conv-handsoff-rx}Instate \ref{A1} to \ref{A4} and suppose $[\nabla^2 f(x^{\star})]^{-
  1} \in \mathcal{P}$, {\osgmhandsoffrx} with $\eta_k = \tfrac{D}{2 L (L D + 1) \sqrt{k}}$ satisfies 
\begin{equation*}
  f (x^{K + 1}) - f^{\star} \leq [f(x^1) - f^{\star}]( \tfrac{C_1}{K} +
     \tfrac{C_2}{\sqrt{K}} )^K,
\end{equation*}

where $C_1 = \tfrac{H^2 \kappa}{4 \mu^2} \big[ K_0 ( 1 - \tfrac{1}{\kappa} + 3 L D (L D + 1) )^{K_0} + 2 \kappa \big], K_0 = \lceil 36 \kappa^2 [L D (L D + 1)]^2 \rceil$ and $C_2 = 3 L D (L D + 1)$.
\end{thm}

\subsection{\osgmhandsoffhx} \label{sec:handsoffhx}
In this section, we assume $f$ is $L$-smooth, optionally $\mu$-strongly convex (\ref{A1}, optionally \ref{A2}), and instantiate {\osgm} with
\[ \ell_x (P) \assign h_x (P), \quad \text{Monotone landscape:~} f(x^{k+1}) \leq \min\{f(x^{k+1/2}), f(x^k)\},\quad \mathcal{A} \assign
   {\ogd} . \]
The algorithm is called {\osgmhandsoffhx} (\Cref{alg:osgmhandsoffhx}).

\begin{algorithm}[H]
{\textbf{input:}  $x^1, P_1 \in \Pcal , \eta_k = \frac{c}{\sqrt{k}} \text{ or } \eta_k \equiv \frac{c}{\sqrt{K}}, c > 0$}\\
\For{$k = 1, 2, \dots$}{
      $x^{k + 1/2} = x^{k} - P_k \nabla f(x^k)$\\
      Choose $x^{k+1}$ that satisfies $f(x^{k+1}) \leq \min\{f(x^{k+1/2}), f(x^k)\}$\\
      $P_{k+1} = \Pi_{\Pcal}[ P_k - \eta_k \nabla r_{x^k}(P_k) ]$
}
\caption{{\osgmhandsoffhx}\label{alg:osgmhandsoffhx}}
\end{algorithm}
The convergence analysis of {\osgmhandsoffhx} is similar to {\osgmhandsonhx}.
\begin{thm}[Global convergence]
  \label{thm:global-conv-handsoff-hx}Under \ref{A1} to \ref{A3}, for any benchmark stepsize $\hat{P} \in \Pcal$, 
  {\osgmhandsoffhx} (\Cref{alg:osgmhandsoffhx}) with $\eta_k \equiv \tfrac{D}{(L D + 1) \sqrt{K}}$
  satisfies
  \begin{align}
 f (x^{K + 1}) - f^{\star} \leq{} &  \min \big\{ \tfrac{\Delta^2}{K \max
    \{ \frac{1}{K} \sum_{k = 1}^K - h_{x^k} (\hat{P}) - \frac{3 D (L D
    + 1)}{\sqrt{K}}, 0 \}}, f (x^1) - f^{\star} \big\} \tag{convex}\\
    f (x^{K + 1}) - f^{\star} \leq{} & [f (x^1) - f^{\star}] ( 1 -
    \textstyle2 \mu \max \{ \tfrac{1}{K} \sum_{k = 1}^K - h_{x^k} (\hat{P}) -
    \tfrac{3 D (L D + 1)}{\sqrt{K}}, 0 \} ). \tag{$\mu$-strongly convex}
  \end{align}
Moreover, {\osgmhandsoffhx} with $\eta_k = \tfrac{D}{(L D +
  1) \sqrt{k}}$ satisfies the same bounds for all $K \geq 1$.
\end{thm}

\begin{thm}[Local adaptivity]
  \label{thm:local-trajectory-conv-handsoff-hx}Under the same assumptions as
  \Cref{thm:global-conv-handsoff-hx}, for any benchmark sequence of stepsizes $\{\hat{P}_k\}$, {\osgmhandsoffhx} with $\eta_k
  \equiv \tfrac{D}{(L D + 1) \sqrt{K}}$ satisfies
  \begin{align}
    f (x^{K + 1}) - f^{\star} \leq{} & \min \big\{ \tfrac{\Delta^2}{K \max
    \{ \frac{1}{K} \sum_{k = 1}^K - h_{x^k} (\hat{P}_k) - \rho_K(\{ \hat{P}_k \}), 0 \}}, f
    (x^1) - f^{\star} \big\}, \tag{convex}\\
     f (x^{K + 1}) - f^{\star} \leq{} & \textstyle [f (x^1) - f^{\star}] ( 1 -
    2 \mu \max \{ \tfrac{1}{K} \sum_{k = 1}^K - h_{x^k} (\hat{P}_k) - \rho_K(\{ \hat{P}_k \}), 0 \} ), \tag{$\mu$-strongly convex}
  \end{align}
  where $\rho_K(\{ \hat{P}_k \}) \assign \tfrac{3 D (L D + 1) (2 D + \mathsf{PL} (\{ \hat{P}_k \}))}{\sqrt{K}}$.
\end{thm}

\begin{thm}[Superlinear convergence]
  \label{thm:local-superlinear-conv-handsoff-hx}Instate the same assumptions
  as \Cref{thm:global-conv-handsoff-rx} and suppose $[\nabla^2 f(x^{\star})]^{-
  1} \in \mathcal{P}$ and $\eta_k = \tfrac{D}{(L D + 1) \sqrt{k}}$.
  {\osgmhandsoffhx} satisfies 
\begin{equation*}
  f (x^{K + 1}) - f (x^k) \leq [f (x^{1}) - f (x^k)]( \tfrac{C_1}{K} + \tfrac{C_2}{\sqrt{K}} )^K,
\end{equation*}
where $C_1 = \tfrac{H^2 \kappa^2}{4 \mu^2} [ K_0 ( 1 -
  \tfrac{1}{\kappa} + 3 D (L D + 1) )^{K_0} + 2 \kappa ], K_0 =
  \lceil 36 \kappa^2 [D (L D + 1)]^2 \rceil$ and $C_2 = 3 D(L D + 1)$.
\end{thm}

\subsection{\texttt{Monotone} {\osgmrx} and \texttt{Monotone Lookahead} {\osgmrx}}

{\osgmrx} does not require monotone landscape action to guarantee convergence, such as {\osgmhandsoffrx} and {\osgmhandsonrx}, 
but it does not hurt to equip both variants with monotone landscape action.
By \Cref{lem:feedback-progress}, the per iteration progress of monotone variants of {\osgmrx} are bounded by
\begin{align*}
   \texttt{Monotone}:{} ~& r_k \leq \min \{ r_{x^k} (P_k), 1 \} \leq r_{x^k} (P_k); \\
   \texttt{Monotone Lookahead}:{} ~& r_k \leq \min \{ r_{x^k} (P_k) - \tfrac{1}{4 L^2} \| \nabla r_{x^k} (P_k) \|_F^2, 1 \} \leq r_{x^k} (P_k) - \tfrac{1}{4 L^2} \| \nabla r_{x^k} (P_k) \|_F^2,
\end{align*}
and the bounds on the right-hand side can be further bounded in the same way for {\osgmhandsoffrx} and {\osgmhandsonrx}.
In other words, the convergence analysis of monotone variants reduces to that of non-monotone variants, {\osgmhandsoffrx} and {\osgmhandsonrx}, respectively.

\subsection{Proof of results in \Cref{sec:hands-off-rx}}

\subsubsection{Proof of Theorem \ref{thm:global-conv-handsoff-rx}}

\input{proofs/proof-global-conv-handsoff-rx.tex}

\subsubsection{Proof of Theorem \ref{thm:local-trajectory-conv-rx-handsoff}}

\input{proofs/proof-local-trajectory-conv-handsoff-rx.tex}

\subsubsection{Proof of Theorem \ref{thm:local-superlinear-conv-handsoff-rx}}

\input{proofs/proof-local-superlinear-conv-handsoff-rx.tex}

\subsection{Proof of results in \Cref{sec:handsoffhx}}

\subsubsection{Proof of Theorem \ref{thm:global-conv-handsoff-hx}}

\input{proofs/proof-global-conv-handsoff-hx.tex}

\subsubsection{Proof of Theorem \ref{thm:local-trajectory-conv-handsoff-hx}}

\input{proofs/proof-local-trajectory-conv-handsoff-hx.tex}

\subsubsection{Proof of Theorem \ref{thm:local-superlinear-conv-handsoff-hx}}

\input{proofs/proof-local-superlinear-conv-handsoff-hx.tex}

%% file: proofs/proof-global-conv-handsoff-rx.tex
We successively deduce that
\begin{align}
  f (x^{K + 1}) - f^{\star} \leq{} & ( \tfrac{1}{K} \textstyle \sum_{k = 1}^K r_k
  )^K  \tag{by \Cref{thm:reduction-ratio}}\\
  ={} & ( \tfrac{1}{K} \textstyle \sum_{k = 1}^K r_{x^k} (P_k) )^K \tag{by \Cref{lem:feedback-progress} + vanilla action}\\
  \leq{} & ( \tfrac{1}{K} \textstyle \sum_{k = 1}^K r_{x^k} (\hat{P}) + \tfrac{3 L D (L
  D + 1)}{\sqrt{K}} )^K \tag{by \Cref{lem:ogd-sqrtK-eta}}. 
\end{align}

For stepsize $\eta_k =\mathcal{O} ( \tfrac{1}{\sqrt{k}} )$, the regret guarantee applies to any $K \geq 1$ and provides anytime convergence.

%% file: proofs/proof-local-trajectory-conv-handsoff-rx.tex
Combining $f (x^{K + 1}) - f^{\star} \leq (\frac{1}{K} \sum_{k = 1}^K r_{x^k} (P_k) )^K$ with \Cref{lem:dynamic-sqrtK} completes the proof.

%% file: proofs/proof-local-superlinear-conv-handsoff-rx.tex
Plugging $\hat{P} = \tfrac{1}{L} I$ into \Cref{thm:global-conv-handsoff-rx}, we
have, for each $k = 1, \ldots, K$, that
\[ f (x^{k + 1}) - f^{\star} \leq [f (x^1) - f^{\star}] ( 1 -
   \tfrac{1}{\kappa} + \tfrac{3 L D (L D + 1)}{\sqrt{k}} )^k \]
and using strong convexity,
\begin{align}
  \| x^k - x^{\star} \|^2 \leq{} & \tfrac{2}{\mu} [f (x^{k + 1}) - f^{\star}] \leq \tfrac{2}{\mu} [f (x^1) - f^{\star}] ( 1 - \tfrac{1}{\kappa}
  + \tfrac{3 L D (L D + 1)}{\sqrt{k}} )^k . \nonumber
\end{align}
and we bound $\textstyle\sum_{k = 1}^K \| x^k - x^{\star} \|^2$ using
\begin{align}
  \textstyle\sum_{k = 1}^K \| x^k - x^{\star} \|^2 \leq{} & \tfrac{2}{\mu} [f (x^1) - f^{\star}] \textstyle\sum_{k = 1}^K ( 1 - \tfrac{1}{\kappa} + \tfrac{3 L D (L D
  + 1)}{\sqrt{k}} )^k \nonumber\\
  \leq{} & \tfrac{2}{\mu} [f (x^1) - f^{\star}] \textstyle\sum_{k = 1}^{\infty} (
  1 - \tfrac{1}{\kappa} + \tfrac{3 L D (L D + 1)}{\sqrt{k}} )^k
  \nonumber\\
  \eqqcolon{} & \tfrac{2}{\mu} [f (x^1) - f^{\star}] \textstyle\sum_{k = 1}^{\infty}
  e_k \nonumber
\end{align}
Let $K_0 = \lceil 36 \kappa^2 [L D (L D + 1)]^2 \rceil$. We have $e_k \leq
( 1 - \tfrac{1}{2 \kappa} )^k$ for all $k \geq K_0$ and
\begin{align}
  \textstyle\sum_{k = 1}^{\infty} e_k ={} & \textstyle\sum_{k = 1}^{K_0} e_k + \textstyle\sum_{k = K_0 +
  1}^{\infty} e_k \leq  K_0 ( 1 - \tfrac{1}{\kappa} + 3 L D (L D + 1) )^{K_0} + 2
  \kappa \nonumber
\end{align}
Using \Cref{thm:global-conv-handsoff-rx} and \Cref{lem:local-opthessian-rx}, we
deduce that
\begin{align}
  f (x^{K + 1}) - f^{\star} \leq{} & [f (x^1) - f^{\star}] (
  \tfrac{1}{K} \textstyle\sum_{k = 1}^K r_{x^k} (\hat{P}) + \tfrac{3 L D (L D +
  1)}{\sqrt{K}} )^K \nonumber\\
  \leq{} & [f (x^1) - f^{\star}] ( \tfrac{H^2 \kappa}{4 \mu^2}
  \tfrac{1}{K} \textstyle\sum_{k = 1}^K \| x^k - x^{\star} \|^2 + \tfrac{3 L D (L D +
  1)}{\sqrt{K}} )^K \nonumber\\
  \leq{} & [f (x^1) - f^{\star}] ( \tfrac{H^2 \kappa}{4 \mu^2 K} [
  K_0 ( 1 - \tfrac{1}{\kappa} + 3 L D (L D + 1) )^{K_0} + 2 \kappa
  ] + \tfrac{3 L D (L D + 1)}{\sqrt{K}} )^K \nonumber\\
  ={} & [f (x^1) - f^{\star}] ( \tfrac{C_1}{K} + \tfrac{C_2}{\sqrt{K}}
  )^K \nonumber
\end{align}

and this completes the proof.

%% file: proofs/proof-global-conv-handsoff-hx.tex
Similar to \Cref{thm:global-conv-handsoff-rx}, chaining \Cref{thm:reduction-hypergrad}, \Cref{lem:feedback-progress} and \Cref{lem:ogd-sqrtK-eta} completes the proof.

%% file: proofs/proof-local-trajectory-conv-handsoff-hx.tex
Similar to \Cref{thm:local-trajectory-conv-rx-handsoff}, chaining \Cref{thm:reduction-hypergrad}, \Cref{lem:feedback-progress} and \Cref{lem:dynamic-sqrtK} completes the proof.

%% file: proofs/proof-local-superlinear-conv-handsoff-hx.tex
Plugging $\hat{P} = \tfrac{1}{L} I$ into \Cref{thm:global-conv-handsoff-hx}, we
have, for each $k = 1, \ldots, K$, that
\[ f (x^{k + 1}) - f^{\star} \leq [f (x^1) - f^{\star}] ( 1 -
   \tfrac{1}{\kappa} + \tfrac{3 D (L D + 1)}{\sqrt{k}} )^k. \]

Using the same argument as \Cref{thm:local-superlinear-conv-handsoff-rx}, we have
\begin{align}
  \| x^k - x^{\star} \|^2 \leq{} & \tfrac{2}{\mu} [f (x^{k + 1}) - f^{\star}]  \leq{} \tfrac{2}{\mu} [f (x^1) - f^{\star}] ( 1 - \tfrac{1}{\kappa}
  + \tfrac{3 D (L D + 1)}{\sqrt{k}} )^k . \nonumber
\end{align}
and we bound $\textstyle\sum_{k = 1}^K \| x^k - x^{\star} \|^2$ using
\begin{align}
  \textstyle\sum_{k = 1}^K \| x^k - x^{\star} \|^2 \leq{} & \tfrac{2}{\mu} [f (x^1) - f^{\star}] \textstyle\sum_{k = 1}^K ( 1 - \tfrac{1}{\kappa} + \tfrac{3 D (L D
  + 1)}{\sqrt{k}} )^k \nonumber\\
  \leq{} & \tfrac{2}{\mu} [f (x^1) - f^{\star}] \textstyle\sum_{k = 1}^{\infty} (
  1 - \tfrac{1}{\kappa} + \tfrac{3 D (L D + 1)}{\sqrt{k}} )^k
  \nonumber\\
  \eqqcolon{} & \tfrac{2}{\mu} [f (x^1) - f^{\star}] \textstyle\sum_{k = 1}^{\infty}
  e_k \nonumber
\end{align}
Let $K_0 = \lceil 36 \kappa^2 [D (L D + 1)]^2 \rceil$. We have $e_k \leq
( 1 - \tfrac{1}{2 \kappa} )^k$ for all $k \geq K_0$ and
\begin{align}
  \textstyle\sum_{k = 1}^{\infty} e_k ={} & \textstyle\sum_{k = 1}^{K_0} e_k + \textstyle\sum_{k = K_0 +
  1}^{\infty} e_k \leq K_0 ( 1 - \tfrac{1}{\kappa} + 3 D (L D + 1) )^{K_0} + 2
  \kappa \nonumber
\end{align}
Similar to the analysis of \Cref{thm:local-superlinear-conv-handson-hx},  using $\sum_{k = 1}^K h_k \leq \sum_{k = 1}^K h_{x^k} ([\nabla^2 f
(x^{\star})]^{- 1}) + \tfrac{3 D (L D + 1)}{\sqrt{K}}$, we further deduce
that
\begin{align}
  f (x^{K + 1}) - f (x^k) \leq{} & [f (x^1) - f^{\star}] ( \tfrac{2
  L}{K} \textstyle \sum_{k = 1}^K [ h_k + \tfrac{f (x^k) - f^{\star}}{\| \nabla
  f (x^k) \|^2} ] )^K \nonumber
\end{align}
and that
\begin{align}
 \textstyle \sum_{k = 1}^K [ h_k + \tfrac{f (x^k) - f^{\star}}{\| \nabla f
  (x^k) \|^2} ] \leq{} & \textstyle \sum_{k = 1}^K [ h_{x^k} ([\nabla^2 f
  (x^{\star})]^{- 1}) + \tfrac{f (x^k) - f^{\star}}{\| \nabla f (x^k)
  \|^2} ] + 3 D (L D + 1) \sqrt{K} \nonumber\\
  \leq{} & \textstyle \tfrac{H^2 \kappa}{8 \mu^3} \sum_{k = 1}^K \| x^k - x^{\star} \|^2 +
  3 D (L D + 1) \sqrt{K} . \nonumber
\end{align}
Plugging the bounds back completes the proof.